\documentclass{article}
\usepackage{amssymb}
\usepackage{graphicx}
\usepackage{indentfirst}
\usepackage{lmodern}
\usepackage[english]{babel}
\usepackage{amsmath}
\numberwithin{equation}{section}  
\usepackage{amsthm}
\usepackage{mathtools}
\usepackage{cancel}
\newcommand{\diam}{\operatorname{diam}}
\usepackage{tikz}
\usepackage{pgfplots}
\pgfplotsset{compat=1.17}
\usetikzlibrary{calc}
\usetikzlibrary{decorations.pathreplacing}
\usepackage{bm}
\usepackage{mathrsfs}
\usepackage{esint}
\usepackage{physics}  
\usepackage{makecell} 
\usepackage{rotating}
\usepackage{fancyhdr}
\usepackage{enumerate}
\usepackage{enumitem}
\usepackage{comment}
\usepackage{textcomp}
\usepackage{colortbl}
\usepackage[dvipsnames]{xcolor}
\usepackage{float}
\usepackage{geometry}
\usepackage{multicol}
\usepackage{hyperref}
\usepackage{cleveref}  
\usepackage{titlesec}

\newcommand{\appendixtitleformat}{
  \titleformat{\section}
    {\normalfont\LARGE\bfseries}
    {Appendix \thesection}
    {1em}{}

  \titleformat{\subsection}
    {\normalfont\Large\bfseries}
    {\thesection.\arabic{subsection}}
    {1em}{}
}
\allowdisplaybreaks
\hypersetup{
    colorlinks=false, 
    pdfborder={0 0 0} 
}
\DeclareSymbolFont{Symbols}{OMS}{cmsy}{m}{n}
\SetSymbolFont{Symbols}{bold}{OMS}{cmsy}{b}{n}
\DeclareMathSymbol{\Setminus}{\mathbin}{Symbols}{"6E}
\renewcommand{\emptyset}{\font\cmsy = cmsy10 at 10pt
 \hbox{\cmsy \char 59}}
\newcommand{\tcr}{\textcolor{red}}

\makeatletter

\newcommand{\@myaddress}{}
\newcommand{\@myemail}{}

\newcommand{\myaddress}[1]{\gdef\@myaddress{#1}}
\newcommand{\myemail}[1]{\gdef\@myemail{#1}}

\newcommand{\printaddress}{%
  \par\bigskip
  \noindent\@myaddress\par
  \noindent\emph{Email address:} \texttt{\@myemail}
}

\makeatother

\newcommand{\ud}[0]{\,\mathrm{d}}

\newcommand{\dist}[0]{\operatorname{dist}}

\newcommand{\Norm}[2]{\|#1\|_{#2}}

\newcommand{\BNorm}[2]{\Big\|#1\Big\|_{#2}}
\newcommand{\pair}[2]{\langle #1,#2 \rangle}

\newcommand{\ave}[1]{\langle #1\rangle}

\newcommand{\supp}[0]{\operatorname{supp}}
\newcommand{\loc}[0]{\operatorname{loc}}

\newcommand{\R}{\mathbb{R}}

\newcommand{\Z}{\mathbb{Z}}

\renewcommand{\S}{\mathbb{S}}

\newcommand{\prob}[0]{\mathbb{P}}

\newcommand{\D}[0]{\mathbb{D}}
\newcommand{\E}[0]{\mathbb{E}}
\renewcommand{\P}[0]{\mathbb{P}}
\newcommand{\eps}[0]{\varepsilon}

\newcommand{\children}[0]{\operatorname{ch}}

\swapnumbers \numberwithin{equation}{section}

\theoremstyle{plain}
\newtheorem{theorem}[equation]{Theorem}
\newtheorem{proposition}[equation]{Proposition}
\newtheorem{corollary}[equation]{Corollary}
\newtheorem{lemma}[equation]{Lemma}

\newtheorem*{FeaturesSlab}{Features of the infinite slab}
\newtheorem*{FeaturesCantor}{Features of the Cantor--type set}
\newtheorem*{dyadicformula}{Dyadic Formula}

\theoremstyle{definition}
\newtheorem{definition}[equation]{Definition}
\newtheorem{notation}[equation]{Notation}

\newtheorem{assumption}[equation]{Assumption}
\newtheorem{remark}[equation]{Remark}

\title{Dyadic Representation of Fractional Integrals on Spaces of Homogeneous Type and Related Boundedness Phenomena}
\author{Lorenzo Zacchini}

\myaddress{Aalto University, Department of Mathematics and Systems Analysis,
P.O. Box 11100, FI-00076 Aalto, Finland}

\myemail{lorenzo.zacchini@aalto.fi}
\date{}

\begin{document}

\maketitle

\begin{abstract}
We establish a dyadic representation formula for fractional integral operators in the setting of spaces of homogeneous type. The formula expresses such operators in terms of simpler dyadic model operators, such as so-called paraproducts and shifts, and provides a unified framework that applies to different notions of fractional integrals arising in this general context.

Motivated by the derivation of this representation, we also investigate the structure of fractional integral operators in spaces of homogeneous type. In this setting, two natural definitions—referred to as volumic and metric—are in general non-equivalent and exhibit markedly different boundedness properties. While the former enjoys well-known $L^p$--$L^q$ estimates, the latter may fail to satisfy any such bounds without additional assumptions; nevertheless, weaker forms of boundedness can still be recovered under suitable conditions.

The representation formula is intended as a tool for the study of Schatten-type properties of commutators $[b,T]$ of fractional integrals and pointwise multipliers, which will be addressed in a forthcoming companion work by T. Hytönen and the author. Beyond this application, the representation itself may be of independent interest and useful in other contexts.
\end{abstract}

 \tableofcontents 

\section{Introduction}

Dyadic representation formulas have become an important tool in modern harmonic analysis, with numerous developments across a variety of settings and operators. Earlier foundational contributions include \cite{Figiel:90,LMPT,NTV:Tb,Pet:00}, while the modern dyadic representation theorem was first established by T. Hytönen \cite{Hyt:A2} for Calderón--Zygmund operators on \(\R^d\). This initial formulation was subsequently streamlined in \cite{Hyt:17}. Other notable works include the representation formula for bi-parameter singular integrals in \cite{Mar:12}, and that for bilinear Calderón–Zygmund operators in \cite{LMOV}. 

In line with the present work, it is important to highlight the approach introduced in \cite{GH:18}, later extended to doubling metric spaces in \cite{Hyt:Enflo}. An additional refinement, developed in the setting of spaces of homogeneous type, can be found in \cite[Theorem 16.1]{Hyt:26}. Following the example of this latter development, the present paper aims to establish an analogous dyadic representation formula for fractional integrals. The original goal of deriving such a formula also motivated the study of these operators in their own right within this paper, whose main ideas and contents we now briefly describe.


The setting we work in is that of a space of homogeneous type \((X, \rho, \mu)\), that is, a set \(X\) equipped with a quasi-distance \(\rho\) and a Borel measure \(\mu\) satisfying the doubling property, meaning that there exists a constant \(c_\mu\geq1\) such that
\begin{equation*}
    0 < \mu\bigl(B(x,2r)\bigr) \leq c_\mu \, \mu\bigl(B(x,r)\bigr) < \infty,
\end{equation*}
for all balls \(B(x,r) := \{ y \in X \mid \rho(x,y) < r \}\). In these spaces, one can construct so-called dyadic cubes (see \cite{Christ:90,David:88,David:91,HK:12}), forming what is known as a \emph{dyadic system}, in a way that is analogous to the standard dyadic cubes in \(\R^d\). Roughly speaking, such a system consists of a collection of measurable sets, organized into levels indexed by \(k \in \mathbb{Z}\), which partition the space at each level, are nested from one level to the next, and resemble balls of radius \(\delta^k\) for some \(\delta \in (0,1)\). These cubes provide a flexible substitute for the classical Euclidean cubes and are a fundamental tool in harmonic analysis on spaces of homogeneous type. Moreover, it is possible to construct these cubes in a random way, that is, to build a probability space \((\Omega,\P)\) such that each element of \(\Omega\) corresponds to a distinct random dyadic system (see \cite{AH:13,Hyt:Enflo,HK:12,Hyt:26,HM:12,NRV,NTV:Tb}).

In this framework, we study the class of fractional integral operators. Such operators hold a central place in analysis thanks to their broad range of applications in potential theory, harmonic analysis, partial differential equations, and Sobolev embeddings. 

The prototype of fractional operators on the Euclidean space $\mathbb{R}^d$ is the Riesz potential, which is extensively treated, for instance, in \cite[Chapter 1]{Grafakos:MFA}. In the present work, we consider fractional integral operators in the more general setting of a space of homogeneous type. More precisely, we consider operators of the form
\[
    Tf(x)=\int_X K(x,y)f(y)\ud\mu(y),
\]
in terms of a kernel $K\in L^1_{\text{loc}}(X\times X)$ satisfying a uniform size estimate of the form
\[
\abs{K(x,y)} \leq C_K\frac{\phi(x,y)}{\mu\bigl(B(x,\rho(x,y))\bigr)}=:C_K\frac{\phi(x,y)}{V(x,y)},
\]
where 
\begin{equation*}
\phi(x,y) :=
\begin{cases}
V(x,y)^{\varepsilon}, & \text{for } \varepsilon \in (0,1) \quad \text{(volumic type)}, \\[4pt]
\rho(x,y)^{\eta},     & \text{for } \eta \in (0,\infty) \quad \text{(metric type)}.
\end{cases}
\end{equation*}
and $C_K\in(0,\infty)$ is a positive constant. For certain results, including the dyadic representation formula, we will also require an additional $\phi$-dependent smoothness estimate.

Given the very general nature of spaces of homogeneous type, the two possible choices of $\phi$ lead to genuinely distinct and non-comparable definitions (see Proposition~\ref{prop:Bessel} and Appendices~\ref{app:Slab} and \ref{app:Cantor}). On the other hand, given that both types of fractional integrals appear in the literature, there is no compelling reason to privilege one definition over the other. 

Volumic fractional integrals have been studied in the context of spaces of homogeneous type, for instance in \cite{GV:90}. In that work, however, the additional assumption that the space be \emph{normal} is imposed, namely that $V(x,y)\approx \rho(x,y)$, in which case the two definitions essentially coincide. These operators are also investigated without this additional assumption, for example in \cite{BC:95, Pan:92}. 

Metric fractional integrals frequently arise in applications, as for instance in \cite[Theorem 3.22 and Section 9.1]{Hei:01}, \cite[Page 295]{HMP:25} and \cite[Eq.~(7)]{PW:01}. As noted in the latter paper, metric-type operators also arise in various examples of potential operators associated with the regularity theory of subelliptic differential equations, such as those arising in \cite{FGW:94,FL:83,Hormander:67}.
 Moreover, a noteworthy example is given by the fractional powers $\Delta_\lambda^{-\alpha/2}$ of the Bessel Laplacian $\Delta_\lambda$, which turn out to be of metric type (see \cite{BDLLW:25} and \cite[Section 12]{HW:25}, as well as Subsection~\ref{subsec:Bessel} of the present paper). 

In the study of the dyadic representation formula, it is in some parts relevant to have a certain degree of $L^p$ control for this class of operators. In the volumic case, boundedness results from $L^p$ to $L^q$, with $p,q\in(1,\infty)$ suitably related, are already well known \cite{BC:95, GV:90}. The same cannot be said for operators of metric type. Without imposing any additional assumptions on the underlying space of homogeneous type, there exist metric fractional integral operators for which some functions in \(L^p\) are mapped to the function identically equal to $+\infty$ (see Appendix~\ref{app:Cantor}). Consequently, no general \(L^p\)-boundedness theory can hold for metric fractional integrals in the full generality of spaces of homogeneous type.

Even in spaces satisfying natural geometric assumptions, such as positive lower and upper dimension, there are metric fractional integral operators that are unbounded from \(L^p\) to \(L^q\) for every \(p,q\in(1,\infty)\). This follows from the example constructed in Appendix~\ref{app:Slab}. However, under the assumption of positive lower dimension, one can still obtain a weaker form of boundedness from $L^p$ into a suitable sum space (see Theorem~\ref{thm:bddfracmet}), which will be sufficient for the purposes of the present work.

With a clear picture of these operators and their differences in place, it becomes natural and feasible to aim for a unified form of dyadic representation formula. In order to give the reader a general overview, we now provide an intuitive sketch of this formula, whose full and rigorous statement can be found in Theorem~\ref{thm:DRT}.

\begin{dyadicformula}
    Let $T$ be a $\phi$-fractional integral operator defined on a space
of homogeneous type $(X,\rho,\mu)$ and satisfying the additional
regularity and summability assumptions in \eqref{eq:CZ1phi} and
\eqref{eq:omegasum}, respectively. Then $T$ admits the following
representation formula:
    \[
T=T_1+\mathbb{E}\Big(\Pi_{T1}+\Pi^*_{T^*1}+\sum\S^\phi\Big),
\]
where:
\begin{itemize}
    \item \(T_1\) is a rank-one operator, appearing only when $\mu(X)<\infty$.
    \item $\Pi_{T1}$ and $\Pi_{T^*1}$ denote the dyadic paraproducts, see Definition~\ref{def:paraT1}, which satisfy a $\phi$-dependent normalization that is useful in the intended applications. 
    \item $\S^\phi$ denote the dyadic shifts, with $\phi$-dependent normalization, see Definition~\ref{def:shifts}. 
    \item The expectation \(\mathbb E\) is taken with respect to the underlying probability space of random dyadic systems.
\end{itemize}
\end{dyadicformula}

It is worth observing that all the operators appearing in the formula admit pointwise domination by fractional integrals (see Propositions~\ref{prop:DomPhiNormParaprod} and \ref{prop:shiftsdom} and Lemma~\ref{lem:Rank1Est}), and hence inherit their \(L^p\)-boundedness properties. In this sense, the individual components of the representation are ``no worse'' than the original operator itself, while they are in other ways ``better'', having an explicit dyadic structure.

In the Euclidean setting, a dyadic representation formula for the operator

$$
I_\lambda f(x):=\int_{\mathbb{R}^n}\frac{f(y)}{|x-y|^\lambda}\,\mathrm{d}y,
\qquad 0<\lambda<n,
$$
is established in \cite[Theorem 1.2]{PS:23}. The representation formula established in the present paper, on the other hand, is developed in the more general setting of spaces of homogeneous type. Moreover, it applies to arbitrary fractional integrals, rather than being restricted to operators whose kernels admit such a specific explicit form. This generality constitutes the main novelty of our result.

In a broader sense, a dyadic \emph{domination} (i.e., a pointwise upper bound, as opposed to an identity as in our representation) for fractional integrals is well known; see \cite[Section 4]{LMPT}. While such domination is sufficient---and often even preferable---for many purposes, a linear representation of the form obtained here turns out to be more suitable for our intended applications. 

More specifically, the primary motivation for establishing this representation arises from its application to Schatten--Lorentz estimates for commutators \([b,T]\) of pointwise multipliers and fractional integrals. These applications will be developed in a forthcoming companion paper with T.~Hytönen, complementing the results of Hytönen and Wu in~\cite{HW:25}. Although this is the main motivation for the present work, the representation formula is of independent interest and may have applications in other contexts; we therefore develop it separately here.

The remainder of the paper is organized as follows.
In Section~\ref{sec:prel} we collect the necessary background on spaces of homogeneous type,
Muckenhoupt weights, and dyadic structures.
In Section~\ref{sec:fracint} we introduce the class of fractional integral operators under consideration
and discuss the distinction between the volumic and metric cases, together with their boundedness properties. Section~\ref{sec:DyadModOp} introduces and explores various properties of the dyadic model operators, which form the building blocks for the representation formula.
Section~\ref{sec:DRT} is devoted to the statement and proof of the dyadic representation formula.
Finally, in Appendices~\ref{app:Slab} and~\ref{app:Cantor} we present concrete examples illustrating the differences between the two notions of fractional integrals
and the possible failure of $L^p$ boundedness in the metric setting.

\section{Preliminaries}\label{sec:prel}
Throughout this work, in order to streamline the formulas and focus on the essential terms, we will make frequent use of the symbols $\gtrsim$, $\lesssim$, and $\approx$. Unless otherwise specified, all such relations are understood up to positive multiplicative constants whose precise values (and dependence) are irrelevant for the estimates under consideration. Accordingly, $\gtrsim$ and $\lesssim$ denote one-sided inequalities up to such constants, while $\approx$ denotes two-sided comparability.
 For example, all constants related to the underlying space (such as the constant in the quasi-triangle inequality or the doubling constant), as well as the uniform constants that arise in the definitions of fractional integrals, will be implicitly absorbed in these symbols.

 In addition to the notation fixed above, and for the sake of clarity, we adopt the standard convention that for $p\in(1,\infty)$ the symbol $p'$ denotes its Hölder conjugate, that is, $p':=\frac{p}{p-1}$.

With these notational conventions in place, we now recall the basic notions and results that will be needed throughout the paper.

\subsection{Spaces of homogeneous type and Muckenhoupt weights}\label{subsec:HomMuck}
We begin by recalling the definition of a space of homogeneous type (sometimes referred to as \textit{homogeneous space}) from \cite[Section 2]{BC:95} and \cite[Section 1]{Christ:90}.

\begin{definition}
    We call a \textit{space of homogeneous type} a triple $(X, \rho, \mu)$, where:
    \begin{itemize}
        \item $X$ is a set.
        \item $\rho$ is a quasi-distance (also called a quasi-metric), meaning a function from $X \times X$ to $[0, \infty)$ that satisfies all the properties of a distance, except for the triangle inequality. The latter is replaced by a relaxed form: there exists a constant $A_0 \in [1, \infty)$ such that
\[
\rho(x, y) \leq A_0 \big( \rho(x, z) + \rho(z, y) \big), \quad \text{for all } x, y, z \in X.
\]
In addition, we require that the balls $B(x, r) := \{ y \in X \mid \rho(x, y) < r \}$ are open in the topology induced by $\rho$.
        \item $\mu$ is a positive Borel measure defined on a $\sigma$-algebra of subsets of $X$ which contains the $\rho$-balls, for which there exists a constant $c_\mu \geq 1$ such that the following \textit{doubling condition} holds:
        \begin{equation*}
    0 < \mu\bigl(B(x,2r)\bigr) \leq c_\mu \, \mu\bigl(B(x,r)\bigr) < \infty,
\end{equation*}
for all \(x \in X\) and \(r > 0\). 
    \end{itemize}
\end{definition}

It is useful to observe that, given a quasi-distance $\rho$, in general the balls are not open in the topology induced by $\rho$. However, by \cite[Theorem 2]{MS:79}, it is always possible to construct a quasi-distance $\rho'$ that is equivalent to $\rho$ and such that the balls defined with respect to $\rho'$ are open.

Since we will frequently refer to the measure of balls, for $x\in X$ we denote
\[
V(x,r) := \begin{cases}
    \mu\bigl(B(x,r)\bigr), \qquad & r>0;\\
    \mu(\{x\}),\qquad  & r=0.
\end{cases}
\]
With a slight abuse of notation, for $x,y\in X$, we also write
\[
V(x,y) := V\bigl(x,\rho(x,y)\bigr)
= \begin{cases}
    \mu\bigl(B(x,\rho(x,y))\bigr), \qquad & x\neq y ;\\
    \mu(\{x\}),\qquad  & x=y.
\end{cases}
\]

A further notation that will be frequently used is the average of a measurable function $f$ over a measurable set $E$ of finite measure, defined as
\[
\ave{f}_E := \fint_E f\,\mathrm{d}\mu := \frac{1}{\mu(E)} \int_E f\,\mathrm{d}\mu.
\]
In the case of the whole space $X$, we define
\[
\ave{f}_X := 
\begin{cases}
0, & \text{if } \mu(X) = \infty, \\
\fint_X f\,\mathrm{d}\mu, & \text{if } \mu(X) < \infty.
\end{cases}
\]

We continue this brief introductory subsection by recalling the notions that quantify the growth of the measure of balls, namely lower and upper dimension, as well as Ahlfors regularity.

\begin{definition}\label{def:DimSpace}
    We say that a quasi-metric measure space $(X,\rho,\mu)$
\begin{itemize}
  \item  has {\em lower dimension} $d>0$ if
\begin{equation}\label{eq:lowerD}
  \Big(\frac{R}{r}\Big)^d\lesssim\frac{V(x,R)}{V(x,r)};
\end{equation}
 \item has {\em upper dimension} $D>0$ if 
\begin{equation}\label{eq:upperD}
  \frac{V(x,R)}{V(x,r)}\lesssim\Big(\frac{R}{r}\Big)^D;
\end{equation}
\item is 
{\em Ahlfors $\Delta$-regular} if there exists $\Delta>0$ such that
it has both lower and upper dimension $\Delta$, that is,
\begin{equation}\label{eq:AhlforsReg}
 V(x,r) \approx r^{\Delta};
\end{equation}
\end{itemize}
in each case uniformly in $x\in X$ and $0<r\leq R\lesssim \operatorname{diam}(X)$.
\end{definition}

\begin{remark}\label{rmk:UpperDim}
    Unlike the positive lower dimension, which a priori may fail to hold (see Appendix~\ref{app:Cantor}), it is well known (see, e.g., \cite[Chapitre III]{CW:71}) that a space of homogeneous type $(X,\rho,\mu)$ always has an upper dimension $D\in(0,\infty)$. Indeed, given $x\in X$ and $0<r\leq R\lesssim \operatorname{diam}(X)$, there exists a unique $k\in\Z$ such that $2^{k-1}\leq\frac{R}{r}<2^k$. Thus, by the doubling property of $\mu$,
    \[
    V(x,R)=V\Big(x,\frac{R}{r}r\Big)\leq V(x,2^kr)\leq (c_\mu)^kV(x,r),
    \]
    and
    \[
    (c_\mu)^k=(2^k)^{\log_2(c_\mu)}\leq \Big(2\frac{R}{r}\Big)^{\log_2(c_\mu)}=c_\mu\Big(\frac{R}{r}\Big)^{\log_2(c_\mu)}.
    \]
    From this, one may therefore deduce that $(X,\rho,\mu)$ has upper dimension $0<D\leq \log_2(c_\mu)$.
\end{remark}

A space of homogeneous type with lower dimension $d\in(0,\infty)$ has no atoms.

\begin{lemma}\label{lem:noatoms}
    Let $(X,\rho,\mu)$ be a space of homogeneous type with lower dimension $d\in(0,\infty)$. Then \[\mu\bigl(\{x_0\}\bigr)=0,\quad \text{ for all $x_0\in X$.}\] 
\end{lemma}
\begin{proof}
    Let $x_0\in X$. By the lower dimension property, in particular for $k\geq1$ we have that 
    \begin{equation*}
        \frac{\mu\bigl(B(x_0,1)\bigr)}{\mu\bigl(B(x_0,2^{-k})\bigr)}\gtrsim \left(\frac{1}{2^{-k}}\right)^d,
    \end{equation*}
    that is,
    \begin{equation*}
        \mu\bigl(B(x_0,2^{-k})\bigr)\lesssim 2^{-kd}\mu\bigl(B(x_0,1)\bigr)<\infty.
    \end{equation*}
    Hence, by the monotone continuity of the measure $\mu$, we have that
    \begin{equation*}
        0\leq\mu\bigl(\{x_0\}\bigr)=\mu\Big(\bigcap_{k\geq1} B(x_0,2^{-k})\Big)=\lim_{k\to\infty} \mu\big(B(x_0,2^{-k})\bigr)\lesssim \mu\bigl(B(x_0,1)\bigr)\lim_{k\to\infty} 2^{-kd}=0,
    \end{equation*}
    which proves the claim.
\end{proof}

Finally, we recall the definition of Muckenhoupt weights.

\begin{definition}[Muckenhoupt weights]\label{def:Muckweight}
Let $(X,\rho,\mu)$ be a quasi-metric measure space. A weight $w : X \to [0,\infty)$ is said to belong to the Muckenhoupt class $A_p$, $1<p<\infty$, if
\[
  [w]_{A_p} := \sup_{B \text{ ball}} 
  \left( \fint_B w \, \ud\mu \right) 
  \left( \fint_B w^{-\frac{1}{p-1}} \, \ud\mu \right)^{p-1} < \infty.
\]
\end{definition}

An equivalent characterization of Muckenhoupt weights can be given in terms of the Hardy--Littlewood maximal operator, as follows.

\begin{theorem}\label{thm:MuckH-L}
    Let $w$ be a locally integrable, non-negative function on a space of homogeneous type $(X,\rho,\mu)$. 
For $1 < p < \infty$, the weight $w$ belongs to the class $A_p$ if and only if there exists a constant $C > 0$ such that for every locally integrable function $f$,
\[
\int_X \bigl(M f(x)\bigr)^p w(x) \, \ud\mu(x) \leq C \int_X |f(x)|^p w(x) \, \ud\mu(x),
\]
where $M$ denotes the Hardy--Littlewood maximal operator defined by
\[
Mf(x) := \sup_{B \ni x} \fint_B |f| \, \ud\mu,
\]
and the supremum is taken over all balls $B$ containing $x$.
\end{theorem}

This result is proved, in the case of $\R^n$ equipped with the Lebesgue measure, in \cite[Chapter 5, Section 2]{Stein:93}. For the general case of a (quasi)-metric measure space equipped with a doubling measure, I refer the reader to \cite[Theorem 9, p.5]{ST:89}, or \cite[Lemma 7 and Theorem 9]{KLVY}. If the reader is interested, moreover, in \cite[Theorem 2.3.1]{GGKK:98} one can also find a slightly more general result, in which Muckenhoupt weights are characterized in terms of the boundedness of the maximal operator on weighted Lorentz spaces.

\subsection{Dyadic and martingale structures}\label{subsec:overcubes}
A fundamental tool, which we will make extensive use of throughout this work, is the concept of dyadic cubes in the context of spaces of homogeneous type. 

A significant part of this subsection is based on \cite[Sections 6 and 17]{Hyt:26}, from which we recall only the material relevant to the present work. In addition, we provide proofs of certain further properties, presumably known, for the sake of completeness and later use.

\subsubsection*{Dyadic cubes}
We begin by clarifying what we mean by a dyadic system. Essentially, this is a system of grids with useful covering and nestedness properties, whose elements—called dyadic cubes, in analogy with the terminology used in $\R^d$—can be treated, at each scale, as ball-like sets with essentially fixed radius.

\begin{definition}\label{def:cubes}
A {\em dyadic system} $\mathscr D$, also called a {\em system of dyadic cubes}, on a space of homogeneous type $(X,\rho,\mu)$ is a family
\begin{equation*}
  \mathscr D=\bigcup_{k\in\Z}\mathscr D_k,
\end{equation*}
satisfying the following conditions:
\begin{enumerate}
  \item\label{it:Dpartition} For each integer $k$, the collection $\mathscr D_k$ forms a measurable partition of $X$;
  \item\label{it:Drefine} Each $\mathscr D_{k+1}$ is a refinement of $\mathscr D_k$;
  \item\label{it:DlikeB} There exist constants $\delta\in(0,1)$ and $0<c_0\leq C_0<\infty$ such that every cube $Q\in\mathscr D_k$ resembles a ball of radius $\delta^k$ in the sense that, for some point $z_Q\in X$,
\begin{equation}\label{eq:cubeVsBalls}
  B(z_Q,c_0\delta^k)\subseteq Q\subseteq B(z_Q,C_0\delta^k).
\end{equation}
\end{enumerate}
The dyadic system is called {\em connected} if it additionally satisfies:
\begin{enumerate}\setcounter{enumi}{3}
  \item\label{it:Dconnected} For every pair of cubes $P,Q\in\mathscr D$, there exists a cube $R\in\mathscr D$ such that $P\cup Q\subseteq R$.
\end{enumerate}
\end{definition}

A construction of such a system for spaces of homogeneous type, building on earlier ideas of David \cite{David:88,David:91}, was originally given in full generality by Christ \cite{Christ:90}, although the partition and refinement properties are ensured only up to sets of measure zero. A version of the construction that exactly meets the conditions and allows the system to be connected can be found in \cite{HK:12}.

These dyadic cubes retain most of the fundamental characteristics of the classical dyadic cubes in $\R^d$, with one important distinction concerning property~\ref{it:Drefine}. In $\R^d$, the refinement is strict: every cube at a finer level is strictly contained in some cube at the coarser level. However, this need not be true in general for spaces of homogeneous type. That is, the same set $Q$ may belong to multiple levels $\mathscr D_k$ of the dyadic system. For this reason, we will sometimes refer to a cube $Q \in \mathscr D$ simply as a set, while other times we will regard it as a pair $(Q, k)$, in order to explicitly include the information about the level $k$.

We now establish a property of dyadic cubes that will be useful later. 
Given a point $x\in X$, the family of all dyadic cubes containing $x$ 
does not necessarily exhaust the entire space $X$ as their scale tends 
to infinity; in principle, it may expand only within a certain 
``quadrant''. Nevertheless, the measures of these cubes always converge 
to the measure of the whole space $X$. In particular, when $\mu(X)=\infty$ 
their measures diverge to infinity, while if $\mu(X)<\infty$ the cubes 
eventually coincide with $X$ itself.

\begin{lemma}\label{lem:CubesCover}
    Let $\mathscr D=\bigcup_{k\in\Z}\mathscr D_k$ be a dyadic system on a space of homogeneous type $(X,\rho,\mu)$. Moreover, let $x\in X$ and, for each $k\in \Z$, let us denote by $Q_k(x)$ the unique cube belonging to $\mathscr D_k$ such that $x\in Q_k(x)$. Then,
    \[
    \lim_{k\to-\infty}\mu\bigl(Q_k(x)\bigr)=\mu(X).
    \]
\end{lemma}

\begin{proof}
We consider the cases $\mu(X)<\infty$ and $\mu(X)=\infty$ separately.

\begin{itemize}[label={}, left=-4.5pt]
\item \textbf{Case \boldmath $\mu(X)<\infty$.}\\
In this case, by \cite[Lemma 5.1]{NY:97}, there exists $r_0>0$ such that $X\subset B(y,r_0)$ for all $y\in X$. Thus, for $k\in \Z$ sufficiently small, by \eqref{eq:cubeVsBalls} we have
\[
X\subset B(z_{Q_k(x)},r_0)\subset B(z_{Q_k(x)},c_0\delta^k)\subset Q_k(x)\subset X,
\]
which implies
\[
Q_k(x)=X.
\]
Therefore, trivially,
\[
\lim_{k\to-\infty}\mu\bigl(Q_k(x)\bigr)=\mu(X).
\]

\item \textbf{Case \boldmath $\mu(X)=\infty$.}\\
To address this case, we begin once again by recalling \eqref{eq:cubeVsBalls}, that is, $B(z_{Q_k(x)},c_0\delta^k)\subset Q_k(x)\subset B(z_{Q_k(x)},C_0\delta^k)$. So, since $x\in Q_k(x)$, we have in particular that
\[
\rho(x,z_{Q_k(x)})\leq C_0\delta^k.
\]
Thus, if $y\in B(x,c_0\delta^k)$, then
\[
\rho(y,z_{Q_k(x)})\leq A_0\bigl(\rho(y,x)+\rho(x,z_{Q_k(x)})\bigr)\leq A_0(c_0+C_0)\delta^k=: \tilde{C}\delta^k,
\]
which implies that
\[
B(x,c_0\delta^k)\subset B(z_{Q_k(x)},\tilde C\delta^k).
\]
Therefore, by the doubling property of $\mu$ and by its continuity from below, we may conclude that
    \[
\mu\bigl(Q_k(x)\bigr)\geq \mu\bigl(B(z_{Q_k(x)},c_0\delta^k)\bigr)
\gtrsim \mu\bigl(B(z_{Q_k(x)},\tilde C\delta^k)\bigr)
\geq \mu\bigl(B(x,c_0\delta^k)\bigr)
\xrightarrow[k\to -\infty]{} \mu(X)=\infty .
\]
    This concludes the proof. \qedhere
\end{itemize}    
\end{proof}

We now introduce some notation that will be used to describe the hierarchical relations between cubes at different levels of a dyadic system. 

Given $Q \in \mathscr D_k$ (where we interpret $Q$ as the pair $(Q,k)$), we define its (dyadic) parent by $Q^{(1)} := (Q,k)^{(1)}$, which denotes the unique cube $R \in \mathscr D_{k-1}$ such that $Q \subseteq R$. Likewise, we define the (dyadic) children of $Q$ as
\begin{equation*}
  \operatorname{ch}(Q) := \operatorname{ch}(Q,k) := \{Q' \in \mathscr D_{k+1} : (Q')^{(1)} = Q\}.
\end{equation*}
We extend these notions recursively by setting $Q^{(j)} := (Q^{(j-1)})^{(1)}$ and
$\operatorname{ch}^{(j)}(Q) := \operatorname{ch}(\operatorname{ch}^{(j-1)}(Q))$.
From the doubling condition, it follows that the number of children satisfies
\begin{equation*}
  M_Q := \#\operatorname{ch}(Q) \leq M < \infty
\end{equation*}
for some constant $M$ that does not depend on the cube $Q \in \mathscr D$.

\subsubsection*{Haar functions}
An important tool that will be used extensively throughout this work is the so-called Haar functions. For this reason, we recall its main properties. On the other hand, it is not necessary for our purposes to present its exact construction, and we refer the reader to \cite[Section 6]{Hyt:26} for an explicit formula.

For each $Q \in \mathscr D$, we define the non-cancellative wavelet as
\begin{equation*}
  h^0_Q := \mu(Q)^{-1/2} \mathbf{1}_Q.
\end{equation*}
On the other hand, if $M_Q>1$, it is possible to construct $M_Q-1$ cancellative wavelets $\{h_Q^i\}_{i=1}^{M_Q-1}$ with the following properties:
\begin{enumerate}
  \item Each $h_Q^i$ satisfies $\|h_Q^i\|_{L^2(\mu)} = 1$ for $i = 0, \ldots, M_Q - 1$.
  \item For each $i = 0, \ldots, M_Q - 1$, we have the pointwise bound
\begin{equation}\label{eq:absHaarBound}
    |h_Q^i| \lesssim \frac{\mathbf{1}_Q}{\sqrt{\mu(Q)}},
\end{equation}
which in turn implies that, for every $p \in [1, \infty]$,
\begin{equation}\label{eq:LpHaar}
    \Norm{h_Q^i}{L^p(\mu)} \lesssim \mu(Q)^{\frac{1}{p} - \frac{1}{2}}.
\end{equation}
In the latter formula, we adopt the convention that if $p = \infty$, then $\frac{1}{p} = 0$.
  \item For $0 \leq i < j \leq M_Q - 1$, the function $h_Q^i$ is constant on the support of $h_Q^j$.
  \item For $j = 1, \ldots, M_Q - 1$, we have $\int h_Q^j \, \ud\mu = 0$.
  \item The family $\{h_Q^i\}_{i=0}^{M_Q-1}$ forms an orthonormal basis of the space
\begin{equation*}
  L^2(\operatorname{ch}(Q)) := \left\{ f \in L^2(\mu) : f \text{ is constant on each } Q' \in \operatorname{ch}(Q), \; \mathbf{1}_{Q^c}f = 0 \right\}.
\end{equation*}
  \item The functions $\{h_Q^i\}_{i=1}^{M_Q-1}$ form an orthonormal basis of the subspace
\begin{equation}\label{eq:L20children}
  L^2_0(\operatorname{ch}(Q)) := \left\{ f \in L^2(\operatorname{ch}(Q)) : \int f \, \ud\mu = 0 \right\}.
\end{equation}  
\end{enumerate}
If $M_Q=1$, then there exists no cancellative wavelet.

\subsubsection*{Conditional expectations and martingale differences}

Define
\begin{equation}\label{eq:EQ}
  \E_Q f := \ave{f}_Q  \mathbf{1}_Q = \langle f, h_Q^0 \rangle h_Q^0,
\end{equation}
and let $\D_Q$ denote the orthogonal projection from $L^2(\mu)$ onto $L^2_0(\operatorname{ch}(Q))$. The latter space is defined in \eqref{eq:L20children}. Then we have the identities:
\begin{equation}\label{eq:DQ}
  \D_Q f = \sum_{Q' \in \operatorname{ch}(Q)} \E_{Q'}f - \E_Q f  = \sum_{i=1}^{M_Q - 1} \langle f, h_Q^i \rangle h_Q^i. 
\end{equation}

Next, define the following operators indexed by level~$k$:
\begin{equation}\label{eq:EkDk}
  \E_k := \sum_{Q \in \mathscr D_k} \E_Q, \qquad \D_k := \E_{k+1} - \E_k = \sum_{Q \in \mathscr D_k} \D_Q.
\end{equation}

In the following lemma, we examine the behavior of the operator $\E_k$ as $k \to \pm\infty$. 

\begin{lemma}\label{lem:convcondi}
   Let $p\in (1,\infty)$ and let $f \in L^p(w)$ for some weight $w\in A_p$. Then the following limits hold in $L^p(w)$:
\begin{align}
    &\lim_{k \to \infty}\E_k f=f, \label{eq:convcondi1}\\ 
    &\lim_{k \to -\infty}\E_kf=\ave{f}_X= \begin{cases} 0, & \mu(X) = \infty, \\ \fint_X f \, \ud\mu, & \mu(X) < \infty. \end{cases},\label{eq:convcondi2}
\end{align}
\end{lemma}

\begin{remark}
    Recall that, in particular, the function $w(x)\equiv 1$ is a weight belonging to $A_p$. Consequently, this result also yields, as a special case, convergence in $L^p(\mu)$.
\end{remark}

\begin{proof}[Proof of Lemma~\ref{lem:convcondi}.]
    Before proving the two convergences, recall that for every $x\in X$ and $k\in\mathbb{Z}$ there exists a unique cube in $\mathscr{D}_k$ containing $x$. We denote this cube by $Q_k(x)$, so that $\{Q_k(x)\}_{k\in\mathbb{Z}}$ is the family of dyadic cubes containing $x$, with $Q_k(x)\in\mathscr{D}_k$ for each $k$.

 Recall also (see, for instance, \cite[Lemma 2.2, p.~396]{GCRF:85}) that, since $w\in A_p$ and $(X,\rho,\mu)$ is a space of homogeneous type, the weighted space $(X,\rho,w(x)\ud\mu(x))$ is also a space of homogeneous type, i.e., the measure given by
 \[
 w(A):=\int_A w\ud\mu
 \]
is doubling.

We now prove the two convergences separately.
\begin{itemize}[label={}, left=-4.5pt]
\item \textbf{Proof of \boldmath\eqref{eq:convcondi1}.}

To begin with, since $Q_k(x)\downarrow\{x\}$ as $k\to\infty$ (by \eqref{eq:Q-BQ}) and $f\in L^p(w)\subset L^1_{\text{loc}}(\mu)$, by the Lebesgue differentiation theorem we have that, for $\mu$-almost every $x\in X$,
\[
\lim_{k\to\infty}\E_k f(x)=\lim_{k \to \infty}\Big(\frac{1}{\mu(Q_k(x))}\int_{Q_k(x)}f\ud\mu \Big)=f(x).
\]
On the other hand, if we denote by $M$ the Hardy--Littlewood maximal operator, we have that
\[
\abs{E_kf(x)-f(x)}^p=\abs{\fint_{Q_k(x)}f\ud\mu-f(x)}^p\lesssim \bigl(Mf(x)+|f(x)|\bigr)^p,
\]
with the function on the right-hand side belonging to $L^1(w)$, by virtue of the characterization of $A_p$ given in Theorem \ref{thm:MuckH-L}. Thus, by Lebesgue's dominated convergence theorem, we get that
\[
\lim_{k\to\infty}\|E_kf-f\|_{L^p(w)}=0.
\]

\item \textbf{Proof of \boldmath\eqref{eq:convcondi2}.}

For this case, we distinguish between the case where $\mu(X) < \infty$ and the one where $\mu(X) = \infty$.
\begin{itemize}
    \item \textbf{Case \boldmath$\mu(X) < \infty$.}

In this case, by \cite[Lemma 5.1]{NY:97}, there exists $r_0>0$ such that $X\subset B(x,r_0)$ for all $x\in X$. In particular, by \eqref{eq:Q-BQ}, this implies that there exists $k_0\in\Z$ such that, for all $k\leq k_0$, the family \(\mathscr{D}_k\) of dyadic cubes of the \(k\)-th generation consists of a single element \(Q = X\). Therefore, for all $x\in X$ and $k\leq k_0$, 
\[
\E_k f(x)=\frac{1}{\mu(Q_k(x))}\int_{Q_k(x)}f\ud\mu=\frac{1}{\mu(X)}\int_{X}f\ud\mu=\ave{f}_X,
\]
which trivially implies that
\[
\lim_{k\to-\infty}\|E_kf-\ave{f}_X\|_{L^p(w)}=0.
\]
    \item \textbf{Case \boldmath$\mu(X) = \infty$.}\\
Since $(X,\rho,\mu)$ is a space of homogeneous type with infinite measure, it necessarily follows that $\operatorname{diam}(X)=\infty$. Indeed, if this were not the case, then $X$ would be entirely contained in a ball of finite radius, which by the very definition of a space of homogeneous type must have finite measure. Therefore, since $(X,\rho,w(x)\ud\mu(x))$ is also a space of homogeneous type, the non-finiteness of the diameter of $X$ implies that $w(X)=\infty$  by \cite[Lemma 5.1]{NY:97}. Moreover, we can apply Lemma~\ref{lem:CubesCover}, which yields
\[
\lim_{k \to -\infty} w(Q_k(x)) = w(X) = \infty,
\]
for every $x \in X$. On the other hand, by Hölder's inequality we get that, for all $k\in\Z$,
\[
\abs{\int_{Q_k(x)}f\ud\mu}\leq \Big(\int_{Q_k(x)}\abs{f}^pw\ud\mu\Big)^{\frac1p}\Big(\int_{Q_k(x)}w^{-\frac{1}{p-1}}\ud\mu\Big)^{\frac{1}{p'}}\leq \Norm{f}{L^p(w)}\Big(\int_{Q_k(x)}w^{-\frac{1}{p-1}}\ud\mu\Big)^{\frac{1}{p'}}.
\]
Therefore, for all $x\in X$, 
\begin{equation*}
    \begin{split}
        \lim_{k \to -\infty}\abs{\E_kf(x)}&\leq \Norm{f}{L^p(w)} \lim_{k \to -\infty}\Bigg[\frac{1}{\mu(Q_k(x))}\Big(\int_{Q_k(x)}w^{-\frac{1}{p-1}}\ud\mu\Big)^{\frac{1}{p'}}\Bigg]\\
        &=\Norm{f}{L^p(w)} \lim_{k \to -\infty}\Bigg[\Bigg(\fint_{Q_k(x)}w\ud\mu\Big(\fint_{Q_k(x)}w^{-\frac{1}{p-1}}\ud\mu\Big)^{p-1}\Bigg)^{\frac1p} \Big(\int_{Q_k(x)}w\ud\mu\Big)^{-\frac1p}\Bigg]\\
        &\leq \Norm{f}{L^p(w)} \Big([w]_{A_p}\Big)^{\frac1p} \lim_{k \to -\infty} w(Q_k(x))^{-\frac1p}=0.
    \end{split}
\end{equation*}
Finally, since $\abs{\E_kf}^p\lesssim (Mf)^p\in L^1(w)$,
we may apply once again Lebesgue's dominated convergence theorem to get that
\[
\lim_{k\to-\infty}\|E_kf\|_{L^p(w)}=0.
\]
This concludes the proof. \qedhere
    \end{itemize}
\end{itemize}
\end{proof}

\subsubsection*{Random dyadic cubes}
Just as in the case of the dyadic formula for singular integrals, we will make use of random cubes also in the setting of fractional integrals. In particular, we will need to recall essentially the same properties stated in \cite[Section 17]{Hyt:26}, where the results of \cite{AH:13} and \cite{Hyt:Enflo} are adapted to this framework. It is worth mentioning that, historically, an initial version in the Euclidean context was given in \cite{NTV:Tb}, followed by a first construction in spaces of homogeneous type in \cite{HM:12}, with later variants appearing in \cite{AH:13,HK:12,NRV}.

To begin with, reference points $x^k_\alpha$ are given for each $k \in \mathbb{Z}$ and $\alpha \in \mathscr{A}_k$, where each $\mathscr{A}_k$ is a countable index set. We consider a parameter space $\Omega = \prod_{k \in \mathbb{Z}} \Omega_k$, where each $\Omega_k$ is a copy of a fixed finite set $\Omega_0$. This defines a canonical probability measure on $\Omega$.

For each $\omega \in \Omega$, we define a partial ordering $\leq_\omega$ on pairs of the form $(k, \alpha)$ and $(\ell, \gamma)$. At the same level, we say $(k, \alpha) \leq_\omega (k, \beta)$ if and only if $\alpha = \beta$. The order relation between levels $k$ and $k+1$ depends only on the component $\omega_k$, and the ordering between arbitrary levels $k \leq \ell$ is obtained by transitivity.

There are also new random reference points $z^k_\alpha = z^k_\alpha(\omega_k)$ such that
\begin{equation}\label{eq:z-x}
  d(z^k_\alpha(\omega), x^k_\alpha) < 2A_0 \delta^k.
\end{equation}

By \cite[Theorem 2.11 and Lemma 2.12]{AH:13}, there exist open and closed sets $\tilde{Q}^k_\alpha(\omega)$ and $\bar{Q}^k_\alpha(\omega)$ that satisfy the inclusions
\begin{equation}\label{eq:Q-BQ}
  B(z^k_\alpha, \tfrac{1}{6} A_0^{-5} \delta^k) \subseteq \tilde{Q}^k_\alpha(\omega) \subseteq \bar{Q}^k_\alpha(\omega) \subseteq B(z^k_\alpha, 6A_0^4 \delta^k) 
  \subseteq B(x^k_\alpha, 8A_0^5 \delta^k),
\end{equation}
where the last inclusion follows directly from \eqref{eq:z-x} and the quasi-triangle inequality. 
In addition, for each fixed $\omega \in \Omega$, the associated families also satisfy some useful disjointness and covering properties.

For the purposes of this paper, however, we will not work directly with these open and closed cubes. Instead, we consider the associated family of \emph{half-open} cubes $Q^k_\alpha(\omega)$, constructed as in \cite[Lemma 2.18]{HK:12}, which satisfy
\begin{equation}\label{eq:hocubes}
  \tilde{Q}^k_\alpha(\omega) \subseteq Q^k_\alpha(\omega) \subseteq \bar{Q}^k_\alpha(\omega),
\end{equation}
and form a genuine dyadic system on $X$ for each fixed $\omega$, in the sense of Definition~\ref{def:cubes}.

It is worth mentioning that, if one follows the procedure of \cite[Lemma 2.18]{HK:12}, the useful property that the cubes at level \(k\) depend only on the tail \((\omega_j)_{j \geq k}\) is lost. In fact, they would depend on \((\omega_j)_{j \geq \min\{k,k_0\}}\) for some fixed \(k_0\). Nevertheless, this dependence on coarser levels affects only the boundaries of the cubes, as points on the boundary are assigned to one cube or another in order to ensure the nestedness property. By \cite[Lemma 17]{Christ:90}, we know that \(\mu(\partial Q^k_\alpha) = 0\) for every cube. Therefore, for the purposes of this work, we may interpret these cubes not as sets in the strict sense, but rather as equivalence classes up to sets of measure zero. Accordingly, we can regard the Haar wavelets as equivalence classes of functions, rather than as pointwise-defined functions. In this sense, we may thus consider this family of random dyadic cubes and treat them, for each fixed \(\omega\), in the same way as the dyadic cubes of any given dyadic system, assigning to them the independence property that will be needed in certain parts of the analysis. 

Another aspect to keep in mind throughout the discussion is that these random dyadic cubes $Q^k_\alpha(\omega)$ can be regarded as balls of radius $\delta^k$, centered either at the deterministic reference points $x_k^\alpha$ or at the random reference points $z_k^\alpha(\omega)$. Indeed, from \eqref{eq:Q-BQ} and \eqref{eq:hocubes}, together with the doubling property of $\mu$, we obtain
\begin{equation}
\mu(Q^k_\alpha(\omega)) \approx V(z_k^\alpha(\omega), \delta^k) \nonumber
\end{equation}
and furthermore, from \eqref{eq:z-x}, again using the doubling property of $\mu$, we also have
\begin{equation}\label{eq:CubesBallsEquiv}
V(x_k^\alpha, \delta^k) \approx V(z_k^\alpha(\omega), \delta^k) \approx \mu(Q^k_\alpha(\omega)).
\end{equation}

From now on, to simplify the notation, we will refer to a generic cube $Q^k_\alpha(\omega)$ simply as $Q$, and we will just call it a \textit{dyadic cube}. We define the collection of dyadic cubes at level $k$ by $\mathscr{D}_k = \mathscr{D}_k(\omega) := \{ Q^k_\alpha(\omega) : \alpha \in \mathscr{A}_k \}$. Moreover, if a cube $Q$ belongs to $\mathscr{D}_k$, we define its sidelength by $\ell(Q) := \delta^k$. 

\section{Fractional integrals}\label{sec:fracint}
In this section, we introduce the notion of fractional integral operator and recall some of its basic properties. We also present the example arising from the Bessel setting, which serves as a motivating case widely studied in the literature. We then turn to the analysis of their boundedness and conclude with estimates for the associated Haar coefficients.

\subsection{Volumic and metric fractional integrals}
First of all, we need to define what we mean by fractional integral operators. Roughly speaking, their definition is obtained by slightly modifying the definition of a singular integral operator (see, for instance, (1.5) and (1.6) in \cite{Hyt:26}), making the kernel ``a bit less singular''. In order to do so, given the general context of spaces of homogeneous type, there is no unique way to define them. In the case of $\R^d$ equipped with the $d$-dimensional Lebesgue measure $\operatorname{Leb}_d$, it is well known that the measure of balls and their respective radii are related in the following way:
\[
\operatorname{Leb}_d \bigl(B(x,r)\bigr)\approx r^d, \quad x\in\R^d, \, r>0.
\]
Thus, in particular,
\[
V(x,y)=\operatorname{Leb}_d \bigl(B(x,|x-y|)\bigr)\approx |x-y|^d.
\]
This property clearly extends to any Ahlfors $d$-regular metric measure space. However, this argument does not apply to an arbitrary space of homogeneous type $(X,\rho,\mu)$, since, in general, we do not have such a proportionality between the radius and the measure of a ball. For this reason, it makes sense to introduce two possible definitions of fractional integral. 

\begin{definition}\label{def:fracint}
Let $(X,\rho,\mu)$ be a quasi-metric measure space, and let $K \in L^1_{\mathrm{loc}}(X \times X)$. Consider the operator $T$ formally given by
\[
Tf(x)=\int_X K(x,y)f(y)\,\ud\mu(y), \qquad x\in X.
\]
Let $\omega:[0,1]\to[0,\infty)$ be such that
\[
\lim_{t\to0^+}\omega(t)=0 \quad \text{and} \quad \omega(0)=0.
\]

\medskip
Let $\phi:X\times X\to[0,\infty)$ be a nonnegative function belonging to one of the admissible classes
\[
\phi(x,y)\in\{\,V(x,y)^{\varepsilon},\ \rho(x,y)^{\eta}\,\},
\qquad \varepsilon\in(0,1),\ \eta\in(0,\infty).
\]

\begin{itemize}

\item \textbf{Fractional kernel and operator.}\\
We say that $K$ is a \textit{$\phi$-fractional kernel}, and that $T$ is a \textit{$\phi$-fractional integral operator}, if $K$ satisfies the \textit{$\phi$-size estimate}
\begin{equation}\label{eq:CZ0phi}
|K(x,y)| \lesssim \frac{\phi(x,y)}{V(x,y)},
\end{equation}
uniformly for every $x,y\in X$.

Moreover, $K$ (and hence $T$) is said to be \textit{$\omega$-regular} if it also satisfies the \textit{$\omega,\phi$-smoothness estimate}
\begin{equation}\label{eq:CZ1phi}
|K(x,y)-K(x',y)|+|K(y,x)-K(y,x')|
\lesssim 
\omega\!\left(\frac{\rho(x,x')}{\rho(x,y)}\right)\frac{\phi(x,y)}{V(x,y)}
\end{equation}
whenever $\rho(x,x')\ll \rho(x,y)$. In both inequalities, we use the convention $\frac{0}{0}:=0$ on the right-hand side.

\vspace{0.3cm}

\item \textbf{Volumic case.}\\
In the case $\phi(x,y)=V(x,y)^{\varepsilon}$ with $\varepsilon\in(0,1)$, we refer to $K$ and $T$ as an \textit{$\varepsilon$-volumic fractional kernel} and an \textit{$\varepsilon$-volumic fractional integral operator}, respectively.

\vspace{0.3cm}

\item \textbf{Metric case.}\\
In the case $\phi(x,y)=\rho(x,y)^{\eta}$ with $\eta\in(0,\infty)$, we refer to $K$ and $T$ as an \textit{$\eta$-metric fractional kernel} and an \textit{$\eta$-metric fractional integral operator}, respectively.

\end{itemize}
\end{definition}

In this definition, we distinguish the two notions by referring to them as \textit{volumic} and \textit{metric}. This terminology reflects the fact that the estimates \eqref{eq:CZ0phi} and \eqref{eq:CZ1phi} modify the singularity \(V(x,y)^{-1}\) of a singular integral kernel through a multiplicative factor which, in the first case, depends on the ``volume'' \(V(x,y)\), and in the second case depends on the (quasi-)metric \(\rho(x,y)\). On the other hand, it should be noted that this choice is not universal, and different, possibly lighter, nomenclature may be adopted outside the present work. Furthermore, observe that in the case of spaces satisfying Ahlfors $d$-regularity for some $d > 0$ (i.e., if $V(x,r)\approx r^d$ uniformly in $x\in X$ and $r>0$), the metric definition is entirely equivalent to the volumic one, provided we restrict ourselves to $\eta \in (0,d)$ and identify $\eta = \varepsilon d$. If, on the other hand, the space under consideration does not satisfy this property, then the two definitions are in general not comparable and give rise to two genuinely different types of operators. This will become clearer to the reader by consulting the concrete examples constructed in Subsection~\ref{subsec:Bessel} and in Appendices~\ref{app:Slab} and~\ref{app:Cantor}.

Although these two classes of operators exhibit even substantially different properties in general, in most of the analysis carried out in this work they can be handled in a unified way. For this reason, it is particularly convenient to introduce a notation that encodes both definitions into a single formalism, as already suggested in Definition~\ref{def:fracint}. To this end, with a slight abuse of notation, we extend the meaning of the symbol $\phi$ as follows.

\begin{notation}\label{not:phi}
From now on, the formalism encoded by $\phi$ will be interpreted in two
different ways, depending on the type of fractional integral under
consideration. The possible interpretations are listed below. In each case,
the first line corresponds to the volumic definition, for some
$\varepsilon\in(0,1)$, while the second corresponds to the metric one,
for some $\eta\in(0,\infty)$.

For every $x,y\in X$, $r>0$, and $Q\in\mathscr{D}$, we set

\[
\begin{aligned}
\phi(x,y) &:= 
\begin{cases}
V(x,y)^\varepsilon\\
\rho(x,y)^\eta
\end{cases}
\qquad
\phi(x,r) :=
\begin{cases}
V(x,r)^\varepsilon\\
r^\eta
\end{cases}
\qquad
\phi(Q) :=
\begin{cases}
\mu(Q)^\varepsilon\\
\ell(Q)^\eta
\end{cases}
\qquad
\begin{aligned}
\text{(volumic)}\\
\text{(metric)}
\end{aligned}
\end{aligned}
\]
    To avoid redundancy, we will often omit explicit mention of the parameters \( \varepsilon \) and \( \eta \), when no confusion arises.
\end{notation}

\begin{remark}
    Notice that \( \phi(x,y) \equiv \phi(x,\rho(x,y)) \). This slight abuse of notation mirrors the analogous convention used for the function \( V(x,y) \equiv V(x, \rho(x,y)) \), and should therefore be familiar to the reader. We also note that, in the metric case, $\phi(x,r)$ actually does not depend on $x$. However, this minor discrepancy does not lead to any significant differences in the treatment, and thus can be ignored in most cases.
\end{remark}

\begin{remark}
    The fact that \(\phi\) may also be applied to a dyadic cube is a further mild abuse of notation. However, this should not cause any ambiguity, since this is the only instance in which \(\phi\) is not applied to a pair of objects, but rather to a single one. Moreover, it is also the only case in which it is applied to a set.
\end{remark}

\begin{remark}
    In the metric case, it is worth recalling that a dyadic cube $Q\in\mathscr D$ may belong to more than one dyadic level $\mathscr D_k$. In this case, we shall interpret $\phi(Q)$ as the application of $\phi$ to the pair $(Q,k)$, rather than to $Q$ alone as a set. Another point to keep in mind is that, when $\mu(X)<\infty$, we have $\mathscr D_k=\{X\}$ for every $k$ below a certain threshold. In this situation, we adopt the convention that $\ell(X)=\delta^k$, where
\[
k = k_{\max}(X):=\sup\{ h\in\Z: X\in\mathscr D_h\}\in\Z.
\]
Finally, we observe that if $(X,\rho,\mu)$ has positive lower dimension, then \cite[Lemma~6.3]{Hyt:26} ensures that each cube $Q\in\mathscr D$ belongs to exactly one level $\mathscr D_k$, so that its side length is uniquely determined.
\end{remark}

\begin{remark}
    To help the reader follow the arguments involving the abstract notation $\phi$, we collect here some of its fundamental properties, which will be used later on. Their proof is a straightforward consequence of the definition of $\phi$, together with the doubling property of the measure $\mu$. 
    \begin{enumerate}
        \item For all $x,y\in X$,
        \begin{equation}\label{eq:RemarkNotation1}
            \phi(x,y)\approx \phi(y,x).
        \end{equation}
        \item For all $x,y\in X$ and $r\gtrsim \rho(x,y)$,
        \begin{equation}\label{eq:RemarkNotation2}
            \phi(x,r)\approx \phi(y,r).
        \end{equation}
        \item For all $x,y\in X$ and $r\approx\rho(x,y)$,
        \begin{equation}\label{eq:RemarkNotation3}
            \phi(x,y)\approx \phi(x,r).
        \end{equation}
        \item For every random dyadic cube $Q=Q^k_\alpha(\omega)$ with deterministic reference point $x_Q=x_k^\alpha$ and random reference point $z_Q=z_k^\alpha(\omega)$, by \eqref{eq:CubesBallsEquiv} we have 
        \begin{equation}\label{eq:RemarkNotation4}
\phi(Q) =
\left\{
\begin{array}{ll}
\mu(Q)^\varepsilon \approx V(z_Q,\ell(Q))^\varepsilon \approx V(x_Q,\ell(Q))^\varepsilon , & \text{(volumic case)}\\[4pt]
\ell(Q)^\eta, & \text{(metric case)}
\end{array}
\right\}
\approx \phi(x_Q,\ell(Q)).
\end{equation}
\item If $P, Q \in \mathscr{D}$ are such that $P \subseteq Q$, that is, if $P \in \operatorname{ch}^{(j)}(Q)$ for some $j$, then
\begin{equation}\label{eq:RemarkNotation5}
    \phi(P)\leq\phi(Q)
\end{equation}
    \end{enumerate}
\end{remark}

In addition to these elementary properties, we prove a slightly less trivial technical lemma that will be useful later on. 

\begin{lemma}\label{lem:sumcubesfrac}
    Let $(X,\rho,\mu)$ be a space of homogeneous type, let $\phi$ be a fractional scaling function as in Notation~\ref{not:phi}, and let $x,y\in X$. In the metric setting, assume additionally that $(X,\rho,\mu)$ has positive lower dimension $d>0$ and that the metric fractional parameter satisfies $\eta\in(0,d)$. Then the following estimates hold:
    \begin{enumerate}[label=(\roman*)]
        \item\label{point:sumcubesxneqy} If $x\neq y$, then
    \begin{equation}\label{eq:sumcubesxneqy}
         \sum_{\substack{R\in\mathscr D\\R\ni x,y}} \frac{\phi(R)}{\mu(R)}\lesssim \frac{\phi(x,y)}{V(x,y)}
        \end{equation}
        \item\label{point:sumcubesx} In the volumic case, if $x=y$ and $\mu(\{x\})>0$, then
    \begin{equation}\label{eq:sumcubesx}
         \sum_{\substack{R\in\mathscr D\\R\ni x}} \frac{\phi(R)}{\mu(R)}\lesssim \frac{1}{\mu(\{x\})^{1-\eps}}=\frac{\phi(x,y)}{V(x,y)}
        \end{equation}
    \end{enumerate}
\end{lemma}

\begin{proof}
    We prove the two points separately.
    \begin{itemize}[label={}, left=-4.5pt]
\item \textbf{Point \ref{point:sumcubesxneqy}.}\\
Let $x\neq y$ and let $P_{x,y}$ be the smallest element of $\mathscr D$ containing both $x$ and $y$. Then,
\[
\sum_{\substack{R\in\mathscr D\\R\ni x,y}} \frac{\phi(R)}{\mu(R)}=\sum_{k=0}^\infty \frac{\phi(P_{x,y}^{[k]})}{\mu(P_{x,y}^{[k]})}.
\]
Now, in the volumic case, by \cite[Formula (6.9)]{Hyt:26}
\begin{equation*}
        \sum_{k=0}^\infty \frac{\phi(P_{x,y}^{[k]})}{\mu(P_{x,y}^{[k]})}=\sum_{k=0}^\infty \frac{1}{\mu(P_{x,y}^{[k]})^{1-\eps}}\lesssim \frac{1}{\mu(P_{x,y})^{1-\eps}}\approx \frac{1}{V(x,y)^{1-\eps}}=\frac{\phi(x,y)}{V(x,y)}.
\end{equation*}
On the other hand, in the metric case, by lower dimension $d>0$ we have
\[
\frac{\mu(P_{x,y}^{[k]})}{\mu(P_{x,y})}\gtrsim (\delta^{-k})^d,
\]
and therefore
\begin{equation*}
        \sum_{k=0}^\infty \frac{\phi(P_{x,y}^{[k]})}{\mu(P_{x,y}^{[k]})}=\sum_{k=0}^\infty \frac{(\delta^{-k})^\eta\ell(P_{x,y})^\eta}{\mu(P_{x,y}^{[k]})}\lesssim \frac{\ell(P_{x,y})^\eta}{\mu(P_{x,y})}\sum_{k=0}^\infty (\delta^{d-\eta})^k\approx \frac{\phi(x,y)}{V(x,y)}.
\end{equation*}
\item \textbf{Point \ref{point:sumcubesx}.}\\
Let $x\in X$ be such that $\mu(\{x\})>0$. Then, by \cite[Theorem 1]{MS:79}, we have that $x$ is an isolated point, i.e., that there exists $r>0$ such that $B(x,r)=\{x\}$. Therefore, if we denote by $P_k(x)$ the only element of $\mathscr D_k$ such that $x\in P_k(x)$, then $P_k(x)=\{x\}$ for all $k\geq k_x$, for some $k_x\in\Z$. Therefore, if we denote by $P_x:=\{x\}$, we have that
\[
\sum_{\substack{R\in\mathscr D\\R\ni x}} \frac{\phi(R)}{\mu(R)}=\sum_{k=0}^\infty \frac{\phi(P_{x}^{[k]})}{\mu(P_{x}^{[k]})}=\sum_{k=0}^\infty \frac{1}{\mu(P_{x}^{[k]})^{1-\eps}}\lesssim \frac{1}{\mu(P_x)^{1-\eps}}=\frac{1}{\mu(\{x\})^{1-\eps}}=\frac{\phi(x,y)}{V(x,y)},
\]
where we have once again applied \cite[Formula (6.9)]{Hyt:26}. This concludes the proof. \qedhere
\end{itemize}
\end{proof}

With the aim of streamlining the exposition and avoiding frequent repetitions, we conclude this subsection by introducing a list of assumptions that we will often require on the function \(\omega\) appearing in Definition~\ref{def:fracint}.

\begin{assumption}\label{ass:omega}
Whenever we invoke this assumption, we mean that the function $\omega \colon [0,1] \to [0,\infty)$ satisfies the following conditions:
\begin{itemize}
    \item $\omega$ is continuous;
    \item $\omega$ is non-decreasing;
    \item $\omega$ is doubling, in the sense that $\omega(2t) \lesssim \omega(t)$ for all $t \in [0,1/2]$;
    \item $\lim_{t \to 0^+} \omega(t) = 0 = \omega(0)$.
\end{itemize}
\end{assumption}

\subsection{The Bessel setting}\label{subsec:Bessel}

Before turning to the general treatment of these operators, we present here the fractional integral operator associated with the Bessel operator, in the Bessel setting, which is a widely studied and discussed example in the literature (see, e.g., \cite{BDLLW:25}). This setting provides a useful example of a space of homogeneous type in which the metric definition can be seen, with an appropriate choice of the parameters, as a special case of the volumic definition, but not vice versa (see Proposition~\ref{prop:Bessel}). 


    We begin by considering the space of homogeneous type given by $X := \R_+ := (0, \infty)$, equipped with the Euclidean metric, that is, the absolute value $|\cdot|$, and the measure $m_\lambda$ defined, for $\lambda > 0$, as follows:
\[
d m_\lambda := x^{2\lambda} \, dx.
\]
The only non-trivial property to verify is that the measure is doubling: this is ensured by the formula (1.2) from \cite{BDLLW:25}. Next, let us consider the Bessel operator $\Delta_\lambda$, which is defined by
\[
\Delta_\lambda := -\frac{d^2}{dx^2} - \frac{2\lambda}{x} \frac{d}{dx}.
\]
For $\alpha \in (0,1)$, the fractional operator $\Delta_\lambda^{-\alpha/2}$ has a kernel $K^\alpha(x,y)$ which, according to Lemmas 5.1 and 5.2 from \cite{BDLLW:25}, satisfies the estimates
\begin{align}
     &K^\alpha(x,y) \approx \frac{\abs{x-y}^\alpha}{m_\lambda(I(x,\abs{x-y}))}, \quad x,y\in \R_+;\label{eq:BesSize}\\
     &\abs{K^\alpha(x,y)-K^\alpha(x',y)} \lesssim \frac{\abs{x-x'}}{\abs{x-y}}\frac{\abs{x-y}^\alpha}{m_\lambda(I(x,\abs{x-y}))}, \quad x,x',y\in \R_+ \text{ with }\abs{x-x'}<\frac{\abs{x-y}}{2},\label{eq:BesSmooth}
   \end{align}
where we denote
\[
I(x,r):=(x-r,x+r)\cap\R_+, \quad x\in R_+, r>0.
\]
Moreover, analogous estimates also hold in the more general setting of $\R^{n+1}_+$, as shown in \cite[Section~12]{HW:25}. Next, by the explicit formula of $K^\alpha$ given in the proof of Lemma 5.2 in \cite{BDLLW:25}, we have in particular that $K^\alpha(x,y)= K^\alpha(y,x)$, so \eqref{eq:BesSmooth} implies that
\begin{equation}\label{eq:BesSmooth2}
    \abs{K^\alpha(x,y)-K^\alpha(x',y)}+\abs{K^\alpha(y,x)-K^\alpha(y,x')} \lesssim \frac{\abs{x-x'}}{\abs{x-y}}\frac{\abs{x-y}^\alpha}{m_\lambda(I(x,\abs{x-y}))},
\end{equation}
for $x,x',y\in \R_+$ such that $\abs{x-x'}<\frac{\abs{x-y}}{2}$. Therefore, $K^\alpha$ is an $\omega$-regular $\alpha$-metric fractional kernel, in the sense of Definition~\ref{def:fracint}, with $\omega(t)=t$. In fact, we can say even more, as \eqref{eq:BesSize} provides an equivalence rather than merely an upper bound. In this sense, the metric definition of fractional integral appears to most accurately reflect the structure of this operator, particularly regarding the size estimate, which is sharp: no genuinely smaller upper bound is attainable.

Interestingly, we can show that $\Delta_\lambda^{-\alpha/2}$ also fits into the volumic definition of fractional integral. This is a consequence of the following result, which shows that in this setting the metric condition is, in some sense, strictly stronger than the volumic one. 

\begin{proposition}\label{prop:Bessel}
    In the Bessel setting, any ($\omega$-regular) $\alpha$-metric fractional integral, with $\alpha\in(0,2\lambda+1)$, is, in particular, an ($\omega$-regular) $\frac{\alpha}{2\lambda + 1}$-volumic fractional integral. However, given a volumic fractional integral, it is not in general possible to view it as a special case of the metric definition, for any choice of parameters.
\end{proposition}
\begin{proof}\
\begin{itemize}[label={}, left=-4.5pt]
        \item \textbf{Metric \boldmath$\Rightarrow$ volumic.}\\
    Let $\alpha\in(0,2\lambda+1)$ and let $K$ be an $\alpha$-metric fractional kernel. Then, from the formula at the beginning of the proof of Proposition 4.2 in \cite{Hyt:26}, we have
\begin{equation}\nonumber
    m_\lambda(I(x,r)) \approx r \max(x, r)^{2\lambda}\geq r^{2\lambda+1}, \quad x \in \R_+, \, r > 0, 
\end{equation}
so that in particular
\begin{equation}\label{eq:BesInterval}
    r\lesssim m_\lambda(I(x,r))^{\frac{1}{2\lambda+1}}, \quad x \in \R_+, \, r > 0.
\end{equation}
Thus, starting with the size estimate \eqref{eq:CZ0phi}, by \eqref{eq:BesInterval} we obtain
\begin{equation}\nonumber
    K(x,y) \lesssim \frac{|x - y|^\alpha}{m_\lambda(I(x, |x - y|))} \lesssim \frac{1}{m_\lambda(I(x, |x - y|))^{1-\frac{\alpha}{2\lambda+1}}}=\frac{1}{V(x,y)^{1 - \frac{\alpha}{2\lambda+1} }}, 
\end{equation}
uniformly in $x,y\in\R_+$. Similarly, the smoothness estimate \eqref{eq:CZ1phi} and \eqref{eq:BesInterval} imply
\begin{equation}\nonumber
    \abs{K(x,y)-K(x',y)}+\abs{K(y,x)-K(y,x')}
    \lesssim \omega\Big(\frac{\abs{x-x'}}{\abs{x-y}}\Big)\frac{1}{V(x,y)^{1 - \frac{\alpha}{2\lambda+1} }},
\end{equation}
for every $x, x', y \in \R_+$ such that $|x - x'| \ll |x - y|$. It follows that, if we choose
\[
\varepsilon := \frac{\alpha}{2\lambda + 1} \in \left(0, \frac{1}{2\lambda + 1}\right) \subset (0,1),
\]
then $K$ is also an $\omega$-regular $\varepsilon$-volumic fractional kernel.
\item \textbf{Volumic \boldmath$\cancel\Rightarrow$ metric.}\\
Suppose for contradiction that there exist a constant $c_1>0$ and $\beta \in (0,\infty)$ such that
\[
V(x,y)\leq c_1 \rho(x,y)^\beta,
\]
uniformly in $x,y\in X$. Then, for every $x,y$ such that $x\geq 2\rho(x,y)$, we would have
\[
c_1\rho(x,y)^\beta \geq  V(x,y)=\int_{x-\rho(x,y)}^{x+\rho(x,y)}u^{2\lambda} \, du\approx x^{2\lambda}\rho(x,y),
\]
which implies
\[
x^{2\lambda}\leq c_2\rho(x,y)^{\beta-1},
\]
for some constant $c_2>0$ independent of $x,y$. This last inequality clearly cannot hold uniformly in $x,y$ with $x \geq 2\rho(x,y)$. For instance, setting $y = x+1$ and letting $x \to \infty$ shows the failure. \qedhere
\end{itemize}
\end{proof}


\subsection{Boundedness}
\label{subsec:fracbdd}
A highly relevant aspect in the study of fractional integrals is their boundedness on $L^p$ spaces. This will therefore be the content of this subsection. We will start by recalling the results already known in the literature regarding the volumic operators, and then shift the focus to those of the metric type. For the latter definition, as we shall see, one can obtain a slightly weaker boundedness result than in the volumic case, and only under the additional assumption of positive lower dimension on the underlying space.

\subsubsection*{Volumic type}
It is well known in the literature \cite{BC:95, GV:90} that volumic fractional integrals are bounded from $L^p$ to $L^q$, where $p$ and $q$ satisfy a specific relation. This result, which we now state, is \cite[Theorem 4.1]{BC:95}.


\begin{theorem}\label{thm:bddfracvol}
    Let $(X, \rho, \mu)$ be a space of homogeneous type and let $T$ be an $\eps$-fractional integral operator for some $\eps \in (0,1)$. Then $T$ is bounded from \( L^p(X,\mu) \) to \( L^q(X,\mu) \) for all $p,q\in(1,\infty)$ such that $1/q=1/p-\eps>0$.
\end{theorem}

\begin{remark}
    In the Bessel setting, illustrated in Subsection~\ref{subsec:Bessel}, the result just stated automatically provides us with some information regarding the boundedness of metric-type operators, including the one associated with the Bessel operator. Indeed, by Proposition~\ref{prop:Bessel}, in this framework any $\alpha$-metric fractional integral, with $\alpha\in(0,2\lambda+1)$, is, in particular, an $\frac{\alpha}{2\lambda + 1}$-volumic fractional integral. Thus, by Theorem~\ref{thm:bddfracvol}, they are bounded from \( L^p(\R_+,m_\lambda) \) to \( L^q(\R_+,m_\lambda) \) for all \( p, q \in (1, \infty) \) such that \( \frac1q = \frac1p - \frac{\alpha}{2\lambda + 1} > 0 \).
\end{remark}

\subsubsection*{Metric type}

We have just seen in the previous paragraph that results concerning the $L^p-L^q$ boundedness of volumic-type fractional integrals are well established in the literature.  
Therefore, since these operators are ``morally'' similar to the metric definition, one might expect that, with some adjustments, a similar result could be proven for the metric-type ones.  
However, in the case of a generic space of homogeneous type, this is not possible: in support of this statement, in Appendix~\ref{app:Slab} we construct a concrete example of a space of homogeneous type (with positive lower and upper dimensions) in which metric fractional integrals are unbounded from $L^p$ to $L^q$ for any choice of $p,q \in (1,\infty)$. This therefore motivates the search for a less conventional form of boundedness, which turns out to be possible under the assumption of a positive lower dimension for the underlying space. The study of this property is the subject of the remainder of this subsection.

Before presenting the main result, we first show that a certain function is a Muckenhoupt weight, according to Definition \ref{def:Muckweight}. We provide here a complete proof, although in the Euclidean setting, for the so-called power weights, the conditions characterizing the membership in this class are well known (see, for instance, \cite[Chapter IV, Section 2, just before Theorem 2.16]{GCRF:85}, or \cite[Example 7.17]{Grafakos:CFA}).

\begin{lemma}\label{lem:Muck}
    Let $(X,\rho,\mu)$ be a space of homogeneous type with lower dimension $d\in(0,\infty)$. Moreover, let $x_0\in X$, $p\in(1,\infty)$, and $\eta\in\left(0,\frac dp\right)$. Then, the function 
    \[
    w(x):=\bigl(1+\rho(x,x_0)\bigr)^{-\eta p}, \quad x\in X,
    \]
    is a Muckenhoupt weight belonging to $A_p$.
\end{lemma}

\begin{proof}
    We want to show that $w$ satisfies
    \[
    \sup_{B \text{ ball}} 
  \left( \fint_B w \, d\mu \right) 
  \left( \fint_B w^{-\frac{1}{p-1}} \, d\mu \right)^{p-1} < \infty.
    \]
    Now, for any ball $B\subseteq X$, let us denote by $x_B\in X$ its center, $r_B>0$ its radius, and $r_0:=\rho(x_B,x_0)$. We consider two cases separately. 
    \begin{itemize}[label={}, left=-4.5pt]
        \item \textbf{Case} $\boldsymbol{\displaystyle{r_B\leq \frac{r_0}{2A_0}}}$.\\
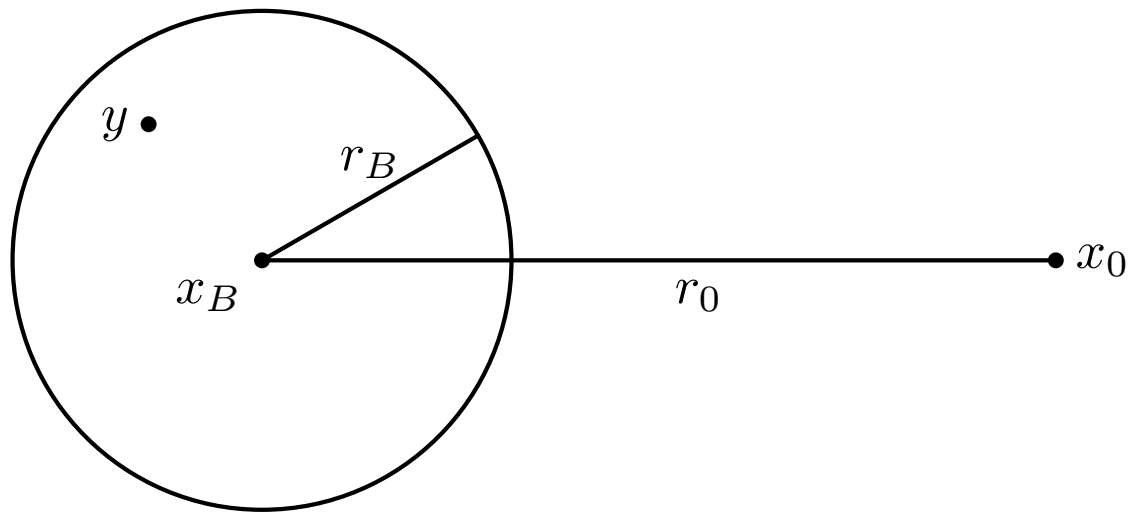
        \begin{figure}[H]
    \centering
    \begin{tikzpicture}[scale=1.5, thick]

        \coordinate (xB) at (0,0);
        \coordinate (x0) at (3.5,0);
        \coordinate (y) at (-0.5,0.6);

        \def\rB{1.1}
        \draw[black] (xB) circle (\rB);

        \path (xB) -- ++(30:\rB) coordinate (rBendpoint);
        \draw[thick] (xB) -- (rBendpoint);
        \node[above] at ($(xB)!0.5!(rBendpoint)$) {$r_B$};

        \draw[black!60!black, thick] (xB) -- (x0);
        \node[black!60!black, below] at ($(xB)!0.55!(x0)$) {$r_0$};

        \fill (xB) circle (1pt) node[below left] {$x_B$};
        \fill (x0) circle (1pt) node[right] {$x_0$};
        \fill (y) circle (1pt) node[left] {$y$};

    \end{tikzpicture}
    \end{figure}
    
    Let $B\subseteq X$ be a ball such that $r_B\leq \frac{r_0}{2A_0}$. Then, for $y\in B$, 
    \[
    \rho(x_B,x_0)\leq A_0\bigl(\rho(x_B,y)+\rho(y,x_0)\bigr),
    \]
    and so
    \[
    \rho(y,x_0)\geq \frac{1}{A_0}\rho(x_B,x_0)-\rho(x_B,y)\geq\frac{1}{A_0}r_0-r_B\geq\frac{r_0}{2A_0}.
    \]
    On the other hand,
    \begin{equation*}
        \rho(y,x_0)\leq A_0\bigl(\rho(y,x_B)+\rho(x_B,x_0)\bigr)\leq A_0(r_B+r_0)\leq \left(\frac{1}{2}+A_0\right)r_0.
    \end{equation*}
    Thus, altogether, for $y\in B$ we have that
    \[
    1+\rho(y,x_0)\approx 1+r_0.
    \]
    Therefore, for all balls $B\subseteq X$ such that $r_B\leq \frac{r_0}{2A_0}$, we have that
    \begin{align*}
        & \fint_B w \ud\mu \approx (1+r_0)^{-\eta p}, \\ 
  &\left( \fint_B w^{-\frac{1}{p-1}} \ud\mu \right)^{p-1}\approx (1+r_0)^{\eta p},
    \end{align*}
    which directly implies
    \[
    \left( \fint_B w \ud\mu \right) 
  \left( \fint_B w^{-\frac{1}{p-1}}\ud\mu\right)^{p-1}\approx 1.
    \]
    \item \textbf{Case} $\boldsymbol{\displaystyle{r_B\geq \frac{r_0}{2A_0}}}$. 
    
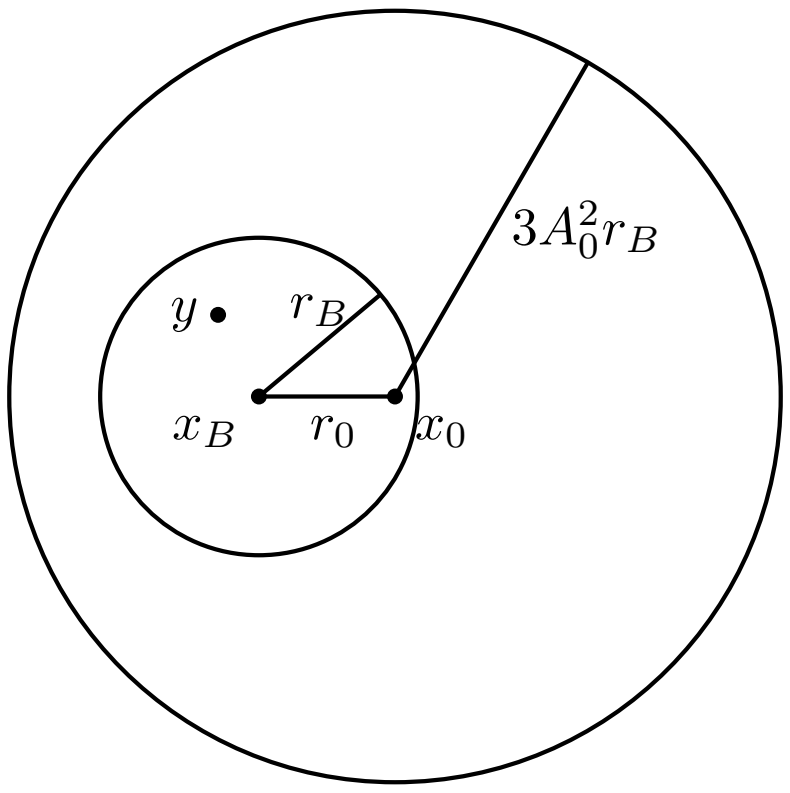
\begin{figure}[H]
    \centering
    \begin{tikzpicture}[scale=1.5, thick]

        \coordinate (xB) at (0.6,0);
        \coordinate (x0) at (1.2,0); 
        \coordinate (y) at (0.42,0.36);

        \def\rB{0.7}
        \def\A{0.9} 
        \def\rBig{3*\A*\A*\rB}

        \draw[black] (xB) circle (\rB);

        \draw[black, thick] (x0) circle (\rBig);

        \path (xB) -- ++(40:\rB) coordinate (rBendpoint);
        \draw[thick] (xB) -- (rBendpoint);
        \node[above] at ($(xB)!0.5!(rBendpoint)$) {$r_B$};

        \draw[black!60!black, thick] (xB) -- (x0);
        \node[black!60!black, below] at ($(xB)!0.55!(x0)$) {$r_0$};

        \fill (xB) circle (1pt) node[below left] {$x_B$};
        \fill (x0) circle (1pt) node[below right] {$x_0$};
        \fill (y) circle (1pt) node[left] {$y$};

        \path (x0) -- ++(60:\rBig) coordinate (rBigendpoint);
        \draw[black, thick] (x0) -- (rBigendpoint);
        \node[black, right] at ($(x0)!0.5!(rBigendpoint)$) {$3A_0^2 r_B$};

    \end{tikzpicture}
\end{figure}

Let now $B\subseteq X$ be a ball such that $r_B\geq \frac{r_0}{2A_0}$. To begin with, let us estimate
\begin{equation}\label{eq:proofmuck5}
        \left( \fint_B w^{-\frac{1}{p-1}} \ud\mu \right)^{p-1}=\left( \fint_B \bigl(1+\rho(x,x_0)\bigr)^{\frac{\eta p}{p-1}}  \ud\mu(x) \right)^{p-1}\lesssim \bigl(1+r_B\bigr)^{\eta p}.
    \end{equation}
Next, let us focus on the other factor $ \fint_B w  \ud\mu$. To this end, we first observe that $B \subseteq B(x_0, 3A_0^2 r_B)$. Indeed, for $y \in B$ we have
\begin{equation*}
    \rho(y,x_0)\leq A_0\bigl(\rho(y,x_B)+\rho(x_B,x_0)\bigr)\leq A_0(r_B+r_0)\leq A_0(r_B+2A_0r_B)\leq 3A_0^2 r_B.
\end{equation*}
Moreover, by the doubling property of $\mu$, we have that
\[
\mu(B)\approx \mu\bigl(B(x_0,r_B)\bigr)\approx \mu\bigl(B(x_0,3A_0^2r_B)\bigr).
\]
So, 
\begin{equation}\label{eq:proofmuck1}
        \fint_B w  \ud\mu\lesssim \frac{1}{\mu\bigl(B(x_0,3A_0^2r_B)\bigr)}\int_{B(x_0,3A_0^2r_B)} \bigl(1+\rho(x,x_0)\bigr)^{-\eta p}\ud \mu(x)=:I. 
\end{equation}

Next, let us concentrate on $I$. Recall from Lemma~\ref{lem:noatoms} that, since we are assuming the space to have lower dimension $d>0$, then there are no atoms. Thus, in particular, $\mu(\{x_0\})=0$. We can therefore write
\begin{equation*}
    \begin{split}
        I&=\fint_{B(x_0,3A_0^2r_B)} \bigl(1+\rho(x,x_0)\bigr)^{-\eta p}\ud \mu(x)\\
        &= \frac{1}{\mu\bigl(B(x_0,3A_0^2r_B)\bigr)}\sum_{k=0}^\infty \int_{B(x_0,2^{-k}3A_0^2r_B)\setminus B(x_0,2^{-k-1}3A_0^2r_B)} \bigl(1+\rho(x,x_0)\bigr)^{-\eta p}\ud \mu(x)\\
        &\lesssim   \frac{1}{\mu\bigl(B(x_0,3A_0^2r_B)\bigr)}\sum_{k=0}^\infty \frac{1}{(1+2^{-k}r_B)^{\eta p}}\mu\bigl(B(x_0,2^{-k}3A_0^2r_B)\bigr).
    \end{split}
\end{equation*}
Next, by lower dimension $d$, we have that for every $k\geq0$
\[
\frac{\mu\bigl(B(x_0,3A_0^2r_B)\bigr)}{\mu\bigl(B(x_0,2^{-k}3A_0^2r_B)\bigr)}\gtrsim \left(\frac{3A_0^2r_B}{2^{-k}3A_0^2r_B}\right)^d=2^{kd}\quad \Leftrightarrow\quad \mu\bigl(B(x_0,2^{-k}3A_0^2r_B\bigr)\lesssim 2^{-kd}\mu\bigl(B(x_0,3A_0^2r_B\bigr),
\]
which yields
\begin{equation}\label{eq:proofmuck4}
    I\lesssim \sum_{k=0}^\infty \frac{2^{-kd}}{(1+2^{-k}r_B)^{\eta p}}.
\end{equation}
Finally, we divide the end of the computation into two subcases.
\begin{itemize}
    \item \textbf{Case \boldmath$r_B\geq1$.}\\
    In this case, we estimate 
    \begin{equation}\label{eq:proofmuck3}
         I\lesssim \sum_{k=0}^\infty \frac{2^{-kd}}{(1+2^{-k}r_B)^{\eta p}}\lesssim r_B^{-\eta p}\sum_{k=0}^\infty 2^{-k(d-\eta p)}\lesssim r_B^{-\eta p},
    \end{equation}
    where in the last inequality we have used the fact that by assumption $d-\eta p>0$. Thus, plugging \eqref{eq:proofmuck3} into \eqref{eq:proofmuck1}, we get
    \begin{equation*}
        \fint_B w  \ud\mu\lesssim r_B^{-\eta p}.
    \end{equation*}
    Hence, combining this estimate with \eqref{eq:proofmuck5}, we get
    \[
    \left( \fint_B w  \ud\mu \right) 
  \left( \fint_B w^{-\frac{1}{p-1}}  \ud\mu \right)^{p-1}\lesssim \left(\frac{1+r_B}{r_B}\right)^{\eta p}\lesssim 1.
    \]
    \textbf{Case \boldmath$r_B<1$.}\\
    In this case, starting again from \eqref{eq:proofmuck4}, we can write
    \begin{equation*}
        I\lesssim \sum_{k=0}^\infty \frac{2^{-kd}}{(1+2^{-k}r_B)^{\eta p}}\leq \sum_{k=0}^\infty 2^{-kd}\lesssim1.
    \end{equation*}
    Thus, we get
    \begin{equation*}
        \fint_B w  \ud\mu\lesssim 1.
    \end{equation*}
   Hence, by \eqref{eq:proofmuck5},
   \[
   \left( \fint_B w  \ud\mu \right) 
  \left( \fint_B w^{-\frac{1}{p-1}}  \ud\mu \right)^{p-1}\lesssim (1+r_B)^{\eta p}\lesssim 1.
   \]
\end{itemize}
\end{itemize}    
The proof is now complete.
\end{proof}

At this point, we are ready to state and prove the boundedness result for fractional integral operators of metric type. 

\begin{theorem}\label{thm:bddfracmet}
Let $(X,\rho,\mu)$ be a space of homogeneous type with lower dimension $d>0$.  
Let $p,q\in(1,\infty)$ and $\eta\in(0,d)$ satisfy the relation $\frac{1}{q}=\frac{1}{p}-\frac{\eta}{d}>0$.  
Then the operator
\[
Sf(x):=\int_X \frac{\rho(x,y)^{\eta}}{V(x,y)}\,f(y)\,\ud\mu(y), \quad x\in X,
\]
is bounded from $L^p(\mu)$ into $L^p(w)+L^q(\mu)$, where $w$ is the Muckenhoupt weight in $A_p$ given by \( w(x):=\bigl(1+\rho(x,x_0)\bigr)^{-\eta p}\) for any fixed $x_0\in X$.
\end{theorem}

\begin{remark}\label{rmk:sumnorm}
The norm on the sum of two normed spaces $X$ and $Y$ contained in a common topological vector space is defined in the standard way as
\[
\|z\|_{X+Y}:=\inf\bigl\{\|x\|_X+\|y\|_Y : z=x+y,\ x\in X,\ y\in Y\bigr\}.
\]
This is the conventional choice in the literature (see, e.g., \cite[Lemma 2.3.1]{BL:76}).
\end{remark}

\begin{remark}
    We point out that the assumption of positive lower dimension is essential in order to obtain a result of this type. Indeed, the example constructed in Appendix~\ref{app:Cantor} shows that, in the absence of this property, a generic metric fractional integral may map functions in $L^p$ to the function identically equal to $+\infty$.
\end{remark}

\begin{proof}[Proof of Theorem~\ref{thm:bddfracmet}]
    Let us fix $x_0\in X$. Then,
    \begin{equation}\label{eq:proofbddmetric1}
        \begin{split}
            Sf(x)&=\int_X \frac{\rho(x,y)^{\eta}}{V(x,y)}\,f(y)\,\ud\mu(y)\\
            &=\int_{\{\rho(x,y)\leq \rho(x,x_0)+1\}} \frac{\rho(x,y)^{\eta}}{V(x,y)}\,f(y)\,\ud\mu(y)+\int_{\{\rho(x,y)> \rho(x,x_0)+1\}} \frac{\rho(x,y)^{\eta}}{V(x,y)}\,f(y)\,\ud\mu(y)\\
            &=:I_1f(x)+I_2f(x).
        \end{split}
    \end{equation}

    Now, if we denote by $r_0:=\rho(x,x_0)+1$ and recall that $\mu(\{x\})=0$ by Lemma~\ref{lem:noatoms}, we get that
    \begin{equation*}
        \begin{split}
            I_1f(x)&=\int_{\{\rho(x,y)\leq r_0\}} \frac{\rho(x,y)^{\eta}}{V(x,y)}\,f(y)\,\ud\mu(y)\\
            &=\sum_{k=0}^\infty\int_{\{2^{-k-1}r_0< \rho(x,y)\leq 2^{-k}r_0\}} \frac{\rho(x,y)^{\eta}}{V(x,y)}\,f(y)\,\ud\mu(y)\\
            &\approx \sum_{k=0}^\infty\int_{\{2^{-k-1}r_0< \rho(x,y)\leq 2^{-k}r_0\}} \frac{(2^{-k}r_0)^{\eta}}{V(x,2^{-k}r_0)}\,f(y)\,\ud\mu(y)\\
            &\leq\sum_{k=0}^\infty(2^{-k}r_0)^{\eta}Mf(x)\\
            &\lesssim r_0^\eta Mf(x),
        \end{split}
    \end{equation*}
    where we denote by $M$ the Hardy--Littlewood maximal operator. Thus, we have shown that
    \begin{equation*}
        I_1f(x)\lesssim \bigl(\rho(x,x_0)+1\bigr)^\eta Mf(x). 
    \end{equation*}
    By the $L^p$-boundedness of the Hardy--Littlewood maximal operator, and if we denote by \( w(x) := \bigl(1+\rho(x,x_0)\bigr)^{-\eta p} \), we get that
    \begin{equation}\label{eq:proofbddmetric2}
        \Norm{I_1f}{L^p(w)}\lesssim \BNorm{x\mapsto \bigl(\rho(x,x_0)+1\bigr)^\eta Mf(x)}{L^p(w)}=\Norm{Mf}{L^p(\mu)}\lesssim \Norm{f}{L^p(\mu)}.
    \end{equation}  
Recall that the weight $w$ belongs to $A_p$ by Lemma~\ref{lem:Muck}.  

Let us now turn our attention to $I_2f$.  
By the assumption of lower dimension $d$ and by the doubling property of $\mu$, we have that, if $\rho(x,y)>\rho(x,x_0)+1$,
\begin{equation*}
        V(x,y)=V(x,\rho(x,y))\approx V(x_0,\rho(x,y))\gtrsim \rho(x,y)^dV(x_0,1)\approx \rho(x,y)^d.
\end{equation*}
Therefore, 
\[
\rho(x,y)^\eta\lesssim V(x,y)^{\frac{\eta}{d}},
\]
which implies that
\begin{equation*}
    \begin{split}
        I_2f(x)&=\int_{\{\rho(x,y)> \rho(x,x_0)+1\}} \frac{\rho(x,y)^{\eta}}{V(x,y)}\,f(y)\,\ud\mu(y)\\
        &\lesssim\int_{\{\rho(x,y)> \rho(x,x_0)+1\}} \frac{V(x,y)^{\frac{\eta}{d}}}{V(x,y)}\,f(y)\,\ud\mu(y)\\
        &\leq\int_X \frac{V(x,y)^{\frac{\eta}{d}}}{V(x,y)}\,f(y)\,\ud\mu(y).
    \end{split}
\end{equation*}
By the boundedness of the volumic-type fractional integrals stated in Theorem~\ref{thm:bddfracvol}, we may conclude that for $q\in(1,\infty)$ such that $\frac1q=\frac1p-\frac{\eta}{d}$
\begin{equation}\label{eq:proofbddmetric3}
    \Norm{I_2f}{L^q(\mu)}\lesssim \Norm{f}{L^p(\mu)}.
\end{equation}
At this point, it only remains to combine \eqref{eq:proofbddmetric1}, \eqref{eq:proofbddmetric2}, and \eqref{eq:proofbddmetric3}, from which we can deduce that
\[
\Norm{Sf}{L^p(w)+L^q(\mu)}\leq \Norm{I_1f}{L^p(w)}+\Norm{I_2f}{L^q(\mu)}\lesssim \Norm{f}{L^p(\mu)}.
\]
This concludes the proof.
\end{proof}

\subsection{Estimates for Haar coefficients}
In the study of the dyadic formula, which will be treated in Section~\ref{sec:DRT}, we will need estimates for the Haar coefficients associated with fractional integrals. To this end, the present subsection is devoted to establishing such estimates. The procedure presented here closely follows, with the necessary adjustments due to the different setting, the one given by T.~Hytönen in \cite[Section~5]{Hyt:Enflo}. 

We clarify that, from now on, we will use the lighter notation $\mathscr D$ and $Q \in \mathscr D$ to refer to the random dyadic systems and cubes, while keeping in mind that they always depend on the choice of the outcome in the sample space \(\Omega\). In particular, it is important to recall that $\{z_Q\}$ denotes random reference points, whereas $\{x_Q\}$ denotes the deterministic reference points, that is, those which do not vary with the different outcomes in \(\Omega\). 

With this necessary premise in place, we begin by stating a result by P. Auscher and E. Routin \cite{AR:13}, which follows as a corollary of their \cite[Lemma 2.4]{AR:13}.

\begin{lemma}\label{lem:AusRou}
    Let $(X,\rho,\mu)$ be a space of homogeneous type and let $P,Q\in \mathscr D$ be disjoint. Then
    \begin{equation*}
        \int_P\int_Q\frac{1}{V(x,y)}\ud\mu(x)\ud\mu(y)\lesssim \mu(P)^{1/s}\mu(Q)^{1/s'},
    \end{equation*}
    for any choice of $s\in(1,\infty)$ with dual exponent $s'$.
\end{lemma}
\begin{proof}
    The statement follows immediately from \cite[Lemma 2.4]{AR:13} by taking $f:=\mathbf{1}_Q$ and $g:=\mathbf{1}_P$.
\end{proof}

Next, let us prove a simple estimate that will be useful later on. 

\begin{lemma}\label{lem:auxprelest12}
     Let $(X,\rho,\mu)$ be a space of homogeneous type, and let $P,Q\in \mathscr D_i$ for some $i\in\Z$, i.e., let $P,Q$ be dyadic cubes belonging to the same generation.
     Then, for any possible choice of $\phi$ in Notation \ref{not:phi}, we have that
     \begin{equation*}
         \phi(x,y)\lesssim \phi\bigl(z_{P},\rho(z_{P},z_{Q})+\ell(P)\bigr).
     \end{equation*}
     for every $x\in P$ and $y\in Q$.
\end{lemma}

\begin{proof}
    Let us consider the two possible instances of $\phi$ separately.
\begin{itemize}[label={}, left=-4.5pt]
\item \textbf{Metric case.}

In this case, we consider $\phi(x,y)=\rho(x,y)^\eta$ and $\phi(x,r)=r^\eta$ for some $\eta\in(0,\infty)$.

First, recall that dyadic cubes $P$ and $Q$ from the same generation, in particular, have equal side lengths, that is, $\ell(P) = \ell(Q)$. Thus, for all $(x,y)\in P\times Q$, we have
\begin{equation}\label{eq:proofauxil1}
    \rho(x,y)\lesssim \ell(P)+\rho(z_{P},z_{Q})+\ell(Q)\approx \rho(z_{P},z_{Q})+\ell(P)
\end{equation}
The claim now follows by raising both sides of this chain of inequalities to the power $\eta$, and then identifying the resulting expression with the appropriate notation involving $\phi$.

\item \textbf{Volumic case.}

In this case, we take $\phi(x,y)=V(x,y)^\eps$ and $\phi(x,r)=V(x,r)^\eps$ for some $\eps\in(0,1)$. 

By \eqref{eq:proofauxil1} and by the doubling property of $\mu$, we get that
\begin{equation*}
    V(x,y)\lesssim V(x,\rho(z_{P},z_{Q})+\ell(P))\lesssim V(z_P,\rho(z_{P},z_{Q})+\ell(P)).
\end{equation*}
As before, the claim follows by raising both sides to the power $\varepsilon$ and making the appropriate identifications using the notation $\phi$. \qedhere
\end{itemize}
\end{proof}

We now turn to the study of how the quadratic form associated with a fractional integral operator acts on the indicator functions of dyadic cubes. It is worth noting that the estimates that follow do not require, in the metric case, the lower dimension $d>0$ of the underlying space; in other words, in the metric case, it is not necessary to rely on the boundedness described in Theorem \ref{thm:bddfracmet}.

\begin{lemma}\label{lem:QuadFracIntCubes}
    Let $(X,\rho,\mu)$ be a space of homogeneous type, and let $T$ be a $\phi$-fractional integral operator as in Definition~\ref{def:fracint}. Then, for every dyadic cube $Q\in \mathscr{D}$, we have that
    \begin{equation*}
        \abs{\langle T\mathbf{1}_Q,\mathbf{1}_Q\rangle}\lesssim \phi(Q)\mu(Q),
    \end{equation*}
    where $\phi(Q)$ is defined as in Notation~\ref{not:phi}.
\end{lemma}

\begin{proof}
    Here too, it is convenient to divide the proof into the two possible cases.
    \begin{itemize}[label={}, left=-4.5pt]
\item \textbf{Volumic case.}

In this case, from Theorem~\ref{thm:bddfracvol}, we know that $T$ is bounded from \( L^p(X,\mu) \) to \( L^q(X,\mu) \) for all $p,q\in(1,\infty)$ such that $1/q=1/p-\eps>0$. Hence, for such $p,q$, we have that
    \begin{align*}
        |\langle T\mathbf{1}_{Q},\mathbf{1}_{Q}\rangle|&\leq\|T\mathbf{1}_Q\|_{L^q(\mu)}\|\mathbf{1}_Q\|_{L^{q'}(\mu)}\\
        &\lesssim\|\mathbf{1}_Q\|_{L^{p}(\mu)}\|\mathbf{1}_Q\|_{L^{q'}(\mu)}=\mu(Q)^{1/p+1/q'} =\mu(Q)^{1+\eps}=\phi(Q)\mu(Q).
    \end{align*}
\item \textbf{Metric case.}

In this case, we begin by observing that for any cube $Q\in\mathscr D$, $x\in Q$, and for some constants $0<C_0<C_1<\infty$ depending only on the structure of the underlying space $X$, we have 
    \begin{equation*}
        \begin{split}
            \abs{T\mathbf{1}_Q(x)}&\lesssim \int_X \frac{\rho(x,y)^\eta}{V(x,y)}\mathbf{1}_Q(y)\ud\mu(y) \\
            &\lesssim \int_{B(z_Q,C_0\ell(Q))} \frac{\rho(x,y)^\eta}{V(x,y)}\ud\mu(y)\\
            &\lesssim \int_{B(x,C_1\ell(Q))} \frac{\rho(x,y)^\eta}{V(x,y)}\ud\mu(y)\\
            &=\sum_{k=0}^\infty \int_{B(x,r_k)\setminus B(x,r_{k+1})} \frac{\rho(x,y)^\eta}{V(x,y)}\ud\mu(y), \qquad r_k := 2^{-k}C_1\ell(Q)\\
            &\lesssim \sum_{k=0}^\infty \Big(2^{-k}C_1\ell(Q)\Big)^\eta \approx \ell(Q)^\eta.
        \end{split}
    \end{equation*}
Thus, 
\begin{equation*}
    \abs{\pair{T\mathbf{1}_Q}{\mathbf{1}_Q}}\leq \int_Q\abs{T\mathbf{1}_Q(x)}\ud\mu(x)\lesssim \ell(Q)^\eta \mu(Q)=\phi(Q)\mu(Q).
\end{equation*}
\qedhere
\end{itemize}
\end{proof}

Using Lemmas~\ref{lem:AusRou}, \ref{lem:auxprelest12} and \ref{lem:QuadFracIntCubes}, we now prove two estimates for the Haar coefficients associated with fractional integrals, which constitute the main goal of this subsection.

\begin{lemma}\label{lem:prelest1}
    Let $(X,\rho,\mu)$ be a space of homogeneous type, and let $T$ be a $\phi$-fractional integral operator as in Definition~\ref{def:fracint}. Furthermore, let $P,Q\in \mathscr D_i$ for some $i\in\Z$, i.e., let $P,Q$ be dyadic cubes belonging to the same generation. Then, for every $0\leq\alpha\leq M_P$ and $0\leq\beta\leq M_Q$ the following estimate holds true:
    \begin{align}
        &|\langle Th_P^\alpha,h_Q^\beta \rangle |\lesssim \phi(z_P, \rho(z_P,z_Q)+\ell(P)) \label{eq:cubes2}.
    \end{align}
\end{lemma}

\begin{remark}
This result remains valid when deterministic reference points are used in place of random ones; in particular, under these assumptions we obtain
\begin{equation}\label{eq:cubes1}
        |\langle Th_P^\alpha,h_Q^\beta \rangle |\lesssim \phi(x_P, \rho(x_P,x_Q)+\ell(P)).
    \end{equation}
In fact, this expression and \eqref{eq:cubes2} are fully equivalent, as their right-hand sides are essentially comparable, for both possible definitions of $\phi$.
\end{remark}

\begin{proof}[Proof of Lemma \ref{lem:prelest1}.]
    To begin with, recall that dyadic cubes $P$ and $Q$ belonging to the same generation in particular satisfy $\ell(P) = \ell(Q)$.

    Next, let us write
    \begin{equation}\label{eq:prelest11}
        \langle Th_P^\alpha,h_Q^\beta\rangle= \sum_{\substack{ P'\in\operatorname{ch}(P) \\ Q'\in\operatorname{ch}(Q)}} \langle T\mathbf{1}_{P'},\mathbf{1}_{Q'}\rangle \langle h_P^\alpha\rangle_{P'} \langle h_Q^\beta\rangle_{Q'}.
    \end{equation}
    Notice that $\langle h_P^\alpha\rangle_{P'}$ and $\langle h_Q^\beta\rangle_{Q'}$ in this case are simply the constant values that the functions $h_P^\alpha$ and $h_Q^\beta$ take on each of the children of the cubes $P$ and $Q$, respectively.

    If $P'\neq Q'$, then, by Lemma~\ref{lem:auxprelest12}, for all $(x,y)\in P'\times Q'\subset P\times Q$ we have that
    \begin{equation*}
        \phi(x,y)\lesssim \phi\bigl(z_{P},\rho(z_{P},z_{Q})+\ell(P)\bigr).
    \end{equation*}
    Thus, by this estimate and by Lemma~\ref{lem:AusRou} we get that
    \begin{align}
        |\langle T\mathbf{1}_{P'},\mathbf{1}_{Q'}\rangle|&=\left|\int_{P'}\int_{Q'} K(x,y)\ud\mu(x)\ud\mu(y)\right|\nonumber \\
        &\lesssim \int_{P'}\int_{Q'} \frac{\phi(x,y)}{V(x,y)}\ud\mu(x)\ud\mu(y) \nonumber \\
        &\lesssim \phi\bigl(z_{P},\rho(z_{P},z_{Q})+\ell(P)\bigr)\mu(P')^{1/2}\mu(Q')^{1/2},\nonumber
    \end{align}
    where in the second step we have used the $\phi$-size estimate~\eqref{eq:CZ0phi} from Definition~\ref{def:fracint}. Thus,
    \begin{align*}
        \Big|&\sum_{\substack{ P'\in\operatorname{ch}(P) \\ Q'\in\operatorname{ch}(Q)\\ P'\neq Q'}} \langle T\mathbf{1}_{P'},\mathbf{1}_{Q'}\rangle \langle h_P^\alpha\rangle_{P'} \langle h_Q^\beta\rangle_{Q'}\Big| \\
        &\lesssim \phi\bigl(z_{P},\rho(z_{P},z_{Q})+\ell(P)\bigr) \left(\sum_{P'} \mu(P')^{1/2}|\langle h_P^\alpha\rangle_{P'}|\right)\left(\sum_{Q'} \mu(Q')^{1/2}|\langle h_Q^\beta\rangle_{Q'}|\right),
    \end{align*}
    and using Cauchy-Schwarz:
    \begin{align*}
        \sum_{P'} \mu(P')^{1/2}|\langle h_P^\alpha\rangle_{P'}| &\leq M^{1/2} \left(\sum_{P'} \mu(P')|\langle h_P^\alpha\rangle_{P'}|^2\right)^{1/2} \\
        &= M^{1/2} \|h_P^\alpha\|_{L^2(\mu)} = M^{1/2},
    \end{align*}
    with a similar bound holding for the sum over $Q'$. Therefore, we get that
    \begin{equation*}
        \Big|\sum_{\substack{ P'\in\operatorname{ch}(P) \\ Q'\in\operatorname{ch}(Q)\\ P'\neq Q'}} \langle T\mathbf{1}_{P'},\mathbf{1}_{Q'}\rangle \langle h_P^\alpha\rangle_{P'} \langle h_Q^\beta\rangle_{Q'}\Big|
        \lesssim \phi\bigl(z_{P},\rho(z_{P},z_{Q})+\ell(P)\bigr).
    \end{equation*}
    
    If $P\neq Q$, then all their children $P'$ and $Q'$ are different, and so, keeping in mind formula \eqref{eq:prelest11}, we may conclude that
    \begin{equation*}
        |\langle Th_P^\alpha,h_Q^\beta \rangle |\lesssim \phi\bigl(z_{P},\rho(z_{P},z_{Q})+\ell(P)\bigr),
    \end{equation*}
    which is what we aimed to prove.

    If instead $P = Q$, we also need to consider the sum 
\begin{equation*}
        \Big|\sum_{P'\in\operatorname{ch}(P)} \langle T\mathbf{1}_{P'},\mathbf{1}_{P'}\rangle \langle h_P^\alpha\rangle_{P'} \langle h_P^\beta\rangle_{P'}\Big|.
\end{equation*}        
    To handle this sum, whose estimate concludes the proof, we use Lemma~\ref{lem:QuadFracIntCubes} and the monotonicity of $\phi$ given by \eqref{eq:RemarkNotation5} to get
    \begin{align*}
         \Big|\sum_{P'\in\operatorname{ch}(P)} &\langle T\mathbf{1}_{P'},\mathbf{1}_{P'}\rangle \langle h_P^\alpha\rangle_{P'} \langle h_P^\beta\rangle_{P'}\Big| \\
         &\lesssim \sum_{P'\in\operatorname{ch}(P)} \phi(P')\mu(P') |\langle h_P^\alpha\rangle_{P'} \langle h_P^\beta\rangle_{P'}| \\
        &\leq \phi(P) \sum_{P'\in\operatorname{ch}(P)} \mu(P') |\langle h_P^\alpha\rangle_{P'} \langle h_P^\beta\rangle_{P'}| \\
        &\leq \phi(P) \left(\sum_{P'\in\operatorname{ch}(P)} \mu(P') |\langle h_P^\alpha\rangle_{P'}|^2\right)^{1/2} \left(\sum_{P'\in\operatorname{ch}(P)} \mu(P') |\langle h_P^\beta\rangle_{P'}|^2\right)^{1/2} \\
        &= \phi(P) \|h_P^\alpha\|_{L^2(\mu)} \|h_P^\beta\|_{L^2(\mu)} \\
        &= \phi(P)\\
        &\approx \phi\big(z_P,\ell(P)\big)\\
        &=\phi\big(z_P,\rho(z_P,z_P)+\ell(P)\big),
    \end{align*}
where in the penultimate step we have used property \eqref{eq:RemarkNotation4}. This concludes the proof.
\end{proof}

\begin{lemma}\label{lem:prelest2}
    Let $(X,\rho,\mu)$ be a space of homogeneous type, and let $T$ be an $\omega$-regular $\phi$-fractional integral operator as in Definition~\ref{def:fracint}, for some function $\omega : [0,1] \to [0,\infty)$ that satisfies Assumption~\ref{ass:omega}. Furthermore, let $P,Q\in \mathscr D$ be such that $\ell(P)=\ell(Q)$, and let $(\alpha,\beta)\neq(0,0)$. Then
    \begin{equation}\label{eq:cubes3}
        |\langle Th_P^\alpha,h_Q^\beta \rangle |\lesssim \omega\left(\frac{\ell(P)}{\rho(z_Q,z_P)+\ell(P)}\right)\frac{\sqrt{\mu(P)\mu(Q)}}{V\bigl(z_Q,\rho(z_Q,z_P)+\ell(P)\bigr)}\phi\bigl(z_Q,\rho(z_Q,z_P)+\ell(P)\bigr).
    \end{equation}
    \end{lemma}
    \begin{remark}
        As in Lemma \ref{lem:prelest1}, this statement also holds when considering deterministic reference points instead of random ones; that is, under these assumptions we also have
\begin{equation}\label{eq:cubes4}
        |\langle Th_P^\alpha,h_Q^\beta \rangle |\lesssim \omega\left(\frac{\ell(P)}{\rho(x_Q,x_P)+\ell(P)}\right)\frac{\sqrt{\mu(P)\mu(Q)}}{V\bigl(x_Q,\rho(x_Q,x_P)+\ell(P)\bigr)}\phi\bigl(x_Q,\rho(x_Q,x_P)+\ell(P)\bigr).
\end{equation}
Indeed, this formula and \eqref{eq:cubes3} are completely equivalent, since their right-hand sides are essentially comparable.
    \end{remark}

\begin{proof}[Proof of Lemma \ref{lem:prelest2}.]
We divide this proof into two cases: when the cubes $P$ and $Q$ are well separated, and when they are not.
\begin{itemize}[label={}, left=-4.5pt]
\item \textbf{Case \boldmath $\rho(P,Q) \gg \ell(P) = \ell(Q)$.}\\
 In this case, assuming for instance that $\alpha \neq 0$ (the case $\beta \neq 0$ is treated symmetrically), we have that
    \begin{align*}
        |\langle Th_P^\alpha,h_Q^\beta \rangle |=&\left|\int_P\int_QK(x,y)h_P^\alpha(y)h_Q^\beta(x)\ud\mu(x)\ud\mu(y)\right|\\
        =&\left|\int_P\int_Q\bigl(K(x,y)-K(x,z_P)\bigr)h_P^\alpha(y)h_Q^\beta(x)\ud\mu(x)\ud\mu(y)\right|,
    \end{align*}
and also, for $(x,y)\in Q\times P$,
\begin{equation*}
    \rho(y,z_P)\lesssim \ell(P)\ll\rho(x,z_P).
\end{equation*}
Hence, we may apply the $\omega$,$\phi$-smoothness estimate~\eqref{eq:CZ1phi} of Definition~\ref{def:fracint} and we get
\begin{equation}
    |\langle Th_P^\alpha,h_Q^\beta \rangle |\lesssim \int_P\int_Q \omega\Big(\frac{\rho(y,z_P)}{\rho(x,z_P)}\Big)\frac{\phi(x,z_P)}{V(x,z_P)}|h_P^\alpha(y)||h_Q^\beta(x)|\ud\mu(x)\ud\mu(y)\label{eq:prelest21}. 
\end{equation}
Now, notice that for $x\in Q$ there holds
\begin{equation*}
    \rho(x,z_P)\lesssim \rho(x,z_Q)+\rho(z_Q,z_P)\lesssim \ell(Q)+\rho(z_Q,z_P)
    = \ell(P)+\rho(z_Q,z_P)\lesssim \rho(x,z_P),
\end{equation*}
which yields
\begin{equation*}
    \rho(x,z_P)\approx \ell(P)+\rho(z_Q,z_P).
\end{equation*}
Hence, applying the properties \eqref{eq:RemarkNotation1}, \eqref{eq:RemarkNotation3}, and \eqref{eq:RemarkNotation2} of the function $\phi$, in this order, we obtain
\begin{equation*}
   \phi(x,z_P)\approx \phi(z_P,x) \approx \phi\bigl(z_P, \ell(P)+\rho(z_Q,z_P)\bigr)\approx \phi\bigl(z_Q, \ell(P)+\rho(z_Q,z_P)\bigr),
\end{equation*}
and, by the doubling property of $\omega$,
\begin{equation*}
    \omega\Big(\frac{\rho(y,z_P)}{\rho(x,z_P)}\Big)\lesssim\omega\Big(\frac{\ell(P)}{\rho(x,z_P)}\Big)\approx \omega\Big(\frac{\ell(P)}{\ell(P)+\rho(z_Q,z_P)}\Big).
\end{equation*}
Therefore, we get that
\begin{align*}
    \text{\eqref{eq:prelest21}}&\lesssim \omega\Big(\frac{\ell(P)}{\ell(P)+\rho(z_Q,z_P)}\Big) \frac{\|h_P^\alpha\|_{L^1(\mu)}\|h_Q^\beta\|_{L^1(\mu)}}{V\bigl(z_Q, \ell(P)+\rho(z_Q,z_P)\bigr)}\phi\bigl(z_Q, \ell(P)+\rho(z_Q,z_P)\bigr) \\
    &\lesssim \omega\Big(\frac{\ell(P)}{\ell(P)+\rho(z_Q,z_P)}\Big) \frac{\sqrt{\mu(P)\mu(Q)}}{V\bigl(z_Q, \ell(P)+\rho(z_Q,z_P)\bigr)}\phi\bigl(z_Q, \ell(P)+\rho(z_Q,z_P)\bigr),
\end{align*}
which concludes the proof for this case.

\item \textbf{Case \boldmath $\rho(P,Q)\lesssim \ell(P)=\ell(Q)$.}\\ 
To address this case, we first recall that, by Lemma~\ref{lem:prelest1}, Formula~\eqref{eq:cubes2}, we have that
\begin{equation}\label{eq:prelest22}
    |\langle Th_P^\alpha,h_Q^\beta \rangle |\lesssim \phi(z_P, \rho(z_P,z_Q)+\ell(P))\approx \phi(z_Q, \rho(z_P,z_Q)+\ell(P)).
\end{equation}
On the other hand, as
\begin{equation*}
    \rho(z_P,z_Q)+\ell(P)\approx \ell(P)=\ell(Q),
\end{equation*}
we get that
\begin{equation*}
    \omega\Big(\frac{\ell(P)}{\ell(P)+\rho(z_Q,z_P)}\Big)\approx\omega(1)\approx1,
\end{equation*}
and
\begin{equation*}
    V\bigl(z_Q, \ell(P)+\rho(z_Q,z_P)\bigr)\approx \mu(Q)\approx\mu(P).
\end{equation*}
Thus,
\begin{equation*}
    \omega\Big(\frac{\ell(P)}{\ell(P)+\rho(z_Q,z_P)}\Big)\frac{\sqrt{\mu(P)\mu(Q)}}{V\bigl(z_Q, \ell(P)+\rho(z_Q,z_P)\bigr)}\approx1,
\end{equation*}
which, combined with \eqref{eq:prelest22}, concludes the proof. \qedhere
\end{itemize}
\end{proof}

\section{Fractional dyadic model operators}\label{sec:DyadModOp}

In this section, we introduce and investigate several properties of the dyadic model operators that serve as the building blocks of the representation formula for fractional integrals presented in Theorem~\ref{thm:DRT}. In particular, we study their properties under a certain fractional normalization, which arises naturally in the study of fractional integral operators $T$. 

\subsection{Dyadic paraproducts and fractional normalization}

Here, we focus on defining the dyadic paraproducts with fractional normalization and relating them to fractional integrals. In particular, we will see that a $\phi$-normalized dyadic paraproduct can be pointwise dominated by a $\phi$-fractional integral, and that it is possible, under certain hypotheses, to make sense of the notion of paraproduct with symbol $T1$, where $T$ is a fractional integral.

\subsubsection*{Definition and pointwise domination by fractional integrals}

Let us begin by defining the dyadic paraproducts and the notion of $\phi$-normalization.


\begin{definition}\label{def:paraprod}
Let $(X,\rho,\mu)$ be a space of homogeneous type, and let
$\{b_Q\}_{Q\in\mathscr D}$ be a family of functions satisfying
\(
b_Q\in L^2_0(\operatorname{ch}(Q)),
\)
(see \eqref{eq:L20children}) or equivalently,
\[
\D_Qb_Q=b_Q,\qquad Q\in\mathscr D.
\]
The \emph{dyadic paraproduct associated with the family}
$\{b_Q\}_{Q\in\mathscr D}$ is the operator formally defined by
\begin{equation}\label{eq:DefParaprod}
\Pi_{\{b_Q\}}
:=
\sum_{Q\in\mathscr D}
b_Q\otimes\frac{\mathbf1_Q}{\mu(Q)}
=
\sum_{k\in\mathbb Z}
\sum_{Q\in\mathscr D_k}
b_Q\otimes\frac{\mathbf1_Q}{\mu(Q)}.
\end{equation}
We say that $\Pi_{\{b_Q\}}$ is \emph{$\phi$-normalized} if
\begin{equation}\label{eq:DefParaprodNorm}
\|\{b_Q\}\|_\phi
:=\sup_{Q\in\mathscr D}
\frac{\|b_Q\|_{L^\infty(\mu)}}{\phi(Q)} =
\sup_{\substack{k\in\Z\\ Q\in\mathscr D_k}}
\frac{\|b_Q\|_{L^\infty(\mu)}}{\phi(Q)} 
<\infty.
\end{equation}
\end{definition}

\begin{remark}
    Notice that, in \eqref{eq:DefParaprod} and \eqref{eq:DefParaprodNorm}, every cube $Q\in\mathscr D$ is counted exactly once in the first expression, while it may appear multiple times in the double sum, according to the number of dyadic levels $\mathscr D_k$ containing it. Nevertheless, the two expressions agree, since $L^2_0(\children(Q))$, and consequently the operator $\D_Q=\D_{(Q,k)}$, is non-trivial only for $k=k_{\max}(Q):=\sup\{h\in\mathbb Z:Q\in\mathscr D_h\}\in\mathbb Z\cup\{\infty\}$.
\end{remark}

We proceed by proving that the $\phi$-normalized dyadic paraproducts admit pointwise domination by fractional integral operators.

\begin{proposition}\label{prop:DomPhiNormParaprod}
    Let $(X,\rho,\mu)$ be a space of homogeneous type, and let $\Pi_{\{b_Q\}}$ be a $\phi$-normalized dyadic paraproduct. In the metric setting, assume additionally that $(X,\rho,\mu)$ has positive lower dimension $d>0$ and that the metric fractional parameter satisfies $\eta\in(0,d)$. Then, for $p\in(1,\infty)$ and $f\in L^p(X)$, the series defining
$\Pi_{\{b_Q\}}f$ converges absolutely for $\mu$-almost every
$x\in X$, and
    \begin{equation}\label{eq:DomPhiNormParaprod}
        \abs{\Pi_{\{b_Q\}} f(x)}\leq  \sum_{Q\in\mathscr D}\abs{b_Q(x) \ave{f}_Q}\lesssim \int_X \frac{\phi(x,y)}{V(x,y)} \abs{f(y)}\ud \mu(y), \quad \text{for $\mu$-a.e. }x\in X.
    \end{equation}
\end{proposition}

\begin{proof}
First, by $\phi$-normalization, we have in particular that for all $Q\in\mathscr D$
\[
\|b_Q\|_{L^\infty(\mu)}\lesssim \phi(Q).
\]
Thus, for $\mu$-almost every $x\in X$, we have
\begin{equation*}
\begin{split}
   \sum_{Q\in\mathscr D}\abs{b_Q(x) \ave{f}_Q}&=\sum_{Q\in\mathscr D}\abs{b_Q(x) \ave{f}_Q\mathbf{1}_Q(x)}\\&\lesssim \sum_{\substack{Q\in\mathscr D\\Q\ni x}}\frac{\phi(Q)}{\mu(Q)}\int_Q\abs{f(y)}\ud\mu(y)\\&=\int_X \Big(\sum_{\substack{Q\in\mathscr D\\Q\ni x,y}}\frac{\phi(Q)}{\mu(Q)}\Big)\abs{f(y)}\ud\mu(y) \\
    &=\int_{X\setminus\{x\}} \Big(\sum_{\substack{Q\in \mathscr D\\ Q\ni x,y}} \frac{\phi(Q)}{\mu(Q)}\Big)\abs{f(y)}\ud\mu(y) + \int_{\{x\}} \Big(\sum_{\substack{Q\in \mathscr D\\ Q\ni x}} \frac{\phi(Q)}{\mu(Q)}\Big)\abs{f(y)}\ud\mu(y).
\end{split}
\end{equation*}
Now, by Lemma~\ref{lem:sumcubesfrac}, for $y\neq x$ we have
\[
\sum_{\substack{Q\in \mathscr D\\ Q\ni x,y}} \frac{\phi(Q)}{\mu(Q)}
\lesssim \frac{\phi(x,y)}{V(x,y)}.
\]
Furthermore, in the metric case, Lemma~\ref{lem:noatoms} ensures that the assumption of positive lower dimension $d>0$ yields
\[
\mu(\{x\})=0.
\]
In the volumic case, instead, if $y=x$ and $\mu(\{x\})>0$, another application of Lemma~\ref{lem:sumcubesfrac} gives
\[
\sum_{\substack{Q\in\mathscr D\\Q\ni x}} \frac{\phi(Q)}{\mu(Q)}
\lesssim \frac{1}{\mu(\{x\})^{1-\eps}}
=\frac{\phi(x,y)}{V(x,y)}.
\]
Combining these observations, we obtain that for $\mu$-a.e. $x\in X$ we have
\[
 \sum_{Q\in\mathscr D}\abs{b_Q(x) \ave{f}_Q}\lesssim\int_X \frac{\phi(x,y)}{V(x,y)}\abs{f(y)}\ud\mu(y).
\]
Finally, recall that we formally defined
\[
\Pi_{\{b_Q\}}f:= \Big(\sum_{Q\in\mathscr D}
b_Q\otimes\frac{\mathbf1_Q}{\mu(Q)}\Big)(f)=\sum_{Q\in\mathscr D}b_Q \ave{f}_Q.
\]
Thus, by triangle inequality and Theorems~\ref{thm:bddfracvol} and \ref{thm:bddfracmet}, the series defining $\Pi_{\{b_Q\}}f$ is absolutely convergent $\mu$-a.e., and \eqref{eq:DomPhiNormParaprod} holds. This completes the proof.
\end{proof}

As an immediate consequence, every $\phi$-normalized dyadic paraproduct enjoys the same $L^p$-boundedness properties as the corresponding $\phi$-fractional integral operators.

\begin{corollary}
    Let $(X,\rho,\mu)$ be a space of homogeneous type. Then:

    \begin{itemize}
        \item If $\Pi_{\{b_Q\}}$ is a volumic $\phi$-normalized dyadic shift associated with $\eps\in (0,1)$, then it defines a bounded operator from \( L^p(X,\mu) \) to \( L^q(X,\mu) \) for every $p,q\in(1,\infty)$ satisfying $\frac1q=\frac1p-\eps>0$.

        \item If $(X,\rho,\mu)$ has lower dimension $d>0$ and $\Pi_{\{b_Q\}}$ is a metric $\phi$-normalized dyadic shift associated with $\eta\in(0,d)$, then it defines a bounded operator from $L^p(\mu)$ to $L^p(w)+L^q(\mu)$ for every $p,q\in(1,\infty)$ satisfying
        \(
            \frac{1}{q}=\frac{1}{p}-\frac{\eta}{d}>0.
        \)
        Here, $w$ denotes the Muckenhoupt $A_p$ weight defined by
        \(
            w(x):=\bigl(1+\rho(x,x_0)\bigr)^{-\eta p},
        \)
        where $x_0\in X$ is an arbitrary fixed point.
    \end{itemize}
\end{corollary}

\begin{proof}
    This follows immediately from Proposition~\ref{prop:DomPhiNormParaprod}, together with Theorems~\ref{thm:bddfracvol} and~\ref{thm:bddfracmet}.
\end{proof}

\subsubsection*{On the definition of \boldmath$\D_P T1$}

As introduced in Definition~\ref{def:fracint}, fractional integral operators are, in general, not well-defined when applied to functions that are constant on the whole space $X$. However, we would like to make sense of a \textit{paraproduct with symbol $T1$}. Consequently, we first rigorously define $\D_P T1$ for $P \in \mathscr{D}$ under certain appropriate hypotheses.

\begin{lemma}\label{lem:basicDecAux}
Let $(X,\rho,\mu)$ be a space of homogeneous type, and let $T$ be an $\omega$-regular $\phi$-fractional integral operator as defined in Definition~\ref{def:fracint}, with $\omega : [0,1] \to [0,\infty)$ satisfying Assumption~\ref{ass:omega}, and let $s\in(1,\infty]$. In the metric case, further suppose that $(X,\rho,\mu)$ has positive lower dimension $d>0$ and that the metric fractional parameter belongs to $\eta \in (0,d)$. Then for every $k\in\Z$ and $P\in\mathscr D_k$,
    \begin{equation}\label{eq:basicDecIneq}
        \sum_{Q\in\mathscr D_k}\frac{\|\D_P T\mathbf{1}_Q\|_{L^s(\mu)}}{\mu(P)^{\frac1s}}\lesssim \sum_{m=0}^\infty \omega(\delta^m)\phi\bigl(z_P,\delta^{-m}\ell(P)\bigr),
    \end{equation} 
    with the convention that $\frac1\infty=0$. Consequently, whenever the right-hand side of this inequality
is finite, the series
$\sum_{Q\in\mathscr D_k}\D_P T\mathbf{1}_Q$
converges absolutely in $L^s(\mu)$. In this case, its limit
$$
\D_P T1:=\sum_{Q\in\mathscr D_k} \D_P T\mathbf{1}_Q
$$
satisfies
\begin{equation}\label{eq:basicDecIneq2}
    \Norm{\D_P T1}{L^\infty(\mu)}\lesssim \sum_{m=0}^\infty \omega(\delta^m)\phi\bigl(z_P,\delta^{-m}\ell(P)\bigr)
\end{equation}
\end{lemma}

\begin{remark}
    It is useful to observe that we are dealing with a definition of the expression $\D_P T1$ as an entity, since $T1$ need not be well-defined by itself, and consequently, we cannot interpret $\D_P T1$ as the application of $\D_P$ to $T1$ in the usual sense.
\end{remark}
 \begin{proof}[Proof of Lemma~\ref{lem:basicDecAux}]
  First, using the expansion of $\D_P$ in terms of the Haar functions $h_P^\alpha$, we get that 
  \[
  \D_PT\mathbf{1}_Q= \sum_{1\leq\alpha<M_P}\sqrt{\mu(Q)}\pair{Th_Q^0}{h_P^\alpha}h_P^\alpha
  \]
  Next, by Lemma~\ref{lem:prelest2}, we have that for all $1\leq\alpha\leq M_P$
\begin{equation*}
    |\langle Th_Q^0,h_P^\alpha \rangle |\lesssim \omega\left(\frac{\ell(P)}{\rho(z_Q,z_P)+\ell(P)}\right)\frac{\sqrt{\mu(P)\mu(Q)}}{V\bigl(z_P,\rho(z_Q,z_P)+\ell(P)\bigr)}\phi\bigl(z_P,\rho(z_Q,z_P)+\ell(P)\bigr),
\end{equation*}
and, by \eqref{eq:absHaarBound},
\[
|h_P^\alpha| \lesssim \frac{1}{\sqrt{\mu(P)}}.
\]
Hence,
\begin{equation}\label{eq:basicDecp1}
    \|\D_P T\mathbf{1}_Q\|_{L^\infty(\mu)}\lesssim \omega\left(\frac{\ell(P)}{\rho(z_Q,z_P)+\ell(P)}\right)\frac{\mu(Q)}{V\bigl(z_P,\rho(z_Q,z_P)+\ell(P)\bigr)}\phi\bigl(z_P,\rho(z_Q,z_P)+\ell(P)\bigr).
\end{equation}
Next, let us write
\begin{equation*}
    \begin{split}
        \sum_{Q\in\mathscr D_k}\|\D_P T\mathbf{1}_Q\|_{L^\infty(\mu)}&= \sum_{\substack{Q\in\mathscr D_k \\ \frac{\rho(z_P,z_Q)}{\ell(P)} < 1}}\|\D_P T\mathbf{1}_Q\|_{L^\infty(\mu)}+\sum_{m=1}^\infty\sum_{\substack{Q\in\mathscr D_k \\ \delta \leq \frac{\rho(z_P,z_Q)}{\delta^{-m}\ell(P)} < 1}}
   \|\D_P T\mathbf{1}_Q\|_{L^\infty(\mu)}\\
   &=: I_0+\sum_{m=1}^\infty I_m.
    \end{split}
\end{equation*}
By \eqref{eq:basicDecp1}, we have that for all $m\geq0$
\begin{equation*}
    \begin{split}
        I_m &\lesssim \sum_{\substack{Q\in\mathscr D_k \\ \delta \leq \frac{\rho(z_P,z_Q)}{\delta^{-m}\ell(P)} < 1}} \omega\left(\frac{\ell(P)}{\rho(z_Q,z_P)+\ell(P)}\right)\frac{\mu(Q)}{V\bigl(z_P,\rho(z_Q,z_P)+\ell(P)\bigr)}\phi\bigl(z_P,\rho(z_Q,z_P)+\ell(P)\bigr) \\
        &\lesssim \sum_{\substack{Q\in\mathscr D_k \\ \delta \leq \frac{\rho(z_P,z_Q)}{\delta^{-m}\ell(P)} < 1}} \omega\left(\delta^m\right)\frac{\mu(Q)}{V\bigl(z_P,\delta^{-m}\ell(P)\bigr)}\phi\bigl(z_P,\delta^{-m}\ell(P)\bigr)\\
        &\lesssim \omega\left(\delta^m\right) \phi\bigl(z_P,\delta^{-m}\ell(P)\bigr),
    \end{split}
\end{equation*}
where the lower bound ``$\delta\leq$'' in the summation condition is omitted when $m=0$. Thus, altogether,
\begin{equation}\label{eq:basicDecp2}
    \sum_{Q\in\mathscr D_k}\|\D_P T\mathbf{1}_Q\|_{L^\infty(\mu)} \lesssim \sum_{m=0}^\infty I_m\lesssim \sum_{m=0}^\infty \omega(\delta^m)\phi\bigl(z_P,\delta^{-m}\ell(P)\bigr). 
\end{equation}
This proves \eqref{eq:basicDecIneq} for $s=\infty$.
Moreover, by the triangle inequality, it yields
\eqref{eq:basicDecIneq2}. Finally, since $\D_P T\mathbf{1}_Q$ is supported in $P$, we also have that, for $s\in(1,\infty)$,
\[
\|\D_P T\mathbf{1}_Q\|_{L^s(\mu)}\lesssim \|\D_P T\mathbf{1}_Q\|_{L^\infty(\mu)}\mu(P)^{\frac1s},
\]
which, combined with \eqref{eq:basicDecp2}, gives \eqref{eq:basicDecIneq} also for the case $s\in(1,\infty)$.
\end{proof}

\subsubsection*{Dyadic paraproducts associated with \boldmath$T1$}

Now that we have all the necessary elements, we give the following definition.

\begin{definition}\label{def:paraT1}
Let $(X,\rho,\mu)$ be a space of homogeneous type, and let $T$ be a fractional integral operator satisfying the assumptions of Lemma~\ref{lem:basicDecAux}. Whenever $\D_Q T1$ is well-defined for every $Q\in\mathscr D$, we define
\[
\Pi_{T1}:=\Pi_{\{\D_QT1\}},
\]
and call it the \emph{dyadic paraproduct with symbol $T1$}.
\end{definition}

\begin{remark}
    In this identification, it is useful to keep in mind that the family $\{\D_Q T\mathbf{1}_P\}_{P\in\mathscr D_k}$ is contained in $L^2_0(\children(Q))$ (see \eqref{eq:L20children}), which is finite-dimensional and therefore closed, and is moreover contained in $L^s(\mu)$ for every $s\in(1,\infty]$. Consequently, the sum of the corresponding series, which we denote by $\D_Q T1$, also belongs to $L^2_0(\children(Q))$. This justifies the identification with the dyadic paraproduct.
\end{remark}

At this point, we show that, under certain assumptions, the paraproducts $\Pi_{T1}$ are $\phi$-normalized. We also show that these same assumptions yield a bilinear absolute summability result for the operator
$$
\sum_{k\in\mathbb Z} \sum_{P,Q\in\mathscr D_k} \D_P T\mathbf{1}_Q \otimes \frac{\mathbf{1}_P}{\mu(P)}.
$$

\begin{lemma}\label{lem:PhiNormParaprod}
Let $(X,\rho,\mu)$ be a space of homogeneous type, and let $T$ be
an $\omega$-regular $\phi$-fractional integral operator in the sense
of Definition~\ref{def:fracint}, where
$\omega : [0,1] \to [0,\infty)$ satisfies
Assumption~\ref{ass:omega}. In the metric setting, additionally
assume that $(X,\rho,\mu)$ has positive lower dimension $d>0$
and that the metric fractional parameter satisfies $\eta\in(0,d)$.

Suppose that
\begin{equation}\label{eq:omegasum}
\begin{cases}
\displaystyle
\sum_{m=0}^{\infty}
\omega(\delta^m)\delta^{-mD\varepsilon}<\infty,
& \text{\rm (volumic case)},\\[1em]
\displaystyle
\sum_{m=0}^{\infty}
\omega(\delta^m)\delta^{-m\eta}<\infty,
& \text{\rm (metric case)},
\end{cases}
\end{equation}
where $0<D<\infty$ is the upper dimension of
$(X,\rho,\mu)$, and $\varepsilon\in(0,1)$ denotes the
volumic fractional parameter. Then the following hold:
\begin{enumerate}[label=(\roman*)]
    \item\label{point:PhiNormParaprod1} The dyadic paraproduct $\Pi_{T1}$ is well-defined and $\phi$-normalized.

 \item\label{point:PhiNormParaprod2}
For every pair $(f,g)$ in the following spaces, we have
\begin{equation}\label{eq:PhiNormParaprod1}
    \sum_{k\in\mathbb Z}
\sum_{P,Q\in\mathscr D_k}
\left|
\big\langle \D_P T\mathbf{1}_Q,g\big\rangle
\langle f\rangle_P
\right|\lesssim
\int_X\int_X
\frac{\phi(x,y)}{V(x,y)}
|f(y)||g(x)|\,d\mu(y)\,d\mu(x)
<\infty.
\end{equation}
Here, the finiteness of the integral on the right-hand side
follows from the following choices of function spaces:
\begin{itemize}
    \item In the volumic case,
    $(f,g)\in L^p(\mu)\times L^{q'}(\mu)$ or $(f,g)\in L^{q'}(\mu)\times L^p(\mu)$, where
    $p,q\in(1,\infty)$ satisfy
    \(
    \frac1q=\frac1p-\varepsilon>0.
    \)

    \item In the metric case,
    $(f,g)\in L^p(\mu)\times
    \bigl(L^{p'}(\tilde w)\cap L^{q'}(\mu)\bigr)$ or $(f,g)\in
    \bigl(L^{p'}(\tilde w)\cap L^{q'}(\mu)\bigr)\times L^p(\mu)$,
    where $p,q\in(1,\infty)$ satisfy
    \(
    \frac1q=\frac1p-\frac{\eta}{d}>0,
    \)
    and $\tilde w\in A_{p'}$ is given by
    \(
    \tilde w(x)
    :=\bigl(1+\rho(x,y_0)\bigr)^{\eta p'}\),
    for a fixed $y_0\in X$.
\end{itemize}
In particular, 
\begin{equation}\label{eq:PhiNormParaprod2}
    \pair{\Pi_{T1} f}{g}=\sum_{k\in\mathbb Z}
\sum_{P,Q\in\mathscr D_k}\big\langle \D_P T\mathbf{1}_Q,g\big\rangle
\langle f\rangle_P
\end{equation}
with absolute convergence.
\end{enumerate}
\end{lemma}

\begin{proof}
To begin with notice that for any $P\in \mathscr D$ we have
\begin{equation}\label{eq:PhiNormParaprodp1}
    \sum_{m=0}^\infty \omega(\delta^m)\phi\bigl(z_P,\delta^{-m}\ell(P)\bigr) \lesssim
\begin{cases}
\displaystyle V(z_P,\ell(P))^{\eps}\underbrace{\sum_{m=0}^{\infty}\omega(\delta^m)\delta^{-mD\varepsilon}}_{<\infty}\lesssim\mu(P)^\eps= \phi(P), & \text{(volumic)}\\[1em]
\displaystyle \ell(P)^\eta\underbrace{\sum_{m=0}^\infty \omega(\delta^m)\delta^{-\eta m}}_{<\infty}\lesssim \ell(P)^\eta=\phi(P). & \text{(metric)}
\end{cases}
\end{equation}
Thus, by Lemma~\ref{lem:basicDecAux}, $\D_P T1$ is well-defined for all $P\in\mathscr D$ and 
\[
\|\D_P T1\|_{L^\infty(\mu)}\lesssim\phi(P).
\]
Thus,
\[
\|\{\D_P T1\}\|_\phi
=\sup_{P\in\mathscr D}
\frac{\|\D_P T1\|_{L^\infty(\mu)}}{\phi(P)} \lesssim 1
\]
which proves point~\ref{point:PhiNormParaprod1}.

It remains to prove point~\ref{point:PhiNormParaprod2}. Let $f$ and $g$ be functions as in the statement. Then
\begin{equation}\label{eq:PhiNormParaprodp2}
\begin{split}
     \sum_{k\in\mathbb Z}
    \sum_{P,Q\in\mathscr D_k}&
    \left|
    \big\langle  \D_P T\mathbf{1}_Q,g\big\rangle
    \langle f\rangle_P
    \right| \\
    &=\sum_{k\in\mathbb Z}
    \sum_{P\in\mathscr D_k}
    \abs{\langle f\rangle_P}
    \sum_{Q\in\mathscr D_k}\left|\big\langle  \D_P T\mathbf{1}_Q,\D_P g\big\rangle
    \right|\\
    &\leq  \sum_{k\in\mathbb Z}
    \sum_{P\in\mathscr D_k}
    \langle \abs{f}\rangle_P
    \sum_{Q\in\mathscr D_k}
     \Norm{\D_P T\mathbf{1}_Q}{L^\infty(\mu)}\Norm{\D_P g}{L^1(\mu)}\\
     &\lesssim  \sum_{k\in\mathbb Z}
    \sum_{P\in\mathscr D_k}
    \langle \abs{f}\rangle_P\Norm{\D_P g}{L^1(\mu)} \Big(\sum_{m=0}^\infty \omega(\delta^m)\phi\bigl(z_P,\delta^{-m}\ell(P)\bigr)\Big)\qquad \text{by Lemma~\ref{lem:basicDecAux}}\\
    &\lesssim \sum_{k\in\mathbb Z}
    \sum_{P\in\mathscr D_k}
   \phi(P) \langle \abs{f}\rangle_P\Norm{\D_P g}{L^1(\mu)} \qquad \text{by~\eqref{eq:PhiNormParaprodp1}}\\
    &= \sum_{P\in\mathscr D}\phi(P)
    \langle \abs{f}\rangle_P\Norm{\D_P g}{L^1(\mu)},
\end{split}
\end{equation}
where, in the last equality, we have used the fact that for all $P\in\mathscr D$ there exists a unique $k\in\Z$ such that $\D_P$ is nonzero. Now, for any $P'\in\children(P)$ we have that
\[
\abs{\ave{g}_{P'}}\leq \ave{\abs{g}}_{P'}\lesssim \ave{\abs{g}}_{P},
\]
and therefore
\begin{equation*}
    \begin{split}
    \Norm{\D_P g}{L^1(\mu)}&=\BNorm{\sum_{P'\in\children(P)}\ave{g}_{P'}\mathbf{1}_{P'}-\ave{g}_{P}\mathbf{1}_{P}}{L^1(\mu)}\\
    &\leq \sum_{P'\in\children(P)}\abs{\ave{g}_{P'}}\mu(P')+\abs{\ave{g}_{P}}\mu(P)\lesssim \ave{\abs{g}}_{P}\mu(P)    
    \end{split}
\end{equation*}
Substituting this estimate into \eqref{eq:PhiNormParaprodp2},
we obtain
\begin{equation*}
    \begin{split}
        \sum_{k\in\mathbb Z}
    \sum_{P,Q\in\mathscr D_k}&
    \left|
    \big\langle  \D_P T\mathbf{1}_Q,g\big\rangle
    \langle f\rangle_P
    \right|\\
    &\lesssim \sum_{P\in\mathscr D}\phi(P)
    \langle \abs{f}\rangle_P\ave{\abs{g}}_{P}\mu(P)\\
    &=\sum_{P\in\mathscr D} \int_X\int_X \frac{\phi(P)}{\mu(P)}\abs{f(y)}\abs{g(x)}\mathbf{1}_P(y)\mathbf{1}_P(x)\ud\mu(y)\ud\mu(x)\\
    &= \int_X\int_X \Big(\sum_{\substack{P\in \mathscr D\\ P\ni x,y}}\frac{\phi(P)}{\mu(P)}\Big)\abs{f(y)}\abs{g(x)}\ud\mu(y)\ud\mu(x)\\ 
    &\lesssim \int_X\int_X \frac{\phi(x,y)}{V(x,y)}\abs{f(y)}\abs{g(x)}\ud\mu(y)\ud\mu(x),
    \end{split}
\end{equation*}
where, in the last step, we have used Lemma~\ref{lem:sumcubesfrac} in a way analogous to what was done in the proof of Proposition~\ref{prop:DomPhiNormParaprod}. This proves formula~\eqref{eq:PhiNormParaprod1}. Moreover, the fact that the quantity on the right-hand side is finite follows from the choice of $f$ and $g$ and from Theorems~\ref{thm:bddfracvol} and~\ref{thm:bddfracmet}.

By the absolute convergence established above, we may
rearrange the series and, for fixed $k\in\Z$ and $P\in\mathscr D_k$,
first sum over $Q\in\mathscr D_k$. We thus obtain
\[
    \sum_{Q\in\mathscr D_k}
    \big\langle  \D_P T\mathbf{1}_Q,g\big\rangle
    \langle f\rangle_P=
    \big\langle \sum_{Q\in\mathscr D_k}  \D_P T\mathbf{1}_Q,g\big\rangle
    \langle f\rangle_P=\big\langle  \D_P T1,g\big\rangle
    \langle f\rangle_P.
\]
In these equalities, we use the fact that, in both the metric
and volumic cases, the series
$\sum_{Q\in\mathscr D_k}\D_P T\mathbf{1}_Q=:\D_P T1$
converges in $L^q(\mu)$ when $g\in L^{q'}(\mu)$, and in
$L^{p'}(\mu)$ when $g\in L^p(\mu)$, as follows from
Lemma~\ref{lem:basicDecAux}. Consequently, summing over $P\in\mathscr D_k$ and $k\in\mathbb{Z}$, we obtain
\[
\sum_{k\in\mathbb Z}
    \sum_{P,Q\in\mathscr D_k}   
    \big\langle \D_P T\mathbf{1}_Q,g\big\rangle
    \langle f\rangle_P=\sum_{k\in\mathbb Z}
    \sum_{P\in\mathscr D_k}
    \big\langle \D_P T1,g\big\rangle
    \langle f\rangle_P=\big\langle \Pi_{T1}f,g\big\rangle,
\]
which concludes the proof.
\end{proof}

\begin{remark}\label{rmk:PhiNormParaprodAdj}
    We observe that analogous considerations apply to the adjoint operator $T^*$, defined by the identity
\[
\langle Tf,g\rangle = \langle f,T^*g\rangle,
\]
for all admissible pairs $(f,g)$ in the corresponding spaces, where, for any appropriate function $\widetilde f$,
\[
T^*\widetilde f(x)=\int_X K(y,x)\widetilde f(y)\ud\mu(y), \quad x\in X.
\]
Notice that the operator $T^*$ is itself an $\omega$-regular $\phi$-fractional integral operator. Indeed, its kernel is given by
$K^*(x,y):=K(y,x)$, and the estimates
\eqref{eq:CZ0phi} and \eqref{eq:CZ1phi} are symmetric
in $x$ and $y$, up to constants. Moreover, $T^*$ also satisfies \eqref{eq:omegasum}, since the function $\omega$ associated with its kernel is the same as that associated with $K$. The same argument then shows that $\Pi_{T^*1}$ is
well-defined and $\phi$-normalized.
Moreover, similar computations yield
\begin{equation}\label{eq:PhiNormParaprodAdj}
    \pair{\Pi_{T^*1}^* f}{g}=\pair{f}{\Pi_{T^*1}g}=\sum_{k\in\mathbb Z}
\sum_{P,Q\in\mathscr D_k}\ave{g}_Q \pair{f}{\D_Q T^*\mathbf{1}_P},
\end{equation}
where the series converges absolutely. Such absolute convergence comes directly from point~\ref{point:PhiNormParaprod2} with $(g,f)$ in place of $(f,g)$.
\end{remark}

\subsection{Dyadic shifts and fractional normalization}

Let us introduce the definition of dyadic shift. 

\begin{definition}\label{def:shifts}
    Let $(X,\rho,\mu)$ be a space of homogeneous type. A \textit{dyadic shift of complexity $(i,j)$} is an operator defined as
\begin{equation}\label{eq:shift}
  \S = \sum_{k\in\mathbb{Z}} \sum_{R\in\mathscr{D}_k} A_{(R,k)}, \qquad
  A_{(R,k)} = \sum_{\substack{P \in \operatorname{ch}^i(R,k) \\ Q \in \operatorname{ch}^j(R,k)}} a_{PQR} \, h_P \otimes h_Q,
\end{equation}
where $h_P$ and $h_Q$ are cancellative Haar functions supported on the cubes $P$ and $Q$, respectively, and the $a_{PQR}$ are scalar coefficients. (If a cancellative Haar function does not exist on a given cube, we adopt the convention that $h_Q = 0$ or $h_P = 0$ as appropriate.) The shift $\S$ is said to be:
\begin{itemize}
\item \emph{$\phi$-normalized}, for $\phi$ as in Notation~\ref{not:phi}, and we may sometimes denote it by $\S^{\phi}$ to emphasize it, if the coefficients satisfy the estimate
\begin{equation}\label{eq:shiftphinorm}
  \abs{a_{PQR}} \lesssim \frac{\sqrt{\mu(P)\mu(Q)}}{\mu(R)}\, \phi(R).
\end{equation}
    \item \emph{$\varepsilon$-volumic normalized}, for some $\eps\in(0,1)$, and we will sometimes denote it as $\S^{(\varepsilon,\mu)}$ to emphasize it, if it is $\phi$-normalized with $\phi(R):=\mu(R)^\eps$.
\item \emph{$\eta$-metric normalized}, for some $\eta \in (0, \infty)$, and we may sometimes denote it by $\S^{(\eta,\rho)}$ to emphasize it, if it is $\phi$-normalized with $\phi(R):=\ell(R)^\eta$.
\end{itemize}
\end{definition}

Before proceeding with the unified treatment, let us examine how the two normalizations relate to each other in the specific case of the Bessel setting.

\begin{remark}
    In the case of the Bessel setting (see Subsection~\ref{subsec:Bessel}), thanks to \eqref{eq:BesInterval}, we obtain that, for every cube $R\in\mathscr D$,
    \begin{equation*}
        \ell(R)\lesssim m_\lambda(I(z_R,\ell(R)))^{\frac{1}{2\lambda+1}}\approx m_\lambda(R)^{\frac{1}{2\lambda+1}}.
    \end{equation*}
    Thus, recalling \eqref{eq:Q-BQ}, if $\S^{(\eta,\rho)}$ is an $\eta$-metric normalized dyadic shift, its coefficients satisfy
    \begin{equation*}
        \abs{a_{PQR}}\lesssim \frac{\sqrt{m_\lambda(P)m_\lambda(Q)}}{m_\lambda(R)} \ell(R)^\eta\lesssim \frac{\sqrt{m_\lambda(P)m_\lambda(Q)}}{m_\lambda(R)}m_\lambda(R)^{\frac{\eta}{2\lambda+1}}.
    \end{equation*}
    Hence, if $\eta\in\left(0,2\lambda+1\right)$, then $\S^{(\eta,\rho)}$ is also an $\frac{\eta}{2\lambda+1}$-volumic normalized dyadic shift.
\end{remark}

A first noteworthy observation is that, just as dyadic paraproducts, $\phi$-normalized dyadic shifts can also be pointwise dominated by fractional integrals.


\begin{proposition}\label{prop:shiftsdom}
    Let $(X,\rho,\mu)$ be a space of homogeneous type, and let $\S^\phi$ be a $\phi$-normalized dyadic shift of complexity $(i,j)$. In the metric setting, assume additionally that $(X,\rho,\mu)$ has positive lower dimension $d>0$ and that the metric fractional parameter satisfies $\eta\in(0,d)$. Then, for $p\in(1,\infty)$ and $f\in L^p(X)$, we have
    \[
    \abs{\S^\phi f(x)}\lesssim \bigl(1+\min(i,j)\bigr) \int_X \frac{\phi(x,y)}{V(x,y)} \abs{f(y)}\ud \mu(y), \quad \text{for $\mu$-a.e. }x\in X.
    \]
    The factor $1+\min(i,j)$ may be omitted whenever $(X,\rho,\mu)$ has positive lower dimension.
\end{proposition}

\begin{proof}
    Let $f\in L^p(\mu)$, for $p\in(1,\infty)$, and consider a $\S^\phi$ be a $\phi$-normalized dyadic shift of complexity $(i,j)$, i.e., $\S^\phi$ is as in \eqref{eq:shift} with coefficients as in \eqref{eq:shiftphinorm}.
Then, for all $(R,k)$, we have that
\begin{equation*}
    \begin{split}
        \abs{A_{(R,k)}f}&\leq \sum_{\substack{P \in \operatorname{ch}^i(R,k) \\ Q \in \operatorname{ch}^j(R,k)}} \abs{a_{PQR}} \abs{\pair{h_Q}{f}} \abs{h_P}\\
        &\lesssim \frac{\phi(R)}{\mu(R)} \sum_{Q \in \operatorname{ch}^j(R,k)}\int_Q \abs{f}\ud\mu \sum_{P \in \operatorname{ch}^i(R,k)}\mathbf{1}_P = \frac{\phi(R)}{\mu(R)} \int_R \abs{f}\ud\mu\,  \mathbf{1}_R
    \end{split}
\end{equation*}
On the other hand, by \cite[Lemma 14.6]{Hyt:26}, for each $R\in\mathscr{D}$, the operator $A_{(R,k)}$ can be non-zero only for indices $k$ in a set $\mathcal{K}_{i,j}(R)\subset\mathbb{Z}$ containing at most $1+\min(i,j)$ elements. Hence, if we denote 
\begin{equation*}
  \mathscr D^*_{i,j}
  :=\Big\{ (R,k):R\in\mathscr D, k\in\mathcal K_{i,j}(R) \Big\},
\end{equation*}
then we get that for $\mu$-almost every $x\in X$
\begin{align}
    \abs{ \S^\phi f(x)} &= \Big|\sum_{k\in\mathbb{Z}} \sum_{R\in\mathscr{D}_k} A_{(R,k)}f(x)\Big|=\Big|\sum_{(R,k)\in \mathscr D^*_{i,j}} A_{(R,k)}f(x)\Big|\nonumber\\
    &\leq \sum_{(R,k)\in \mathscr D^*_{i,j}}\frac{\phi(R)}{\mu(R)} \int_R \abs{f}\ud\mu\,  \mathbf{1}_R(x)\leq \bigl(1+\min(i,j)\bigr)\sum_{R\in \mathscr D}\frac{\phi(R)}{\mu(R)} \int_R \abs{f}\ud\mu\,  \mathbf{1}_R(x)\nonumber\\
    &=\bigl(1+\min(i,j)\bigr)\sum_{R\in \mathscr D} \int_X \frac{\phi(R)}{\mu(R)} \mathbf{1}_R(x)\mathbf{1}_R(y)\abs{f(y)}\ud\mu(y)\nonumber\\
    &=\bigl(1+\min(i,j)\bigr)\int_X \Big(\sum_{\substack{R\in \mathscr D\\ R\ni x,y}} \frac{\phi(R)}{\mu(R)}\Big)\abs{f(y)}\ud\mu(y)\nonumber\\
    &\lesssim\bigl(1+\min(i,j)\bigr) \int_X \frac{\phi(x,y)}{V(x,y)} \abs{f(y)}\ud \mu(y), \nonumber
\end{align}
where, in the last step, we have used Lemma~\ref{lem:sumcubesfrac} in a way analogous to what was done in the proof of Proposition~\ref{prop:DomPhiNormParaprod}. Finally, if $(X,\rho,\mu)$ has positive lower dimension,
the factor $1+\min(i,j)$ can be omitted, since each cube
$R\in\mathscr D$ belongs to a unique generation
$\mathscr D_k$. In particular, the index $k$ in
$(R,k)\in\mathscr D^*_{i,j}$ is uniquely determined by
$R$, and hence $\#\mathcal K_{i,j}(R)=1$. This concludes the proof.
\end{proof}

This result immediately implies that every $\phi$-normalized dyadic shift inherits the same $L^p$-boundedness properties as $\phi$-fractional integrals.

\begin{corollary}
    Let $(X,\rho,\mu)$ be a space of homogeneous type. Then:
    
    \begin{itemize}
        \item If $\S$ is a $\varepsilon$-volumic normalized dyadic shift of complexity $(i,j)$, for $\eps\in (0,1)$, then it is bounded from \( L^p(X,\mu) \) to \( L^q(X,\mu) \) for all $p,q\in(1,\infty)$ such that $\frac1q=\frac1p-\eps>0$. Moreover,
        \[
        \|\S\|_{L^p(\mu)\to L^q(\mu)}
        \lesssim 1+\min(i,j),
        \]
        where the implicit constant is independent of the
        complexity $(i,j)$. If $(X,\rho,\mu)$ has positive
        lower dimension, the right-hand side can be replaced
        by a constant independent of $(i,j)$.
        \item If $(X,\rho,\mu)$ has lower dimension $d>0$ and $\S$ is an $\eta$-metric normalized dyadic shift of complexity $(i,j)$, for $\eta\in(0,d)$, then it is is bounded from $L^p(\mu)$ to $L^p(w)+L^q(\mu)$ for all $p,q\in(1,\infty)$ such that $\frac{1}{q}=\frac{1}{p}-\frac{\eta}{d}>0$, with
        \[
        \|\S\|_{L^p(\mu)\to L^p(w)+L^q(\mu)}
        \lesssim 1,
        \]
        where the implicit constant is independent of the
        complexity $(i,j)$. Here, $w$ is the Muckenhoupt weight in $A_p$ given by \( w(x):=\bigl(1+\rho(x,x_0)\bigr)^{-\eta p}\) for any fixed $x_0\in X$.
    \end{itemize}
\end{corollary}

\begin{proof}
    The claim follows directly by combining Proposition~\ref{prop:shiftsdom} and Theorems~\ref{thm:bddfracvol} and~\ref{thm:bddfracmet}.
\end{proof}

\section{Dyadic formula}\label{sec:DRT}

In this section we prove a dyadic representation formula for fractional integrals. As already discussed in more detail in the introductory section, several dyadic formulas appear in the literature. In our case, however, we draw inspiration from the one established in \cite[Theorem 16.1]{Hyt:26} for Calderón--Zygmund operators, adapting the techniques developed there to the setting of fractional integrals. 

\subsection{Statement of the formula}

Let us begin by stating the formula. 

\begin{theorem}\label{thm:DRT}
Let $(X,\rho,\mu)$ be a space of homogeneous type, and let $T$ be an $\omega$-regular $\phi$-fractional integral operator in the sense of Definition~\ref{def:fracint}, where $\omega : [0,1] \to [0,\infty)$ satisfies Assumption~\ref{ass:omega} and condition \eqref{eq:omegasum}. In the metric setting, assume in addition that $(X,\rho,\mu)$ has positive lower dimension $d>0$, and that the metric fractional parameter satisfies $\eta \in (0,d)$. Then $T$ admits the following representation:
\begin{equation}\label{eq:DRT}
\begin{split}
  T=T_1 &+\mathbb{E}\Big(\Pi_{T1}+\Pi^*_{T^*1}\Big) \\
  &+\mathbb{E}\left(\sum_{m=0}^\infty\omega(\delta^m)\Big\{
  \S_{m,m}^{\phi}+  \sum_{i=0}^{m-1}(\S_{m,i}^{\phi}+\S_{i,m}^{\phi})+  \sum_{i=1}^m(\S_{0,i}^{\phi}+\S_{i,0}^{\phi})\Big\}\right),
\end{split}
\end{equation}
where:
\begin{enumerate}
  \item $T_1$ is an operator of rank at most one, explicitly given by
\begin{equation*}
  T_1=\begin{cases} 
  \mu(X)^{-2}\langle T1, 1 \rangle \, 1 \otimes 1, & \text{if } \mu(X) < \infty, \\ 
  0, & \text{if } \mu(X) = \infty; 
  \end{cases}
\end{equation*} 
  \item $\Pi_{T1}$ and $\Pi_{T^*1}$ denote $\phi$-normalized dyadic paraproducts with symbol $T1$ and $T^*1$, in the sense of Definition~\ref{def:paraT1}.
\item each $\S_{m,n}^{\phi}$ denotes a (sum of $O(1)$ many) $\phi$-normalized dyadic shift(s) of complexity $(m,n)$. 

\item The expectation \(\mathbb E\) is taken with respect to the underlying probability space of random dyadic systems.
\end{enumerate}
This representation formula should be interpreted as follows, depending on the two types of operators encoded by $\phi$:
\begin{itemize}
    \item In the case where $T$ is a volumic fractional integral associated with some $\varepsilon \in (0,1)$, it is to be understood as an equality of bilinear forms on \( L^p(\mu) \times L^{q'}(\mu) \), for \( p, q \in (1,\infty) \) such that \( \frac{1}{q} = \frac{1}{p} - \varepsilon > 0 \).
    
    \item In the case where $T$ is a metric fractional integral associated with some \( \eta \in (0,d) \), it is to be understood as an equality of bilinear forms on \( L^p(\mu) \times \big(L^{p'}(\tilde w) \cap L^{q'}(\mu)\big) \), for \( p, q \in (1,\infty) \) such that \( \frac{1}{q} = \frac{1}{p} - \frac{\eta}{d} > 0 \). Here, \( \tilde w \in A_{p'}\) is the Muckenhoupt weight given by the explicit formula $ \tilde w(x):=\bigl(1+\rho(x,y_0)\bigr)^{\eta p'}$, for some fixed $y_0\in X$.
\end{itemize}
\end{theorem}

\begin{remark}
Given two normed spaces \(X\) and \(Y\) within a common topological vector space, the norm on their intersection is defined in the usual way by
\[
\|z\|_{X\cap Y} := \max \bigl\{ \|z\|_X ,\|z\|_Y \bigr\}.
\]
This definition is the standard one in the literature (see, for instance, \cite[Lemma 2.3.1]{BL:76}). Thus, in particular, this is how we interpret the space \( L^{p'}(\tilde w) \cap L^{q'}(\mu) \) in the statement of Theorem~\ref{thm:DRT} and throughout the subsequent discussion.
\end{remark}

\begin{remark}
In applications, the function $\omega$ most commonly takes the form $\omega(t):=t^\beta$, for some $\beta>0$. We thus examine how condition \eqref{eq:omegasum} specializes in this case.

\begin{itemize}
    \item \textbf{Metric case.}\\
Here, $\phi(x_0,\delta^{-m})=\delta^{-m\eta}$, and condition \eqref{eq:omegasum} reduces to
\[
\sum_{m=0}^{\infty} \delta^{m\beta}\delta^{-m\eta}<\infty 
\quad \Leftrightarrow \quad 
\beta>\eta.
\]

For instance, in the one-dimensional Bessel setting considered in Subsection~\ref{subsec:Bessel}, the operator $\Delta_\lambda^{-\eta/2}$ is $\omega$-regular with $\omega(t)=t$ (i.e., $\beta=1$) and $\eta$-metric with $\eta\in(0,1)$, so the above condition is automatically satisfied. 

    \item \textbf{Volumic case.}\\
Here, condition \eqref{eq:omegasum} becomes
\[
\sum_{m=0}^{\infty} \delta^{m\beta}\delta^{-mD\varepsilon}<\infty 
\quad \Leftrightarrow \quad 
\beta>D\varepsilon.
\]
\end{itemize}
In summary, in both settings, condition \eqref{eq:omegasum} encodes a compatibility between the smoothness exponent $\beta$ of $\omega$ and the fractional parameter governing $\phi$, ensuring that kernel regularity dominates its fractional scaling. In the volumic case, the upper dimension also enters this balance.
\end{remark}

\begin{remark}
     Notice that \( \tilde w \in A_{p'}\) is the dual weight of the $A_p$ Muckenhoupt weight given by $ w(x):=\bigl(1+\rho(x,y_0)\bigr)^{-\eta p}$, for some fixed $y_0\in X$, appearing in Theorem~\ref{thm:bddfracmet}.
\end{remark}

Before proceeding with the proof of the formula, let us state a useful Lemma which, as already done for paraproducts and shifts, provides an upper bound for the operator $T_1$ appearing in \eqref{eq:DRT}.

\begin{lemma}\label{lem:Rank1Est}
    If $\mu(X)<\infty$, then
    \begin{equation}\label{eq:Rank1Est1}
        \abs{\langle T1, 1 \rangle}\lesssim \phi(X)\mu(X),
    \end{equation}
    and, for any $p\in(1,\infty)$ and $f\in L^p(\mu)$,
    \begin{equation}\label{eq:Rank1Est2}
        \abs{T_1 f(x)}\lesssim \int_X \frac{\phi(x,y)}{V(x,y)}\abs{f(y)}\ud\mu(y), \quad x\in X.
    \end{equation}
\end{lemma}

\begin{proof}
    First off, since $\mu(X) < \infty$, we have $X \in \mathscr{D}$. Thus, \eqref{eq:Rank1Est1} is simply the special case of Lemma~\ref{lem:QuadFracIntCubes} where $Q = X$. 
    
    Now, in the volumic case, for $x\in X$ we therefore have 
    \begin{equation*}
    \begin{split}
        \abs{T_1 f(x)}&\lesssim \int_X \frac{\phi(X)}{\mu(X)} \abs{f(y)} \ud\mu(y)= \int_X \frac{1}{\mu(X)^{1-\eps}} \abs{f(y)} \ud\mu(y)\\&\leq \int_X \frac{1}{V(x,y)^{1-\eps}} \abs{f(y)} \ud\mu(y) =\int_X \frac{\phi(x,y)}{V(x,y)}\abs{f(y)}\ud\mu(y)
    \end{split}
    \end{equation*}
    On the other hand, in the metric case
    \begin{equation}\label{eq:Rank1Estp1}
        \abs{T_1 f(x)}\lesssim \int_X \frac{\phi(X)}{\mu(X)} \abs{f(y)} \ud\mu(y)=\int_X \frac{\ell(X)^\eta}{\mu(X)} \abs{f(y)} \ud\mu(y)\lesssim \int_X \frac{\diam(X)^\eta}{\mu(X)} \abs{f(y)} \ud\mu(y).
    \end{equation}
    By lower dimension $d>0$ and Lemma~\ref{lem:noatoms} we have $\mu(\{x\})=0$, and, for $x\neq y$, 
    \[
    \frac{\mu(X)}{V(x,y)}=\frac{V(x,2\diam(X))}{V(x,y)}\gtrsim \Big(\frac{\diam(X)}{\rho(x,y)}\Big)^d.
    \]
    Thus,
    \[
    \frac{\diam(X)^\eta}{\mu(X)}\lesssim \frac{\rho(x,y)^d\diam(X)^{\eta-d}}{V(x,y)}\leq \frac{\rho(x,y)^\eta}{V(x,y)}=\frac{\phi(x,y)}{V(x,y)}.
    \]
    Substituting this estimate into \eqref{eq:Rank1Estp1} completes the proof.
\end{proof}

\subsection{Proof of the formula}\label{subsec:proofDRT}

The proof of formula \eqref{eq:DRT} requires several intermediate steps, and most of the main ideas follow the example of those presented in \cite[Section 18]{Hyt:26}. 

Throughout this section, whenever we refer
to a dyadic cube $Q$, we implicitly mean the pair $(Q,k)$, where $k$
denotes the level of $Q$. 

\subsubsection*{An initial decomposition}

First, recall the definitions of $\E_Q$, $\D_Q$, $\E_k$ and $\D_k$ from \eqref{eq:EQ}, \eqref{eq:DQ} and \eqref{eq:EkDk}. By Lemma~\ref{lem:convcondi}, for $f\in L^p(w)$, with $p\in(1,\infty)$ and $v\in A_p$, we have that 
\begin{equation*}
   \E_k f\to\begin{cases} f, & k\to\infty, \\ \ave{f}_X, & k\to-\infty,\end{cases}\qquad
   \ave{f}_X:=\begin{cases} 0, & \mu(X)=\infty, \\ \fint_X f\ud\mu, & \mu(X)<\infty,\end{cases}
\end{equation*}
and such convergences are to be understood in the $L^p(v)$ norm.

Now, let $T$ be a fractional integral operator, and, in the metric case, suppose that $(X,\rho,\mu)$ has positive lower dimension $d>0$, and that $\eta \in (0,d)$, as stated in Theorem \ref{thm:DRT}. We know that $T$, under both possible definitions (volumic and metric), satisfies certain boundedness properties. More precisely, we have:

\begin{itemize}
    \item If we regard $T$ as a volumic fractional integral associated with some $\varepsilon \in (0,1)$, then by Theorem~\ref{thm:bddfracvol} we know that $T$ is bounded from \( L^p(\mu) \) to \( L^q(\mu) \) for all \( p,q \in (1,\infty) \) such that \( \frac1q = \frac1p - \varepsilon > 0 \).
    
    \item If we regard $T$ as a metric fractional integral associated with some \( \eta \in (0,d) \), then by Theorem~\ref{thm:bddfracmet} $T$ is bounded from \( L^p(\mu) \) to \( L^p(w) + L^q(\mu) \), for \( p,q \in (1,\infty) \) satisfying \( \frac{1}{q} = \frac{1}{p} - \frac{\eta}{d} > 0 \), for \( w(x):=\bigl(1+\rho(x,y_0)\bigr)^{-\eta p}\in A_p \).
\end{itemize}

Next, recall that the dual of the sum of spaces \( L^p(w) + L^q(\mu) \) is given by
\[
\big(L^p(w) + L^q(\mu)\big)^* = L^p(w)^* \cap L^q(\mu)^* = L^{p'}(\tilde w) \cap L^{q'}(\mu),
\]
where \( p', q' \in (1,\infty) \) are the Hölder conjugates of \( p \) and \( q \), respectively, and 
\[
\tilde w(x) = w(x)^{-\frac{1}{p-1}}=\bigl(1+\rho(x,y_0)\bigr)^{\eta p'}\in A_{p'}
\]
is the dual weight of $w$.

Therefore, with suitable choices of \( p \) and \( q \), we obtain that for every \( f \in L^p(\mu) \) and every \( g \) belonging to \( L^{q'}(\mu) \) in the volumic case, or to \( L^{p'}(\tilde w) \cap L^{q'}(\mu) \) in the metric case, thanks to Lemma \ref{lem:convcondi}, the following convergences hold:
\begin{equation*}
  \pair{T\E_k f}{\E_k g} \to \begin{cases} \pair{Tf}{g}, & k \to \infty, \\ \ave{f}_X \ave{g}_X \pair{T1}{1}, & k \to -\infty. \end{cases}
\end{equation*}
Thus, in the sense of bilinear forms on \( L^p(\mu) \times L^{q'}(\mu) \) in the volumic case, and on \( L^p(\mu) \times \big(L^{p'}(\tilde w) \cap L^{q'}(\mu)\big) \) in the metric case, we can write
\begin{equation}\label{eq:T0T1}
  T = \lim_{\substack{b \to \infty \\ a \to -\infty}} \left( \E_b T \E_b - \E_a T \E_a \right) + \pair{T1}{1} \, \frac{1}{\mu(X)} \otimes \frac{1}{\mu(X)} =: T_0 + T_1,
\end{equation}
where the term $T_1$ corresponds to that in Theorem~\ref{thm:DRT}; note that $T_1$, and thus also $T_0 = T - T_1$, is independent of the choice of dyadic system used in this decomposition.

\subsubsection*{Identification of the terms with the dyadic paraproducts and shifts}

We begin by stating and proving the analogue of \cite[Lemma~18.2]{Hyt:26}.

\begin{lemma}\label{lem:basicDec}
Let $T_0$ be defined as in \eqref{eq:T0T1}. Then
\[
\begin{split}
    T_0
    &= \Pi_{T1}+\Pi_{T^*1}^*+T_{00},\\
    T_{00}
    &=
    \sum_{k\in\mathbb Z}
    \sum_{P,Q\in\mathscr D_k}
    \left(
        \D_PT\D_Q
        +\D_{P,Q}^*T\D_Q
        +\D_PT\D_{Q,P}
    \right),
\end{split}
\]
where $\Pi_{T1}$ and $\Pi_{T^*1}$ denote the dyadic paraproducts associated,
respectively, with the families
$\{\D_QT1\}_{Q\in\mathscr D}$ and
$\{\D_QT^*1\}_{Q\in\mathscr D}$ in the sense of
Definition~\ref{def:paraprod}, and
\begin{equation}\label{eq:DQP}
  \D_{Q,P}
  :=
  \mathbf1_Q\otimes
  \left(
      \frac{\mathbf1_Q}{\mu(Q)}
      -
      \frac{\mathbf1_P}{\mu(P)}
  \right).
\end{equation}
These identities are to be understood in the sense of bilinear forms,
according to the interpretation described in the final part of the statement
of Theorem~\ref{thm:DRT}.
\end{lemma}

\begin{proof}
    Starting from the definition given in \eqref{eq:T0T1}, we apply the telescoping sum identity to expand the difference:
\begin{equation*}
\begin{split}
  \E_b T \E_b - \E_a T \E_a
  &= \sum_{k=a}^{b-1} \left( \E_{k+1} T \E_{k+1} - \E_k T \E_k \right) \\
  &= \sum_{k=a}^{b-1} \left( \D_k T \D_k + \E_k T \D_k + \D_k T \E_k \right) \qquad\qquad\text{[as }\E_{k+1}=\E_k+\D_k] \\
  &= \sum_{k=a}^{b-1} \sum_{P, Q \in \mathscr D_k}
  \left( \D_P T \D_Q + \E_P T \D_Q + \D_P T \E_Q \right).
\end{split}
\end{equation*}
As a consequence, we obtain the representation
\begin{equation}\label{eq:T0exp}
  T_0 = \sum_{k \in \Z} \sum_{P, Q \in \mathscr D_k}
  \left( \D_P T \D_Q + \E_P T \D_Q + \D_P T \E_Q \right).
\end{equation}

In order to handle the localized operators, we proceed by further decomposing the term $\D_P T \E_Q$ as follows:
\begin{equation}\label{eq:LocDec}
\begin{split}
  \D_P T \E_Q
  &= \D_P T \mathbf{1}_Q \otimes \frac{\mathbf{1}_Q}{\mu(Q)} \\
  &= \D_P T \mathbf{1}_Q \otimes \left( \frac{\mathbf{1}_Q}{\mu(Q)} - \frac{\mathbf{1}_P}{\mu(P)} \right)
  + \D_P T \mathbf{1}_Q \otimes \frac{\mathbf{1}_P}{\mu(P)} \\
  &= \D_P T \D_{Q,P} + \D_P T \mathbf{1}_Q \otimes \frac{\mathbf{1}_P}{\mu(P)}.
\end{split}
\end{equation}


Now, by Lemma \ref{lem:PhiNormParaprod}, point \ref{point:PhiNormParaprod2}, for $f \in L^p(\mu)$ and $g$ in the appropriate dual space we have that
$$
\sum_{k\in\mathbb Z} \sum_{P,Q\in\mathscr D_k} \D_P T\mathbf{1}_Q \otimes \frac{\mathbf{1}_P}{\mu(P)}(f,g)=\sum_{k\in\mathbb Z}
\sum_{P,Q\in\mathscr D_k}\big\langle \D_P T\mathbf{1}_Q,g\big\rangle
\langle f\rangle_P=\pair{\Pi_{T1} f}{g}
$$
Consequently, in the sense of bilinear forms, this yields the identity
\begin{equation}\label{eq:FinLocDec1}
  \sum_{k\in\mathbb Z}\sum_{P, Q \in \mathscr D_k} \D_P T \E_Q
  = \Pi_{T1}+\sum_{k\in\mathbb Z}\sum_{P, Q \in \mathscr D_k} \D_P T \D_{Q,P}.
\end{equation}

We now turn to the analysis of the operators $\E_P T \D_Q$. Again, we cannot argue directly at the level of operators, but at the level of bilinear forms. For every $P,Q \in \mathscr D_k$, by the same computation made in \eqref{eq:LocDec}, we have
\begin{equation}\label{eq:LocDec2}
    \langle \E_P T \D_Q f, g \rangle
= \langle f, \D_Q T^* \E_P g \rangle=\pair{f}{\D_QT^*\D_{P,Q}g}+ \ave{g}_Q \pair{f}{\D_Q T^*\mathbf{1}_P}=:I_{P,Q}+J_{P,Q}
\end{equation}
Now,
\begin{equation*}
    I_{P,Q}=\pair{\D_{P,Q}^*T\D_Qf}{g}.
\end{equation*}
On the other hand, by Remark~\ref{rmk:PhiNormParaprodAdj},
\begin{equation*}
    \sum_{k\in\mathbb Z}
\sum_{P,Q\in\mathscr D_k}J_{P,Q}= \sum_{k\in\mathbb Z}
\sum_{P,Q\in\mathscr D_k}\ave{g}_Q \pair{f}{\D_Q T^*\mathbf{1}_P}= \pair{\Pi_{T^*1}^* f}{g}.
\end{equation*}
Therefore, in the sense of bilinear forms, we get
\begin{equation}\label{eq:FinLocDec2}
    \sum_{k\in\mathbb Z}\sum_{P,Q \in \mathscr D_k}\E_P T \D_Q=\Pi_{T^*1}^*+\sum_{k\in\mathbb Z}\sum_{P,Q \in \mathscr D_k}\D_{P,Q}^*T\D_Q.
\end{equation}

At this point, by combining \eqref{eq:FinLocDec1} and \eqref{eq:FinLocDec2}, we get that
\begin{equation}\label{eq:extractParap}
\begin{split}
  \sum_{k\in\mathbb Z}\sum_{P,Q\in\mathscr D_k}
  &\Big(\D_PT\D_Q+\E_PT\D_Q+\D_PT\E_Q\Big) \\
  &=\Pi_{T1}+\Pi_{T^*1}^*+\sum_{k\in\mathbb Z}\sum_{P,Q\in\mathscr D_k}\Big(\D_PT\D_Q+\D_{P,Q}^*T\D_Q+\D_PT\D_{Q,P}\Big).
\end{split}
\end{equation}
Finally, putting together \eqref{eq:T0exp} and \eqref{eq:extractParap} leads to the identities stated in the lemma.
\end{proof}

Now that we have established this result, the next step is to recall \cite[Lemma~18.8]{Hyt:26}, which remains valid in our setting with identical proof. For this reason, we refer the reader to the original source for the details. This lemma allows us to rearrange the cubes appearing in the decomposition obtained in Lemma~\ref{lem:basicDec}, and to compute an average over all possible choices of random dyadic systems.

\begin{lemma}\label{lem:underCommon}
Let
\begin{equation*}
  T_{P,Q}\in\left\{\D_P T \D_Q,\quad \D_{P,Q}^* T \D_Q,\quad \D_P T \D_{Q,P}\right\}.
\end{equation*}
Then the expectation $\E$ taken over a random dyadic system $\mathscr D$ satisfies the identity
\begin{equation}\label{eq:underCommon}
  \E \sum_{k\in\Z} \sum_{P,Q\in\mathscr D_k} T_{P,Q}
  = \E \sum_{m=m_0}^\infty \sum_{k\in\Z} \sum_{\substack{R\in\mathscr D_k \\ P,Q \in \operatorname{ch}^m(R)}}
  \theta_m(P,Q) \, T_{P,Q},
\end{equation}
where $m_0$ is the constant, depending only on the underlying space, appearing in \cite[Lemma~17.7]{Hyt:26}, and $\theta_m(P,Q)$ are coefficients belonging to $\{0\} \cup [1,2]$ that depend on the probabilistic structure of the dyadic grids. Furthermore, for all $m > m_0$, we have the implication
\begin{equation}\label{eq:thetaNon0}
  \theta_m(P,Q) \neq 0 \quad \Rightarrow \quad \rho(x_P,x_Q) \gtrsim \delta^{-m} \ell(P).
\end{equation}
\end{lemma}

The results that follow allow us to identify the terms appearing on the right-hand side of \eqref{eq:underCommon} with the $\phi$-normalized dyadic shifts, as defined in Definition~\ref{def:shifts}, for each of the possible cases of $T_{P,Q}$ appearing in Lemma~\ref{lem:underCommon}. 

We begin by addressing the case \( T_{P,Q} = \D_P T \D_Q \). This result is the counterpart of \cite[Lemma~18.13]{Hyt:26}.

\begin{lemma}\label{lem:P-Q-shiftsPsi}
Let $m \geq m_0$ and let $\theta_m(P,Q)$ be as defined in Lemma~\ref{lem:underCommon}. Then the following identity holds:
\begin{equation*}
  \sum_{k\in\Z}\sum_{\substack{R\in\mathscr D_k \\ P,Q\in\operatorname{ch}^m(R)}}\theta_m(P,Q)\D_PT\D_Q
  =\omega(\delta^m)\S_{m,m}^{\phi},
\end{equation*}
where $\S_{m,m}^{\phi}$ denotes a sum of $O(1)$ many $\phi$-normalized dyadic shifts of complexity $(m,m)$.
\end{lemma}

\begin{proof}
    To begin with, let us expand $\D_P$ and $\D_Q$ in terms of the Haar functions, so that
    \begin{equation*}
\theta_m(P,Q)\D_PT\D_Q=\sum_{\substack{1\leq\alpha<M_P \\ 1\leq\beta<M_Q}} \theta_m(P,Q)
   \pair{h_P^\alpha}{Th_Q^\beta} h_P^\alpha\otimes h_Q^\beta.
    \end{equation*}
    Now, by Lemma~\ref{lem:underCommon}, if $m>m_0$ and $\theta_m(P,Q)\neq0$, then for $R\in\mathscr D$ such that $P,Q\in\operatorname{ch}^m(R)$ we have
    \begin{equation*}
\ell(R)= \delta^{-m} \ell(P)\lesssim \rho(x_P, x_Q) \leq \rho(x_P, x_Q)+\ell(P) \lesssim \ell(R).
\end{equation*}
In the special case where $m = m_0$, we just have
\[
\ell(R)\gtrsim \ell(P) + \rho(x_P, x_Q) \geq \ell(P) \approx \ell(R).
\]
where the implicit constants absorb the dependence on $m_0$. Thus, in both cases, we get that
\begin{equation*}
    \rho(x_P, x_Q)+\ell(P)\approx \ell(R).
\end{equation*}
Hence, 
    \begin{equation*}
  \phi\bigl(x_P,\ell(P)+\rho(x_P,x_Q)\bigr)=\left\{
\begin{aligned}
& V\bigl(x_P,\ell(P)+\rho(x_P,x_Q)\bigr)^\eps\approx  V\bigl(x_R,\ell(R)\bigr)^\eps\approx \mu(R)^\eps\\
& \bigl(\ell(P)+\rho(x_P,x_Q)\bigr)^\eta\approx \ell(R)^\eta 
\end{aligned}
\right\}\approx \phi(R).
    \end{equation*}
    Thus, combining these observations with  \eqref{eq:cubes4}, we see that in each term where $\theta_m(P,Q)\neq 0$ we have
    \begin{equation}\label{eq:HaarEstEps1}
\begin{split}
    |\pair{h_P^\alpha}{Th_Q^\beta}| &\lesssim 
    \omega\left( \frac{\ell(P)}{\rho(x_Q,x_P)+\ell(P)} \right)
    \frac{\sqrt{\mu(P)\mu(Q)}}{V(x_P,\rho(x_Q,x_P)+\ell(P))} \phi\bigl(x_P,\ell(P)+\rho(x_P,x_Q)\bigr)\\
    &\lesssim 
     \omega(\delta^m)\frac{\sqrt{\mu(P)\mu(Q)}}{\mu(R)}\,\phi(R).
\end{split}
\end{equation}
    Since both \( M_P \) and \( M_Q \) are bounded above by \( M \), the double sum  
\[
\sum_{P,Q\in\operatorname{ch}^m(R)}\theta_m(P,Q)\D_PT\D_Q
\]  
can be expressed as a finite sum of series of the form  
\[
\omega(\delta^m)\sum_{P,Q\in\operatorname{ch}^m(R)} a_{PQR}h_P\otimes h_Q,
\]  
where \( h_P = h_P^\alpha \) and \( h_Q = h_Q^\beta \) are Haar functions, and the coefficients satisfy  
\[
\abs{a_{PQR}} \lesssim \frac{\sqrt{\mu(P)\mu(Q)}}{\mu(R)}\, \phi(R).
\]  
This structure matches exactly that of a component \( A_R \) of a $\phi$-normalized dyadic shift of complexity \( (m,m) \). Therefore, we conclude that  
\[
\sum_{k\in\Z}\sum_{R\in\mathscr D_k}\sum_{P,Q\in\operatorname{ch}^m(R)}\theta_m(P,Q) \D_P T \D_Q 
= \omega(\delta^m)\S_{m,m}^{\phi},
\]  
which proves the desired identity.
\end{proof}

Next, let us consider the other possible cases of $T_{P,Q}$. Since the modifications we made in Lemma~\ref{lem:P-Q-shiftsPsi}, compared to the corresponding result in the singular integral setting of \cite{Hyt:26}, are essentially the same as those needed to prove the next two results, we will not reproduce their proofs in full. Instead, we will simply highlight the necessary changes.

\begin{lemma}\label{lem:P-QP-shiftsPsi}
Let $m \geq m_0$ and let $\theta_m(P,Q)$ be as in Lemma \ref{lem:underCommon}. Then the following identity holds:
\begin{equation*}
  \sum_{k\in\Z}\sum_{R\in\mathscr D_k}
  \sum_{P,Q\in\operatorname{ch}^m(R)}\theta_m(P,Q)
  \D_PT\D_{Q,P}
  =\omega(\delta^m)\Bigg(\sum_{i=0}^{m-1}\S_{m,i}^{\phi} + \sum_{j=1}^{m}\S_{j,0}^{\phi}\Bigg),
\end{equation*}
where each $\S_{a,b}^{\phi}$ is a sum of $O(1)$ many $\phi$-normalized dyadic shifts of complexity $(a,b)$.
\end{lemma}

\begin{proof}
    The proof of this statement is identical to that of \cite[Lemma 18.16]{Hyt:26}, which is the corresponding result in the case of singular integrals, with the only difference being that we use formula~\eqref{eq:HaarEstEps1} in place of formula \cite[(18.14)]{Hyt:26}, which is used in the original source.
\end{proof}

Similarly, we also obtain the following result, in analogy to the original case. It corresponds to \cite[Lemma 18.22]{Hyt:26} and is essentially a symmetric counterpart of Lemma~\ref{lem:P-QP-shiftsPsi}.

\begin{lemma}\label{lem:PQ-Q-shiftsPsi}
Let $m \geq m_0$ and let $\theta_m(P,Q)$ be as defined in Lemma \ref{lem:underCommon}. Then,
\begin{equation}\label{eq:PQ-Q-shifts}
  \sum_{k\in\Z}\sum_{\substack{R\in\mathscr D_k \\ P,Q\in\operatorname{ch}^m(R)}} \theta_m(P,Q)\D_{P,Q}^*T\D_Q
  = \omega(\delta^m)\Bigg(\sum_{i=0}^{m-1} \S_{i,m}^{\phi} + \sum_{j=1}^m \S_{0,j}^{\phi}\Bigg),
\end{equation}
where each $\S_{a,b}^{\phi}$ represents a sum of $O(1)$ many $\phi$-normalized dyadic shifts of complexity $(a,b)$.
\end{lemma}

At this point, we have everything we need to prove Theorem~\ref{thm:DRT}. It will be enough to combine all the parts shown in the previous lemmas.

\begin{proof}[Proof of Theorem~\ref{thm:DRT}]
First of all, by \eqref{eq:T0T1} and Lemma \ref{lem:basicDec}, we have that for any fixed dyadic system $\mathscr D$ there holds
\begin{equation*}
  T = T_1 + \Pi_{T1} + \Pi_{T^*1}^* + \sum_{\substack{k \in \Z \\ P,Q \in \mathscr D_k}}
  \Big(\D_PT\D_Q + \D_{P,Q}^*T\D_Q + \D_PT\D_{Q,P}\Big).
\end{equation*}
Thus, we can take the expectation of both sides under a random choice of $\mathscr D$, and by doing so it is important to observe that $T$ and $T_1$ are independent of $\mathscr D$. So, we get that
\begin{equation*}
\begin{split}
  T &- T_1 - \E(\Pi_{T1} + \Pi_{T^*1}^*) = \E\sum_{\substack{k \in \Z \\ P,Q \in \mathscr D_k}}
  \Big(\D_PT\D_Q + \D_PT\D_{Q,P} + \D_{P,Q}^*T\D_Q\Big) \\
  &= \E\sum_{m=m_0}^\infty \sum_{\substack{k \in \Z, R \in \mathscr D_k \\ P,Q \in \operatorname{ch}^m(R)}}
  \theta_m(P,Q)\Big(\D_PT\D_Q + \D_PT\D_{Q,P} + \D_{P,Q}^*T\D_Q\Big) \\
  &\quad\quad \text{by Lemma \ref{lem:underCommon}} \\
  &= \E\sum_{m=m_0}^\infty \omega(\delta^m)\Big(\S_{m,m}^{\phi}
  +[\sum_{i=0}^{m-1}\S_{m,i}^{\phi} + \sum_{j=1}^m\S_{j,0}^{\phi}]
  +[\sum_{i=0}^{m-1}\S_{i,m}^{\phi} + \sum_{j=1}^m\S_{0,j}^{\phi}]\Big) \\
  &\quad\quad \text{by Lemmas \ref{lem:P-Q-shiftsPsi}, \ref{lem:P-QP-shiftsPsi}, and \ref{lem:PQ-Q-shiftsPsi}}.
\end{split}
\end{equation*}
This almost concludes the proof, as moving the terms $T_1$ and $\E(\Pi_{T1} + \Pi_{T^*1}^*)$ to the right-hand side gives a formula very similar to the one stated in Theorem \ref{thm:DRT}: the only difference is that in the original statement, the summation starts from $m=0$ rather than from $m=m_0$. Therefore, to obtain exactly that formula, it suffices to set the extra shifts to be identically zero.

Finally, the $\phi$-normalization property of the paraproducts $\Pi_{T1}$ and $\Pi_{T^*1}$ follows directly from Lemma~\ref{lem:PhiNormParaprod}, point~\ref{point:PhiNormParaprod1}.
\end{proof}

\appendix
\appendixtitleformat
\renewcommand{\thesection}{\Alph{section}} 
\renewcommand{\thesubsection}{\thesection.\arabic{subsection}} 
\numberwithin{equation}{subsection} 

\section{The infinite slab}\label{app:Slab}
In this appendix we present a concrete and rather simple example of a space
of homogeneous type, namely an infinite slab embedded in $\R^n$.
The construction and analysis of this space, carried out rigorously below,
reveal several significant features of fractional integral operators,
summarized in the following.

\begin{FeaturesSlab}
The infinite slab, as defined in Definition~\ref{def:infslab},
has the following properties:
    \begin{enumerate}
    \item\label{item:slab1} It is a space of homogeneous type
    with lower dimension $d=1$ and upper dimension $D=n$.
    
    \item\label{item:slab2} In this space, the metric definition of fractional integrals is not
    a special case of the volumic one, for any choice of the parameters.
        
    \item\label{item:slab3} In this space, for every $p,q\in(1,\infty)$,
    metric fractional integrals are, in general, unbounded from $L^p$ to $L^q$.
\end{enumerate}
\end{FeaturesSlab}

In view of these properties, we conclude that even in spaces that are not particularly pathological---such as spaces with positive lower and upper dimensions---one cannot in general expect classical boundedness from $L^p$ to $L^q$ for metric fractional integrals. This, in turn, provides some justification for results such as
Theorem~\ref{thm:bddfracmet}, which may at first glance appear somewhat more unconventional (and weaker) than the corresponding statement for volumic operators.

Beyond boundedness considerations, property~\ref{item:slab2} complements the observation made in Subsection~\ref{subsec:Bessel}, where the opposite
phenomenon occurs: in the Bessel setting, volumic fractional integrals cannot be regarded as a special case of the metric ones. Combining the conclusions drawn from these two examples, we deduce that the two notions of fractional
integral are in general not comparable. This confirms the relevance of considering both definitions, without privileging one a priori.

Having summarized its properties, we now proceed to a detailed study of this space.

\subsection{Construction of the space and its properties}
We begin by defining what we mean by an infinite slab contained in $\mathbb{R}^n$.

\begin{definition}\label{def:infslab}
Let $n \geq 2$. We call an \textit{infinite slab} the set
\[
X := \mathbb{R} \times [0,1]^{n-1} \subset \mathbb{R}^n,
\]
equipped with the Euclidean distance $\rho(x,y) := \|x - y\|$ and with the $n$-dimensional Lebesgue measure $m_n$ restricted to $X$.
\end{definition}

Below, we present two figures illustrating the construction of this space in the cases $n = 2, 3$.

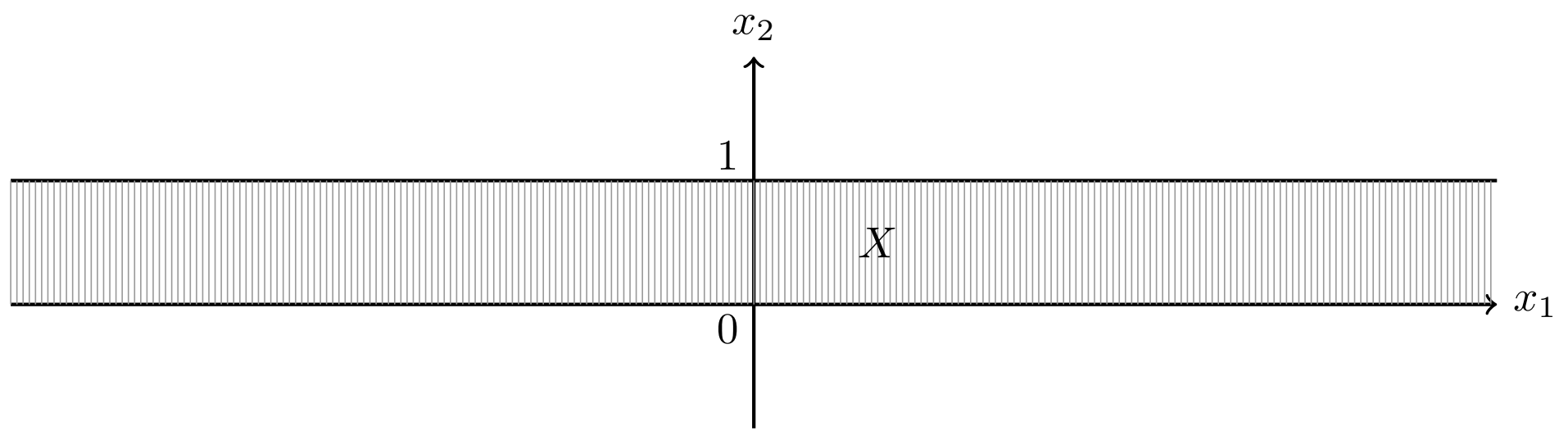
\begin{figure}[H]
\centering
\begin{tikzpicture}[scale=1]

\draw[thick,->] (-6,0) -- (6,0) node[right] {$x_1$};
\draw[thick,] (-6,1) -- (6,1);
\draw[thick,->] (0,-1) -- (0,2) node[above] {$x_2$};

\node[left] at (0,-0.2) {$0$};
\node[left] at (0,1.2) {$1$};

\begin{scope}
    \draw (-4,0) -- (4,0);
    \draw (-4,1) -- (4,1);

    \foreach \t in {-6,-5.95,...,6}
        \draw[black!40, line width=0.3pt] (\t,0) -- (\t,1 );
\end{scope}

\node at (1,0.5) {$X$};

\end{tikzpicture}
\caption{The infinite slab in $\R^2$.}
\end{figure}

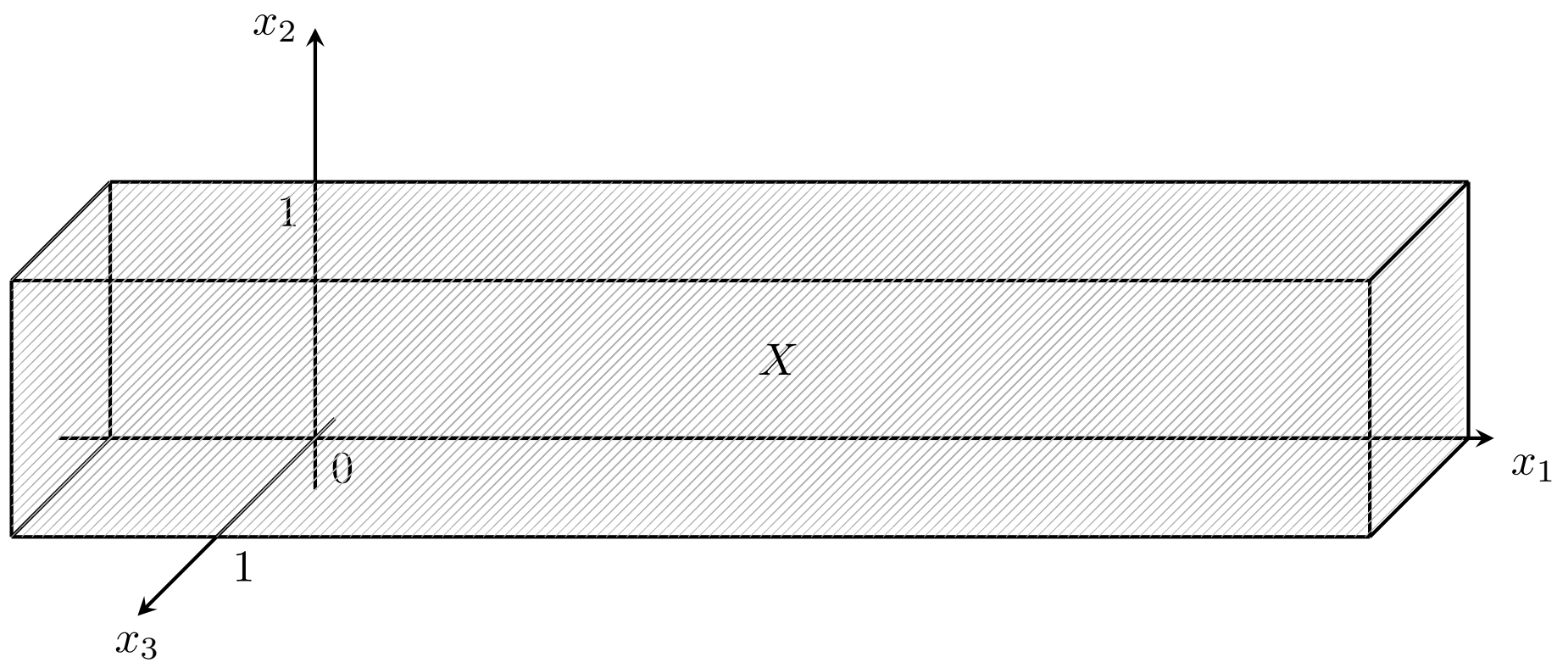
\begin{figure}[H]
\centering
\begin{tikzpicture}[scale=2,>=stealth]

\definecolor{hatchcolor}{RGB}{180,180,180}

\draw[thick,->] (-1,0,0) -- (4.6,0,0) node[below right] {$x_1$};
\draw[thick,->] (0,-0.2,0) -- (0,1.6,0) node[left] {$x_2$};
\draw[thick,->] (0,0,-0.2) -- (0,0,1.8) node[below] {$x_3$};

\node[below right] at (0,0,0) {$0$};
\node[below left] at (0,1,0) {$1$};
\node[below right] at (0,0,1) {$1$};

\def\xmin{-0.8} 
\def\xmax{4.5} 

\draw[thick,black] (\xmin,0,0) -- (\xmax,0,0);
\draw[thick,black] (\xmin,0,1) -- (\xmax,0,1);
\draw[thick,black] (\xmin,0,0) -- (\xmin,0,1);
\draw[thick,black] (\xmax,0,0) -- (\xmax,0,1);

\draw[thick,black] (\xmin,1,0) -- (\xmax,1,0);
\draw[thick,black] (\xmin,1,1) -- (\xmax,1,1);
\draw[thick,black] (\xmin,1,0) -- (\xmin,1,1);
\draw[thick,black] (\xmax,1,0) -- (\xmax,1,1);

\draw[thick,black] (\xmin,0,0) -- (\xmin,1,0);
\draw[thick,black] (\xmax,0,0) -- (\xmax,1,0);
\draw[thick,black] (\xmin,0,1) -- (\xmin,1,1);
\draw[thick,black] (\xmax,0,1) -- (\xmax,1,1);

\foreach \x in {-0.8,-0.76,...,4.5} {
    \foreach \y in {0,0.04,...,1} {
        \draw[hatchcolor, line width=0.2pt] (\x,\y,0) -- (\x,\y,1);
    }
}

\node at (2,0.5,0.5) {$X$};

\end{tikzpicture}
\caption{The infinite slab in $\R^3$.}
\end{figure}

From now on, throughout the rest of this appendix, we fix the letter $X$ to denote this infinite slab in $\mathbb{R}^n$, endowed with the appropriate distance and measure. We define
\[
V(x,r) := m_n\bigl(B(x,r)\cap X\bigr), \quad x \in X,\ r>0,
\]
and
\[
V(x,y) := V\bigl(x,\|x-y\|\bigr), \quad x,y \in X.
\]

Intuitively, what makes this space useful for our purposes is the fact that, if we zoom into it, we obtain an $n$-dimensional object. On the other hand, if we zoom out and observe this set from afar, it appears as a one-dimensional object. We now translate this observation into rigorous terms.

\begin{lemma}\label{lem:slabfund}
    Let $x\in X$. Then
    \[
    V(x,r)\approx\begin{cases}
        r^n, \quad &r\leq1,\\
        r,\quad &r\geq1.
    \end{cases}
    \]
\end{lemma}

The proof of this fact is rather straightforward and is therefore left to the reader.

Next, we prove the doubling property and the positive lower and upper dimensions of this space. 

\begin{lemma}
    The infinite slab $X$, as in Definition~\ref{def:infslab}, is a doubling space
    with lower dimension $d=1$ and upper dimension $D=n$.
\end{lemma}
\begin{proof}
We prove the properties separately.
    \begin{itemize}[label={}, left=-4.5pt]
\item \textbf{Doubling property.}\\
Let $x\in X$, $r>0$. Then, we essentially have two cases:
\begin{enumerate}[label=(\alph*)]
    \item If $r\leq1$, then by Lemma~\ref{lem:slabfund} we have
    \[
    V(x,2r)\leq m_n\bigl(B(x,2r)\bigr)\approx 2^n r^n\approx V(x,r).  
    \]
    \item If $r\geq1$, then, again by Lemma~\ref{lem:slabfund}, we have
    \[
    V(x,2r)\approx 2r \approx V(x,r).  
    \]
\end{enumerate}
\item \textbf{Dimensions.}\\
We prove explicitly only the lower dimension, since the upper dimension can be established analogously.

Let $x\in X$, $0<r\leq R$. Then, keeping in mind Lemma~\ref{lem:slabfund}, we have that
\[
\frac{V(x,R)}{V(x,r)}\approx\begin{cases}
    \displaystyle\frac{R}{r},\qquad &R\geq r\geq1;\\
    \\
    \displaystyle\frac{R}{r^n}= \frac{1}{r^{n-1}}\frac{R}{r}\geq \frac{R}{r}, \qquad &R\geq 1\geq r;\\
    \\
    \displaystyle\frac{R^n}{r^n}= \Big(\frac{R}{r}\Big)^{n-1}\frac{R}{r}\geq \frac{R}{r}, \qquad &1\geq R\geq r.
\end{cases}
\]

Thus, the proof is complete, since in each case we get that 
\[
\frac{V(x,R)}{V(x,r)}\gtrsim \frac{R}{r}.
\]
\qedhere
\end{itemize}
\end{proof}

\subsection{Metric vs. volumic fractional integrals}
We now continue the discussion in this appendix by showing that, in this space, the metric
definition of fractional integrals cannot be regarded as a special case of the volumic one, for any
choice of parameters. This fact is a consequence of the following result.

\begin{lemma}
    There do not exist constants $c\in(0,\infty)$ and $\alpha\in (0,\infty)$ such that
    \[
    \|x-y\|^\alpha\leq c V(x,y)
    \]
    for all $x,y\in X$.
\end{lemma}
\begin{proof}
    Assume for contradiction that there exist $c\in(0,\infty)$ and $\alpha\in (0,\infty)$ such that
    \begin{equation}\label{eq:ProofMetVolSlab1}
        \|x-y\|^\alpha\leq c V(x,y)
    \end{equation}
    for all $x,y\in X$. Such a uniform estimate, due to the structure of this space, implies two necessary conditions, formulated as follows.
    \begin{itemize}
        \item If $\|x-y\|\geq1$, then by \eqref{eq:ProofMetVolSlab1} and Lemma~\ref{lem:slabfund} we get
        \[
        \|x-y\|^\alpha\lesssim \|x-y\|.
        \]
        Sending $\|x-y\|$ to infinity, we obtain that necessarily \[\alpha \le 1.\]
        \item If $\|x-y\|\leq1$, then by \eqref{eq:ProofMetVolSlab1} and Lemma~\ref{lem:slabfund} we get
        \[
        \|x-y\|^\alpha\lesssim \|x-y\|^n.
        \]
        Sending $\|x-y\|$ to $0$, we obtain that necessarily \[\alpha \ge n.\]
    \end{itemize}
    The two necessary conditions obtained are clearly contradictory, and thus this concludes the proof.
\end{proof}

    \subsection{Failure of \texorpdfstring{$L^p$--$L^q$}{Lp--Lq} boundedness for metric fractional integrals}
In this final section, we show that the metric fractional integrals are not bounded from $L^p$ to $L^q$. In order to prove such unboundedness, it will be sufficient to obtain two necessary conditions, which will turn out to be contradictory. These conditions are established in the following two lemmas.

\begin{lemma}\label{lem:NecCondSlab1}
    Let $\alpha\in(0,\infty)$ and let $I_\alpha$ be the operator with integral representation
    \[
    I_\alpha f(x):=\int_X \frac{\|x-y\|^\alpha}{V(x,y)}f(y) \ud y, 
    \]
    where $\ud y$ denotes integration with respect to the $n$-dimensional Lebesgue measure restricted to $X$. Moreover, let $p,q\in(1,\infty)$. Then, a necessary condition for $I_\alpha$ to be bounded from $L^p(X)$ to $L^q(X)$ is that
    \[
    \frac1q\geq \frac{1}{p}-\frac{\alpha}{n}.
    \]
\end{lemma}
\begin{proof}
    Assume that $I_\alpha$ is bounded from $L^p(X)$ to $L^q(X)$, that is, 
    \begin{equation*}
        \Norm{I_\alpha f}{L^q(X)}\lesssim \Norm{f}{L^p(X)}, \quad \forall f\in L^p(X).
    \end{equation*}
    Then, in particular, such an estimate holds for functions $f \in L^p(X)$ whose support is contained in $[0,1]^n$. Next, it is clearly true that
    \[
    \Norm{\mathbf{1}_{[0,1]^n}I_\alpha f}{L^q(X)}\lesssim \Norm{I_\alpha f}{L^q(X)},
    \]
    where $\mathbf{1}_{[0,1]^n}$ denotes the characteristic function of the set $[0,1]^n$. Therefore, so far we have obtained the estimate
    \[
    \Norm{\mathbf{1}_{[0,1]^n}I_\alpha f}{L^q(X)}\lesssim \Norm{f}{L^p(X)}, \quad \forall f\in L^p(X), \, \supp f \subset [0,1]^n.
    \]
    Now, for such functions $f$ and for $\lambda>0$,  consider the dilated function $f_\lambda(x) := f(\lambda x)$. To give meaning to this definition for every $\lambda>0$, note that we are implicitly extending $f$ to be $0$ outside of $X$. Moreover, we observe that
\[
f_\lambda(x) \neq 0 \ \Rightarrow\ \lambda x \in [0,1]^n \ \Leftrightarrow\ x \in \Big[0,\frac{1}{\lambda}\Big]^n,
\]
and therefore, if $\lambda \ge 1$, we have
\[
\supp f_\lambda \subset \Big[0,\frac{1}{\lambda}\Big]^n \subset [0,1]^n.
\]
With this restriction on $\lambda$, we thus obtain that
\[
    \Norm{\mathbf{1}_{[0,1]^n}I_\alpha f_\lambda}{L^q(X)}\lesssim \Norm{f_\lambda}{L^p(X)}.
    \]
Now, by Lemma~\ref{lem:slabfund},
\begin{equation*}
\begin{split}
\Norm{\mathbf{1}_{[0,1]^n}I_\alpha f_\lambda}{L^q(X)}
&= \Bigg( \int_{\R^n} \mathbf{1}_{[0,1]^n}(x) \Big| \int_{\R^n} \frac{f(\lambda y)}{\|x-y\|^{n-\alpha}} \, \ud y \Big|^q \, \ud x \Bigg)^{1/q} \\
&= \Bigg( \int_{\R^n} \mathbf{1}_{[0,1]^n}\Big(\frac{x}{\lambda}\Big) \Big| \int_{\R^n} \frac{f(y)}{\bigl\|\frac{x}{\lambda}-\frac{y}{\lambda}\bigr\|^{n-\alpha}} \frac{\ud y}{\lambda^n} \Big|^q \frac{\ud x}{\lambda^n} \Bigg)^{1/q} \\
&= \lambda^{-\alpha - n/q} \Bigg( \int_{\R^n} \underbrace{\mathbf{1}_{[0,\lambda]^n}(x)}_{\geq \mathbf{1}_{[0,1]^n}(x)} |I_\alpha f(x)|^q \, \ud x \Bigg)^{1/q} \\
&\geq \lambda^{-\alpha - n/q} \Norm{\mathbf{1}_{[0,1]^n} I_\alpha f}{L^q(X)}.
\end{split}
\end{equation*}
On the other hand,
\[
\Norm{f_\lambda}{L^p(X)}=\Bigg(\int_{\R^n}\abs{f(\lambda x)}^p\ud x\Bigg)^{1/p}=\lambda^{-n/p}\Norm{f}{L^p(X)}.
\]
Thus, altogether, we may conclude that
\[
\lambda^{-\alpha-n/q}\Norm{\mathbf{1}_{[0,1]^n} I_\alpha f}{L^q(X)}\lesssim \lambda^{-n/p}\Norm{f}{L^p(X)},
\]
uniformly in $\lambda\geq 1$ and $f \in L^p(X)$ whose support is contained in $[0,1]^n$. In order for this to hold, particularly as $\lambda$ tends to infinity, it is necessary that
\[
\alpha+\frac nq\geq \frac np \quad \Leftrightarrow \quad \frac1q\geq \frac1p-\frac \alpha n.
\]
This concludes the proof.
\end{proof}

\begin{lemma}\label{lem:NecCondSlab2}
    Let $\alpha\in(0,\infty)$ and let $I_\alpha$ be the operator with integral representation
    \[
    I_\alpha f(x):=\int_X \frac{\|x-y\|^\alpha}{V(x,y)}f(y) \ud y, 
    \]
    where $\ud y$ denotes integration with respect to the $n$-dimensional Lebesgue measure restricted to $X$. Moreover, let $p,q\in(1,\infty)$. Then, a necessary condition for $I_\alpha$ to be bounded from $L^p(X)$ to $L^q(X)$ is that
    \[
    \frac1q\leq \frac{1}{p}-\alpha.
    \]
\end{lemma}

\begin{proof}
    Assume that $I_\alpha$ is bounded from $L^p(X)$ to $L^q(X)$, that is,
    \begin{equation}\label{eq:NecCondSlab2Proof1}
        \Norm{I_\alpha g}{L^q(X)}\lesssim \Norm{g}{L^p(X)}, \quad \forall g\in L^p(X).
    \end{equation}
    We now aim to derive a necessary condition from this relation by constructing a specific function, given by a weighted sum of smooth bump functions, on which we impose this condition.
    
    Take $v:=\left(0,\frac12,\dots,\frac12\right)\in X$, and let $\Phi\in C^\infty_0\bigl(B(v,1/4)\bigr)$. Next, let us take a sequence $(a_k)_{k\in\Z}$, with $a_k\geq0$, to be determined later, and let us set
    \[
    f(y):=\sum_{k\in \Z} a_k \Phi\bigl(y-(k,0,\dots,0)\bigr), \quad y\in X.
    \]
    Now, we examine how the boundedness assumption on $I_\alpha$ translates when applied to this function $f$. For any $l\in\Z$ and $x\in B\bigl((l,\frac12,\dots,\frac12),\frac14\bigr)\subset X$, we have
    \begin{equation*}
        \begin{split}
            I_\alpha f(x)&=\sum_{k\in \Z} a_k\int_X \frac{\|x-y\|^\alpha}{V(x,y)} \Phi\bigl(y-(k,0,\dots,0)\bigr)\ud y\\
            &\geq \sum_{\substack{k\in \Z\\k\neq l}} a_k\int_X \frac{\|x-y\|^\alpha}{V(x,y)} \Phi\bigl(y-(k,0,\dots,0)\bigr)\ud y\\
            &\approx\sum_{\substack{k\in \Z\\k\neq l}} a_k\abs{l-k}^{\alpha-1},
        \end{split}
    \end{equation*}
    where, in the last step, we used the compact support of $\Phi$, the doubling property, and Lemma~\ref{lem:slabfund} to estimate in each term of the sum
    \begin{equation*}
        \begin{cases}
            \|x-y\|^\alpha\approx |l-k|^\alpha,\\
            V(x,y)\approx V\Bigg(\Big(l,\frac12,\dots,\frac12\Big),\Big(k,\frac12,\dots,\frac12\Big) \Bigg)=V\Bigg(\Big(l,\frac12,\dots,\frac12\Big),\abs{l-k} \Bigg)\approx \abs{l-k}.
        \end{cases}
    \end{equation*}
    Thus, 
    \begin{equation*}
        \Norm{I_\alpha f}{L^q(X)} \geq \Big(\sum_{l\in\Z}\int_{B\bigl((l,\frac12,\dots,\frac12),\frac14\bigr)}\abs{I_\alpha f(x)}^q\ud x\Big)^{\frac1q} \gtrsim \Big(\sum_{l\in\Z}\Big(\sum_{\substack{k\in \Z\\k\neq l}} a_k\abs{l-k}^{\alpha-1}\Big)^q\Big)^{\frac1q}.
    \end{equation*}
    On the other hand,
    \begin{equation*}
        \Norm{f}{L^p(X)}=\Big(\int_X \sum_{k\in\Z}a_k^p\Phi\bigl(x-(k,0,\dots,0)\bigr)^p \ud x\Big)^{\frac1p}\approx \left(\sum_{k\in\Z}a_k^p\right)^{\frac1p}.
    \end{equation*}
    Therefore, combining these two estimates with \eqref{eq:NecCondSlab2Proof1}, we obtain
    \begin{equation}\label{eq:NecCondSlab2Proof2}
        \Big(\sum_{l\in\Z}\Big(\sum_{\substack{k\in \Z\\k\neq l}} a_k\abs{l-k}^{\alpha-1}\Big)^q\Big)^{\frac1q}\lesssim \left(\sum_{k\in\Z}a_k^p\right)^{\frac1p}.
    \end{equation}
    At this stage, we choose a suitable sequence $(a_k)$. For $N\geq1$, take
    \[
    a_k:=\abs{k}^{-\frac1p}\mathbf{1}_{\{0<\abs{k}\leq N\}}.
    \]
    Then, we have that
    \begin{equation}\label{eq:NecCondSlab2Proof3}
    \left(\sum_{k\in\Z}a_k^p\right)^{\frac1p}\approx \left(\sum_{k=1}^N k^{-1}\right)^{\frac1p}\sim \bigl(\log N\bigr)^{\frac1p}, \quad \text{as $N\to\infty$},
    \end{equation}
    where the symbol $\sim$ denotes asymptotic equivalence. On the other hand, for $l\in\Z$ such that $\abs{l}\leq N$, we have
    \begin{equation*}
        \sum_{k\neq l}a_k\abs{l-k}^{\alpha-1}=
        \sum_{\substack{k\neq l\\0<\abs{k}\leq N}}\abs{l-k}^{\alpha-1}\abs{k}^{-\frac1p}\gtrsim \sum_{0<\abs{k}\leq \frac12 \abs{l}}\abs{l}^{\alpha-1}\abs{k}^{-\frac1p}\gtrsim \abs{l}^{\alpha-\frac1p},
    \end{equation*}
    where, in the last step, we have used the fact that
    \[
    \sum_{k=1}^M k^{-\frac1p}\sim \frac{M^{1-\frac1p}}{1-\frac1p}=\frac{p}{p-1}M^{1-\frac1p}, \quad \text{as }M\to\infty.
    \]
    Hence,
    \[
    \Big(\sum_{l\in\Z}\Big(\sum_{\substack{k\in \Z\\k\neq l}} a_k\abs{l-k}^{\alpha-1}\Big)^q\Big)^{\frac1q}\gtrsim \Big(\sum_{0<\abs{l}<N}\bigl(\abs{l}^{\alpha-\frac1p}\bigr)^q\Big)^{\frac1q},
    \]
    and, if
    \[
    \Big(\alpha-\frac1p\Big)q>-1 \quad \Leftrightarrow\quad \frac1q+\alpha-\frac1p>0,
    \]
    then
    \[
    \Big(\sum_{0<\abs{l}<N}\bigl(\abs{l}^{\alpha-\frac1p}\bigr)^q\Big)^{\frac1q}\sim c \Big(N^{1+\big(\alpha-\frac1p\big)q}\Big)^{\frac1q}\approx N^{\frac1q+\alpha-\frac1p}, \quad \text{as }N\to\infty.
    \]
    By \eqref{eq:NecCondSlab2Proof2} and \eqref{eq:NecCondSlab2Proof3}, this is a contradiction. Hence, we may conclude that it must be
    \[
    \Big(\alpha-\frac1p\Big)q\leq -1 \quad \Leftrightarrow \quad \frac1q\leq \frac{1}{p}-\alpha,
    \]
    which is what we wanted to prove.
\end{proof}

At this point, we have everything we need to conclude that the fractional integral operators on this space are not bounded from $L^p$ to $L^q$.

\begin{lemma}\label{lem:UnbddFracMetSlab}
    Let $\alpha\in(0,\infty)$ and let $I_\alpha$ be the operator with integral representation
    \[
    I_\alpha f(x):=\int_X \frac{\|x-y\|^\alpha}{V(x,y)}f(y) \ud y.
    \]
   Moreover, let $p,q\in(1,\infty)$. Then $I_\alpha$ is not bounded from $L^p(X)$ to $L^q(X)$.
\end{lemma}

\begin{proof}
    Thanks to the results obtained in Lemmas \ref{lem:NecCondSlab1} and \ref{lem:NecCondSlab2}, necessary conditions for $I_\alpha$ to be bounded from $L^p(X)$ to $L^q(X)$ are the following:
\begin{equation*}
    \begin{cases}
        \displaystyle\frac1q\geq \frac{1}{p}-\frac{\alpha}{n} & \text{(Lemma~\ref{lem:NecCondSlab1})},\\
        \\
        \displaystyle\frac1q\leq \frac{1}{p}-\alpha& \text{(Lemma~\ref{lem:NecCondSlab2})}.
    \end{cases}
\end{equation*}
By combining these two conditions, we obtain as a consequence that
\[
\frac{1}{p}-\frac{\alpha}{n}\leq \frac{1}{p}-\alpha \quad \Leftrightarrow \quad \frac{\alpha}{n}\geq \alpha.
\]
Since $\alpha>0$ and $n\geq 2$, this last inequality is clearly false. Consequently, $I_\alpha$ is not bounded from $L^p(X)$ to $L^q(X)$.
\end{proof}

The result just obtained is precisely what we had set as our initial goal, namely, to show that metric-type fractional integral operators are not, in general, bounded from $L^p$ to $L^q$, not even on spaces of homogeneous type with positive lower and upper dimensions. This therefore shows the necessity, for this class of operators, to consider a more unconventional notion of boundedness, such as the one established in Theorem~\ref{thm:bddfracmet}.

\section{A Cantor--type set}\label{app:Cantor}


In Appendix~\ref{app:Slab} we observed, through a concrete example, that even for relatively simple spaces possessing lower and upper dimensional properties, one cannot obtain classical $L^p$--$L^q$ boundedness for metric fractional integrals.

In this second appendix, our goal is to illustrate what can truly go wrong for metric fractional integrals when we consider the even more general setting in which no positive lower dimension is assumed. To this end, we construct a Cantor--type subset of the real line exhibiting several noteworthy properties, summarized below.

\begin{FeaturesCantor}
The Cantor--type set, constructed in Subsection~\ref{subsec:ConstrCantor}, has the following properties:
    \begin{itemize}
\item It is a space of homogeneous type.
\item It does not admit a lower dimension $d>0$.
\item In this space, the metric definition of fractional integrals cannot be regarded as a special case of the volumic one, for any choice of parameters.
\item In this space, for every $p \in (0,\infty)$ a generic metric fractional integral may map functions in $L^p$ to the function identically equal to $+\infty$.
\end{itemize}
\end{FeaturesCantor}

The last property is particularly significant. Indeed, this implies that without additional assumptions on the underlying space of homogeneous type, no boundedness from $L^p$ to any conventional function space can be obtained for these operators. In particular, this establishes the minimality of the positive lower dimension assumption in Theorem~\ref{thm:bddfracmet}.

In what follows, we study in detail the construction of this space and show the features listed above.

\subsection{Construction}\label{subsec:ConstrCantor}

For this construction, let us consider an increasing sequence $a=(a_n)_{n=1}^\infty$ such that $a_n\in(0,\infty)$ for all $n\in\{1,2,\dots\}$. Moreover, we denote the shift of a set $A \subseteq \R$ by $h \in \R$ as
\[
A+h:=\{a+h\mid a\in A\}.
\]
Set $I_{0,1}:=[1,2]$. Inductively, suppose $I_{n,1}$ has been constructed and put
\[
\Delta_n:=\operatorname{diam}(I_{n,1})+a_{n+1}=\bigl(\sup I_{n,1}-\inf I_{n,1}\bigr)+a_{n+1}.
\]
Define
\[
I_{n+1,1}:=I_{n,1}\cup\bigl(I_{n,1}+\Delta_n\bigr).
\]
Then the left endpoint of the translated copy $I_{n,1}+\Delta_n$ equals $\sup I_{n,1}+a_{n+1}$,
so the gap between the rightmost point of $I_{n,1}$ and the leftmost point of the copy is exactly $a_{n+1}$. In other words,
\begin{equation}\label{eq:distCantor}
    \operatorname{dist}(I_{n,1},I_{n,1}+\Delta_n)=a_{n+1}, \quad n\geq0.
\end{equation}
Finally, set
\[
X_a:=\bigcup_{n\ge0} I_{n,1}.
\]

Since we will later need to identify precisely the various components of this recursive structure,  
so as to encode both their position within the structure and their hierarchical level,  
we adopt the notation summarized as follows:
\begin{align*}
    I_{0,1}&:=[1,2], \\
    I_{1,1}&:=I_{0,1}\cup(I_{0,1}+\Delta_0)\\
    &=:I_{0,1}\cup I_{0,2}\\
    I_{2,1}&:=I_{1,1}\cup(I_{1,1}+\Delta_1)\\
    &=:I_{1,1}\cup I_{1,2}\\
    &=:I_{0,1}\cup I_{0,2}\cup I_{0,3}\cup I_{0,4}\\[0.5em]
    \dots \\[0.5em]
   \text{For $n\geq3$:}\quad I_{n,1}&:=I_{n-1,1}\cup(I_{n-1,1}+\Delta_{n-1})\\
    &=:I_{n-1,1}\cup I_{n-1,2}\\
    &=:I_{n-2,1}\cup I_{n-2,2}\cup I_{n-2,3}\cup I_{n-2,4}\\
    &\dots\\
    &=: \bigcup_{i=1}^{2^{n-k}} I_{k,i}\qquad(0\leq k \leq n)\\
    &\dots\\
    &=: \bigcup_{i=1}^{2^{n}} I_{0,i}.
\end{align*}
With this notation, the components of the construction $\{I_{i,j}\}_{i,j=1}^\infty$ are identified by two indices: the first encodes the hierarchical level, while the second specifies its location. Note also that, for a fixed $n$, the family $\{I_{n,j}\}_{j=1}^\infty$ consists of infinitely many translated copies of $I_{n,1}$, and that, in light of \eqref{eq:distCantor}, we get
\[
\operatorname{dist}(I_{n,1},I_{n,2})=a_{n+1}, \quad n\geq0.
\]

In the hope that this helps the reader better understand the construction of this set, we include two figures illustrating the first three steps and the $n$-th step in the construction of $X_{a}$. In such figures, the dotted line represents the real axis, while the solid black segments denote the building blocks of the Cantor-type set. The hierarchical structure is emphasized by the braces below the figures together with the labels $I_{i,j}$. The braces above the figures indicate the distances between the components, encoded by the sequence $(a_n)_{n\ge 1}$.
\vspace{3em}
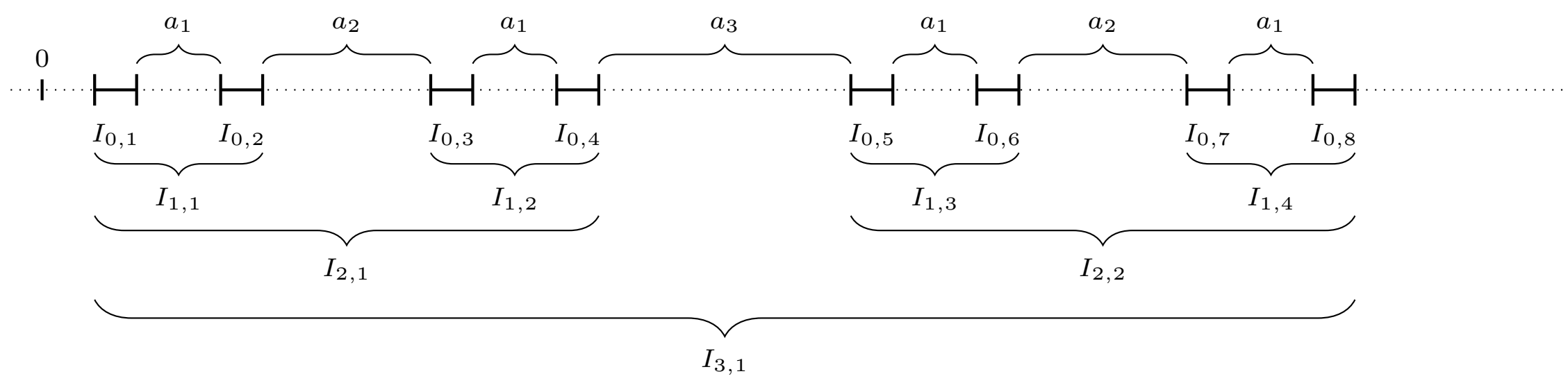
\begin{figure}[H]
  \centering

\begin{tikzpicture}[scale=1.0,decoration={brace}]
\def\unitlen{0.4}
\def\scalefactor{0.4}
\def\aone{2*\scalefactor}
\def\atwo{4*\scalefactor}
\def\athree{6*\scalefactor}
\def\afour{8*\scalefactor}

\draw[dotted] (-0.8,0) -- (14,0);
\draw[thick] (-0.5,-0.1)--(-0.5,0.1);
\node[above] at (-0.5,0.1){\scriptsize $0$};

\newcommand{\interval}[3]{
  \draw[thick] (#1,0)--(#1+#2,0);
  \draw[thick] (#1,0.15)--(#1,-0.15);
  \draw[thick] (#1+#2,0.15)--(#1+#2,-0.15);
  \node[below] at (#1+0.5*#2,-0.2){\scriptsize #3};
}

\interval{0}{\unitlen}{$I_{0,1}$}
\interval{\unitlen+\aone}{\unitlen}{$I_{0,2}$}

\draw[decorate,decoration={brace,amplitude=5pt}]
  (\unitlen,0.25) -- (\unitlen+\aone,0.25)
  node[midway,above=5pt]{\scriptsize $a_1$};

\draw[decorate,decoration={brace,amplitude=6pt,mirror}]
  (0,-0.6)--(\unitlen+\aone+\unitlen,-0.6)
  node[midway,below=6pt]{\scriptsize $I_{1,1}$};

\def\offsetA{\unitlen+\aone+\unitlen+\atwo}
\interval{\offsetA}{\unitlen}{$I_{0,3}$}
\interval{\offsetA+\unitlen+\aone}{\unitlen}{$I_{0,4}$}
\draw[decorate,decoration={brace,amplitude=5pt}]
  (\unitlen+\aone+\unitlen,0.25) -- (\offsetA,0.25)
  node[midway,above=5pt]{\scriptsize $a_2$};

\draw[decorate,decoration={brace,amplitude=5pt}]
  (\unitlen+\aone+\unitlen+\atwo+\unitlen,0.25) -- (\unitlen+\aone+\unitlen+\atwo+\unitlen+\aone,0.25)
  node[midway,above=5pt]{\scriptsize $a_1$};  
\draw[decorate,decoration={brace,amplitude=6pt,mirror}]
  (\offsetA,-0.6)--(\offsetA+\unitlen+\aone+\unitlen,-0.6)
  node[midway,below=6pt]{\scriptsize $I_{1,2}$};

\draw[decorate,decoration={brace,amplitude=8pt,mirror}]
  (0,-1.2)--(\offsetA+\unitlen+\aone+\unitlen,-1.2)
  node[midway,below=8pt]{\scriptsize $I_{2,1}$};

\def\offsetB{\offsetA+\unitlen+\aone+\unitlen+\athree}
\interval{\offsetB}{\unitlen}{$I_{0,5}$}
\interval{\offsetB+\unitlen+\aone}{\unitlen}{$I_{0,6}$}
\draw[decorate,decoration={brace,amplitude=5pt}]
  (\offsetA+\unitlen+\aone+\unitlen,0.25)--(\offsetB,0.25)
  node[midway,above=5pt]{\scriptsize $a_3$};

\draw[decorate,decoration={brace,amplitude=5pt}](\offsetA+\unitlen+\aone+\unitlen+\athree+\unitlen,0.25) -- (\offsetA+\unitlen+\aone+\unitlen+\athree+\unitlen+\aone,0.25)node[midway,above=5pt]{\scriptsize $a_1$}; 
  
\draw[decorate,decoration={brace,amplitude=6pt,mirror}]
  (\offsetB,-0.6)--(\offsetB+\unitlen+\aone+\unitlen,-0.6)
  node[midway,below=6pt]{\scriptsize $I_{1,3}$};

\def\offsetC{\offsetB+\unitlen+\aone+\unitlen+\atwo}
\interval{\offsetC}{\unitlen}{$I_{0,7}$}
\interval{\offsetC+\unitlen+\aone}{\unitlen}{$I_{0,8}$}
\draw[decorate,decoration={brace,amplitude=5pt}]
  (\offsetB+\unitlen+\aone+\unitlen,0.25)--(\offsetC,0.25)
  node[midway,above=5pt]{\scriptsize $a_2$};

\draw[decorate,decoration={brace,amplitude=5pt}](\offsetB+\unitlen+\aone+\unitlen+\atwo+\unitlen,0.25) -- (\offsetB+\unitlen+\aone+\unitlen+\atwo+\unitlen+\aone,0.25)node[midway,above=5pt]{\scriptsize $a_1$}; 
  
\draw[decorate,decoration={brace,amplitude=6pt,mirror}]
  (\offsetC,-0.6)--(\offsetC+\unitlen+\aone+\unitlen,-0.6)
  node[midway,below=6pt]{\scriptsize $I_{1,4}$};

\draw[decorate,decoration={brace,amplitude=8pt,mirror}]
  (\offsetB,-1.2)--(\offsetC+\unitlen+\aone+\unitlen,-1.2)
  node[midway,below=8pt]{\scriptsize $I_{2,2}$};

\draw[decorate,decoration={brace,amplitude=10pt,mirror}]
  (0,-2.0)--(\offsetC+\unitlen+\aone+\unitlen,-2.0)
  node[midway,below=10pt]{\scriptsize $I_{3,1}$};
\end{tikzpicture}

\caption{The first three steps of the construction of this Cantor-type set.}
  \label{fig:FirstThreeCantor}
\end{figure}

\vspace{3em}

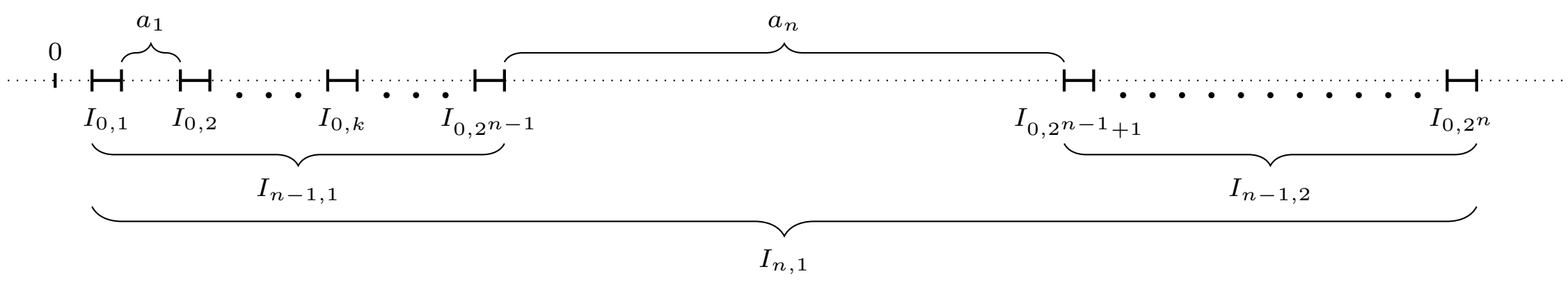
\begin{figure}[H]
  \centering
  \begin{tikzpicture}[scale=0.7,decoration={brace}]
\def\unitlen{0.4}
\def\scalefactor{0.4}
\def\aone{2*\scalefactor}
\def\atwo{4*\scalefactor}
\def\athree{6*\scalefactor}
\def\afour{8*\scalefactor}

\def\blocklen{\unitlen+\afour+\atwo+\unitlen}

\draw[dotted] (-1.1428,0) -- (20,0);
\draw[thick] (-0.5,-0.1)--(-0.5,0.1);
\node[above] at (-0.5,0.1){\scriptsize $0$};

\newcommand{\interval}[3]{
  \draw[thick] (#1,0)--(#1+#2,0);
  \draw[thick] (#1,0.15)--(#1,-0.15);
  \draw[thick] (#1+#2,0.15)--(#1+#2,-0.15);
  \node[below] at (#1+0.5*#2,-0.2){\scriptsize #3};
}

\interval{0}{\unitlen}{$I_{0,1}$}
\interval{\unitlen+\aone}{\unitlen}{$I_{0,2}$}

\draw[fill] (\unitlen+\aone+0.8,-0.2) circle (0.9pt);
\draw[fill] (\unitlen+\aone+1.2,-0.2) circle (0.9pt);
\draw[fill] (\unitlen+\aone+1.6,-0.2) circle (0.9pt);

\draw[decorate,decoration={brace,amplitude=5pt}]
  (\unitlen,0.25) -- (\unitlen+\aone,0.25)
  node[midway,above=5pt]{\scriptsize $a_1$};

\draw[decorate,decoration={brace,amplitude=6pt,mirror}]
  (0,-0.8571)--(\blocklen,-0.8571)
  node[midway,below=6pt]{\scriptsize $I_{n-1,1}$};

\def\offsetA{\unitlen+\unitlen+\atwo}
\interval{\offsetA+\aone}{\unitlen}{$I_{0,k}$}
\interval{\offsetA+\unitlen+\aone+\atwo}{\unitlen}{$I_{0,2^{n-1}}$}

\draw[fill] (\offsetA+\aone+0.8,-0.2) circle (0.9pt);
\draw[fill] (\offsetA+\aone+1.2,-0.2) circle (0.9pt);
\draw[fill] (\offsetA+\aone+1.6,-0.2) circle (0.9pt);

\def\offsetB{\offsetA+\offsetA+\offsetA+\aone+\athree+\unitlen+\athree}
\interval{\offsetB}{\unitlen}{$I_{0,2^{n-1}+1}$}
\interval{\offsetB+\blocklen-\unitlen}{\unitlen}{$I_{0,2^n}$}

\draw[fill] (\offsetB+\unitlen+0.4,-0.2) circle (0.9pt);
\draw[fill] (\offsetB+\unitlen+0.8,-0.2) circle (0.9pt);
\draw[fill] (\offsetB+\unitlen+1.2,-0.2) circle (0.9pt);
\draw[fill] (\offsetB+\unitlen+1.6,-0.2) circle (0.9pt);
\draw[fill] (\offsetB+\unitlen+2.0,-0.2) circle (0.9pt);
\draw[fill] (\offsetB+\unitlen+2.4,-0.2) circle (0.9pt);
\draw[fill] (\offsetB+\unitlen+2.8,-0.2) circle (0.9pt);
\draw[fill] (\offsetB+\unitlen+3.2,-0.2) circle (0.9pt);
\draw[fill] (\offsetB+\unitlen+3.6,-0.2) circle (0.9pt);
\draw[fill] (\offsetB+\unitlen+4.0,-0.2) circle (0.9pt);
\draw[fill] (\offsetB+\unitlen+4.4,-0.2) circle (0.9pt);

\draw[decorate,decoration={brace,amplitude=5pt}]
  (\offsetA+\unitlen+\athree+\unitlen,0.25)--(\offsetB,0.25)
  node[midway,above=5pt]{\scriptsize $a_n$};

\draw[decorate,decoration={brace,amplitude=6pt,mirror}]
  (\offsetB,-0.8571)--(\offsetB+\blocklen,-0.8571)
  node[midway,below=6pt]{\scriptsize $I_{n-1,2}$};

\draw[decorate,decoration={brace,amplitude=8pt,mirror}]
  (0,-1.7142)--(\offsetB+\blocklen,-1.7142)
  node[midway,below=8pt]{\scriptsize $I_{n,1}$};

\end{tikzpicture}
  \caption{The $n$-th step of the construction of the Cantor-type set.}
\end{figure}

Now that we have constructed the set $X_{a}$, we equip it with the $1$-dimensional Lebesgue measure $\mu := m_1\vert_{X_{a}}$ restricted to $X_{a}$, and we take the distance $\rho(x,y) := |x-y|$. 

\subsection{Homogeneous type property}

The idea behind this construction is to choose a sequence $a$ whose growth is sufficiently fast so that the growth of the measure of balls is significantly slower than the corresponding growth of their radius. On the other hand, we want the set to retain the properties of a space of homogeneous type. More precisely, we require that the doubling property is preserved. These features are summarized in the following lemma.

\begin{lemma}\label{lem:Cantorprop}
The space $(X_{a},|\cdot-\cdot|,m_1\vert_{X_{a}})$ just constructed satisfies the following properties:
    \begin{enumerate}[label=(\roman*)]
        \item\label{Cantorpoint:i} For any choice of the sequence $a=(a_n)_{n=1}^\infty$, for every $n\in\{0,1,2,\dots\}$ and $i\in\{1,2,3,\dots\}$, 
        \begin{equation}\label{eq:Cantori}
            m_1(I_{n,i})=2^n.
        \end{equation}
        \item\label{Cantorpoint:ii} For every $n\in\{0,1,2,\dots\}$ and $i\in\{1,2,3,\dots\}$, 
        \begin{equation}\label{eq:Cantorii1}
            \dist(I_{n,i},I_{n,i+1})\geq a_{n+1},
        \end{equation}
        and
        \begin{equation}\label{eq:Cantorii2}
            \operatorname{diam}(I_{n,i})=2^n+\sum_{k=1}^n2^{n-k} a_{k}.
        \end{equation}
        \item\label{Cantorpoint:iii} If 
        \begin{equation}\label{eq:seqgrow}
            a_n>2\Big(2^{n-1}+\sum_{k=1}^{n-1}2^{n-1-k} a_{k}\Big)=2\operatorname{diam}(I_{n-1,1}), \quad n\geq1,
        \end{equation}
        then $m_1\vert_{X_{a}}$ is doubling with constant $2$.
        \item\label{Cantorpoint:iv} If we define
        \begin{equation*}
            a_n:=2^{2^n}, \quad n\geq1,
        \end{equation*}
        then $a=(a_n)_{n=1}^\infty$ satisfies condition \eqref{eq:seqgrow}.
    \end{enumerate}
\end{lemma}
\begin{remark}
    In this lemma, our goal is to identify sufficient conditions under which such a space is doubling while also exhibiting the useful property of slow growth of the measures of balls, achieved through a rapid growth of $a$, which we will use later. However, we are not claiming (and it does not seem relevant for our purposes to investigate) that these conditions (i.e., formula \eqref{eq:seqgrow}) are necessary.
\end{remark}

\begin{proof}[Proof of Lemma \ref{lem:Cantorprop}]
We prove each point of the lemma separately.
    \begin{itemize}[label={}, left=-4.5pt]
        \item \textbf{Point~\ref{Cantorpoint:i}.}\\
        Fix $n \in \{0,1,2,\dots\}$. Since the Lebesgue measure is translation-invariant, and each $I_{n,i}$, for $i \in \{1,2,3,\dots\}$, is merely a translated copy of $I_{n,1}$, it is sufficient to show formula~\eqref{eq:Cantori} for $i=1$. Using again the translation invariance of $m_1$ and its disjoint additivity, as well as the structure of $I_{n,1}$, we obtain that
        \[
    m_1(I_{n,1})=m_1\left(\bigcup_{i=1}^{2^{n}} I_{0,i}\right)=\sum_{i=1}^{2^{n}} m_1(I_{0,i})=2^nm_1(I_{0,1})=2^n.
        \]
        
        \item \textbf{Point~\ref{Cantorpoint:ii}.}\\
        To prove formula \eqref{eq:Cantorii1}, fix $n\ge 0$ and $i\ge 1$. By construction, the interval $I_{n,i+1}$
is obtained from $I_{n,i}$ by a translation of the form
\[
I_{n,i+1}=I_{n,i}+\operatorname{diam}(I_{n,i})+a_k,
\qquad\text{for some } k\ge n+1,
\]
depending on the position of the two intervals within the construction. Since the sequence $(a_h)_{h\ge 0}$ is increasing, we have
$\operatorname{diam}(I_{n,i})+a_k\ge\operatorname{diam}(I_{n,i})+a_{n+1}$.
Hence,
\[
\dist(I_{n,i},I_{n,i+1})
   =(\operatorname{diam}(I_{n,i})+a_k)-\operatorname{diam}(I_{n,i})\geq a_{n+1},
\]
This proves the first claim.

        We now turn to formula \eqref{eq:Cantorii2}. To this end, first observe that, as in point \ref{Cantorpoint:i}, it is sufficient to prove the claim for $i=1$. We establish this assertion by induction.
\begin{itemize}
    \item \textbf{Base case.} For $n=0$, we have
    \begin{align*}
        &\operatorname{diam}(I_{0,1})=1, 
    \end{align*}
    \item \textbf{Inductive step.} Assume that \eqref{eq:Cantorii2} is valid for $n$. Then,
    \[
    \operatorname{diam}(I_{n+1,1})=2\operatorname{diam}(I_{n,1})+a_{n+1}=2^{n+1}+\sum_{i=1}^n\Big(2^{n+1-i} a_{i}\Big)+a_{n+1}=2^{n+1}+\sum_{i=1}^{n+1}\Big(2^{n+1-i} a_{i}\Big).
    \]
\end{itemize}
        \item \textbf{Point~\ref{Cantorpoint:iii}.}\\
        Let $a$ be a sequence that satisfies \eqref{eq:seqgrow}. We want to prove that there exists a constant $C>0$ such that, for all $x\in X_{a}$ and $r>0$,
        \[
        m_1\bigl(B_{X_{a}}(x,2r)\bigr)\leq Cm_1\bigl(B_{X_{a}}(x,r)\bigr),
        \]
        where we denote
        \[
        B_{X_{a}}(x,r):=(x-r,x+r)\cap X_{a}.
        \]
        Let us take $x\in X_{a}$. Then there exists a unique sequence of positive integers $i_x=\bigl(i_x(k)\bigr)_{k=0}^\infty$ such that \[x\in \bigcap_{k=0}^\infty I_{k,i_x(k)}.\]
        Next, take $r\in(0,\infty)$. Notice that, by assumption \eqref{eq:seqgrow}, in particular $a_m>2a_{m-1}$ for all $m\geq 2$. Thus, we may write
        \[
        (0,\infty)=\Big(0,\frac{a_1}{2}\Big)\cup \bigcup_{m=1}^\infty \Big[\frac{a_m}{2},\frac{a_{m+1}}{2}\Big).
        \]
        This leads us to examine two possible cases.
        \begin{itemize}
            \item \textbf{Case 1: there exists \boldmath$m_r\in\{1,2,\dots\}$ such that $r\in\Big[\frac{a_{m_r}}{2}, \frac{a_{m_r+1}}{2}\Big)$.}
            
In  this case, by \eqref{eq:seqgrow}, we have that
\[
r\geq \frac{a_{m_r}}{2}>\operatorname{diam}\bigl(I_{m_r-1,i_x(m_r-1)}\bigr).
\]
On the other hand, by definition of the sequence $i_x$, we have that $x\in I_{m_r-1,i_x(m_r-1)}$. Thus,
            \[
            I_{m_r-1,i_x(m_r-1)}\subseteq (x-r,x+r)\cap X_{a}=:B_{X_{a}}(x,r),
            \]
            which implies, by \eqref{eq:Cantori},
            \begin{equation}\label{eq:proofCantor1}
                m_1\bigl(B_{X_{a}}(x,r)\bigr)\geq m_1\bigl(I_{m_r-1,i_x(m_r-1)}\bigr)=2^{m_r-1}.
            \end{equation}
            Next, by \eqref{eq:Cantorii1}, 
            \[
            2r\leq a_{m_r+1}\leq \begin{cases}
                \dist \bigl(I_{m_r,i_x(m_r)},I_{m_r,i_x(m_r)+1}\bigr),\\
                \dist \bigl(I_{m_r,i_x(m_r)},I_{m_r,i_x(m_r)-1}\bigr), \quad \text{if }i_x(m_r)>1,
            \end{cases}
            \]
            so that, keeping in mind that $x\in I_{m_r,i_x(m_r)}$, we get
            \[
            B_{X_{a}}(x,2r)\subseteq I_{m_r,i_x(m_r)}.
            \]
            Therefore, making use of \eqref{eq:proofCantor1}, we get that
            \begin{equation*}
            m_1\bigl(B_{X_{a}}(x,2r)\bigr)\leq m_1\bigl(I_{m_r,i_x(m_r)}\bigr)= 2^{m_r}\leq 2m_1\bigl(B_{X_{a}}(x,r)\bigr).
            \end{equation*}

\item \textbf{Case 2: \boldmath$r\in\Big(0,\frac{a_1}{2}\Big)$.}\\
In this case, the two balls $B_{X_{a}}(x,r)$ and $B_{X_{a}}(x,2r)$ can only contain points belonging to the length-$1$ interval to which $x$ belongs, namely $I_{0,i_x(0)}$. Indeed, by \eqref{eq:Cantorii1},
\[
r<2r<a_1\leq \begin{cases}
                \dist \bigl(I_{0,i_x(0)},I_{0,i_x(0)+1}\bigr),\\
                \dist \bigl(I_{0,i_x(0)},I_{0,i_x(0)-1}\bigr), \quad \text{if }i_x(0)>1.
            \end{cases}
\]
Therefore, in this case, everything reduces to verifying that a length-$1$ interval in the real line, equipped with the $1$-dimensional Lebesgue measure, is doubling. This verification is straightforward (the constant obtained is, again, $2$) and is left to the reader.
        \end{itemize}
In conclusion, we obtain that the space under consideration is doubling with constant $C=2$. 

        \item \textbf{Point~\ref{Cantorpoint:iv}.}\\
        Let $a_n:=2^{2^n}$, for $n\geq1$, and let us show that such choice of the sequence $a$ satisfies condition \eqref{eq:seqgrow}, by induction.
        \begin{itemize}
            \item \textbf{Base case.}\\
For $n=1$ we have
\[
a_1=4>2=2(2^{1-1}+0).
\]
            \item \textbf{Inductive step.}\\
            Assume that \eqref{eq:seqgrow} is satisfied for $n$. Then,
            \begin{equation*}
                \begin{split}
                    2 \Big(2^{n}+\sum_{k=1}^{n}2^{n-k} a_{k}\Big)&=2\Big(\underbrace{2\cdot 2^{n-1}+2\sum_{k=1}^{n-1}2^{n-1-k} a_{k}}_{<a_n}+a_n\Big)\\&<2(a_n+a_n)=4a_n=2^{2+2^n}\leq 2^{2\cdot2^n}=2^{2^{n+1}}=a_{n+1}
                \end{split}
            \end{equation*}
\qedhere

            
        \end{itemize}

\end{itemize}
\end{proof}

\subsection{Failure of a positive lower dimension}
Here, we show that the Cantor-type space $X_a$ does not satisfy the lower dimension $d$ property (see Definition \ref{def:DimSpace}), for any choice of $d>0$. In what follows, we denote
\[
V_{X_a}(x,r):=m_1\bigl(B_{X_a}(x,r)\bigr)=m_1\Big((x-r,x+r)\cap X_a\Big), \quad x\in X_a, \, r\in(0,\infty),
\]
and 
\[
V_{X_a}(x,y):=V_{X_a}(x,|x-y|), \quad x,y\in X_a.
\]
\begin{lemma}\label{lem:CantorNotLow}
     Let $a_n:=2^{2^n}$, for $n\geq1$, and consider the space $X_{a}$ that we have constructed in this appendix. Then there do not exist constants $c\in(0,\infty)$ and $d\in(0,\infty)$, such that
    \[
    \frac{V_{X_{a}}(x,R)}{V_{X_{a}}(x,r)}\geq c\Big(\frac{R}{r}\Big)^d,
    \]
    for all $x\in X_{a}$ and $0<r\leq R<\infty$. In other words, the space $(X_{a},|\cdot-\cdot|,m_1\vert_{X_{a}})$ does not satisfy the lower dimension $d$ property, for any choice of $d\in(0,\infty)$.
\end{lemma}
\begin{proof}
    Assume for contradiction that there exist $c\in(0,\infty)$ and $d\in(0,\infty)$ such that
    \begin{equation}\label{eq:proofCantor2}
        \frac{V_{X_{a}}(x,R)}{V_{X_{a}}(x,r)}\geq c\Big(\frac{R}{r}\Big)^d,
    \end{equation}
    uniformly in $x\in X_{a}$ and $0<r\leq R<\infty$. Next, for $n\geq 1$, choose a point $x\in I_{n-1,1}\subset I_{n,1}$, a smaller radius $r=a_n=2^{2^n}$, and a bigger radius $R=a_{n+1}=2^{2^{n+1}}$. Then, by \eqref{eq:Cantori}, \eqref{eq:Cantorii1}, and since $a$ satisfies \eqref{eq:seqgrow}, we get that
    \begin{equation*}
        V_{X_{a}}(x,R)=m_1(I_{n,1})=2^n, \quad V_{X_{a}}(x,r)= m_1(I_{n-1,1})=2^{n-1}.
    \end{equation*}
    Thus, by \eqref{eq:proofCantor2}, we get that, for any $n\geq1$, there holds
    \[
    \frac{2^n}{2^{n-1}}\geq c\Big(\frac{2^{2^{n+1}}}{2^{2^n}}\Big)^d\quad \Leftrightarrow \quad 2^{d2^n}\leq \frac2c,
    \]
    which is clearly a contradiction, since $\lim_{n\to\infty}2^{d2^n}=\infty$.
\end{proof}

\subsection{Metric vs. volumic fractional integrals}
As in the case of the infinite slab constructed in Appendix~\ref{app:Slab}, also in this setting it is not possible to view the metric fractional integrals as a special case of the volumic ones. This is a consequence of the following result.

\begin{lemma}
Let $a_n := 2^{2^n}$ for $n\ge 1$, and let $X_a$ be the space constructed in this appendix. 
Then there do not exist constants $c>0$ and $\alpha>0$ such that
\[
    |x-y|^\alpha \le c\, V_{X_a}(x,y)
\]
for all $x,y\in X_a$.
\end{lemma}

\begin{proof}
    Assume for contradiction that there exist $c\in(0,\infty)$ and $\alpha\in(0,\infty)$, such that
    \[
    |x-y|^\alpha\leq cV_{X_{a}}(x,y),
    \]
    uniformly in $x,y\in X_{a}$. In particular, for $n\geq1$, let us take $x\in I_{n-1,1}$ and $y\in I_{n-1,2}$. By \eqref{eq:Cantorii1} and since $a$ satisfies \eqref{eq:seqgrow}, we have that
    \begin{equation*}
        2^{\alpha2^n}=(a_n)^\alpha\leq |x-y|^\alpha\leq cV_{X_{a}}(x,y)\leq cm_1({I_{n,1}})=c2^n.
    \end{equation*}
    Therefore, we get that, for any $n\geq 1$, 
    \[
    2^{\alpha2^n-n}\leq c,
    \]
    which is clearly a contradiction, since $\lim_{n\to\infty}2^{\alpha2^n-n}=\infty$.
\end{proof}

\subsection{Failure of any \texorpdfstring{$L^p$}{Lp} theory for metric fractional integrals}\label{subsec:CantorUnbdd}

In what follows, we show that the structure of this Cantor-type set, associated with a suitable increasing sequence, makes it completely impossible to obtain any $L^p$ theory for the metric fractional integrals.

\begin{lemma}\label{lem:unbddfracmet}
    Let $X_a$ be Cantor-type set we constructed, associated with the sequence $a=(a_n)_{n=1}^\infty$ defined as $a_n:=2^{2^n}$. Let us equip such set with the $1$-dimensional Lebesgue measure $\mu := m_1\vert_{X_a}$ restricted to $X_a$, and the distance $\rho(x,y) := |x-y|$. Next, let $\eta\in(0,\infty)$ and consider the operator
    \[
    Sg(x):=\int_{X_a} \frac{|x-y|^\eta}{V_{X_a}(x,y)}g(y)\ud m_1(y), \quad x\in X_a.
    \]
    Then, for any $p\in (0,\infty)$, there exists a function $f_p\in L^p(X_a)$ such that $Sf_p(x)=\infty$ for all $x\in X_a$.
\end{lemma}

\begin{proof}
    Let $p\in (0,\infty)$ and let 
    \[
    f_p(x):=\sum_{k=1}^\infty \frac{1}{2^{k/p}} \mathbf{1}_{I_{0,2^k+1}}.
    \]
    Essentially, this function takes the value $1/2^{k/p}$ on the first length-$1$ interval belonging to the block $I_{k,2}$, for every $k \geq 1$. Therefore, we have that
\[
\Norm{f_p}{L^p(X_a)}=\left(\int_{X_a}\abs{f_p}^p\ud m_1\right)^{1/p}=\left(\sum_{k=1}^\infty \frac{1}{2^k}\right)^{1/p}=1<\infty.
\]

Now, let $x\in X_a=\bigcup_{k\geq 0} I_{0,k}$. Then, there exists $k_x\geq0$ such that $x\in I_{0,k_x}$. Next, let us take $k_0\geq0$ such that $k_x \leq 2^{k_0}$, and write
\begin{equation*}
    \begin{split}
        Sf_p(x)&=\int_{X_a} \frac{|x-y|^\eta}{V_{X_a}(x,y)}f_p(y)\ud m_1(y)\\
        &=\sum_{k=1}^\infty \int_{I_{0,2^k+1}}\frac{|x-y|^\eta}{V_{X_a}(x,y)}\frac{1}{2^{k/p}}\ud m_1(y)\\
        &\geq \sum_{k=k_0}^\infty \int_{I_{0,2^k+1}}\frac{|x-y|^\eta}{V_{X_a}(x,y)}\frac{1}{2^{k/p}}\ud m_1(y)=:J.
    \end{split}
\end{equation*}
By the choice of $k_0$, we have $x \in I_{k_0,1}$. Moreover, for every $k \geq k_0$, it follows that $x \in I_{k,1}$ and $x \notin I_{k,2}$. On the one hand, as already mentioned before, for every $k \geq k_0$, the interval $I_{0,2^k+1}$ is the first length-$1$ interval contained in $I_{k,2}$. Therefore, by \eqref{eq:Cantorii1}, for $y\in I_{0,2^k+1}$ we have that
\begin{equation*}
    \abs{x-y}\geq a_{k+1}=2^{2^{k+1}}.
\end{equation*}
On the other hand, such $x$ and $y$ both belong to $I_{k+1,1}$. Thus, by \eqref{eq:Cantori}, \eqref{eq:Cantorii1} and \eqref{eq:seqgrow}, we have that
\begin{equation*}
    V_{X_a}(x,y)\leq m_1\bigl(I_{k+1,1}\bigr)=2^{k+1}.
\end{equation*}
Altogether, we get that
\[
J\geq \sum_{k=1}^\infty \frac{2^{\eta 2^{k+1}}}{2^{(1+1/p)k+1}}=\infty.
\]
This concludes the proof.
\end{proof}

This result is noteworthy, as for every $p \in (0,\infty)$ we construct a function $f_p \in L^p$ such that its image under the operator $S$ is identically $+\infty$. Consequently, not only can one not hope to obtain an $L^p$--$L^q$ boundedness as in the volumic case, but one cannot even expect a weaker boundedness, such as that obtained in Theorem~\ref{thm:bddfracmet} for metric fractional integrals on spaces with positive lower dimension. In light of this example, the latter property turns out to be indispensable, and essentially minimal, in order to ensure such boundedness.

\section*{Acknowledgements}
This paper is part of the author's forthcoming doctoral thesis written under the supervision of Professor Tuomas Hytönen. The author would like to sincerely thank him for his valuable advice, for the many helpful conversations, and for suggesting the study of these topics, which the author found extremely interesting. Indeed, without his suggestion, the author would likely have never become aware of them nor studied them. The author is supported by the Research Council of Finland through the Finnish Centre of Excellence in Randomness and Structures ``FiRST'' (grant no. 364208).

\section*{Use of Artificial Intelligence}

In the preparation of this manuscript, the author made use of the free versions of three publicly available AI assistants: ChatGPT (OpenAI), Gemini (Google), and Claude (Anthropic). No paid subscriptions or premium models were used. Given the length of time over which this manuscript was developed, no specific model versions are listed.

Artificial intelligence was used to assist the author with the LaTeX code, particularly in generating the TikZ code for the figures included herein. It was also used to translate certain passages from Italian into English. Although the author is proficient in written and spoken English, it was sometimes faster and more effective to formulate certain ideas in his native language, allowing him to express the intended meaning and nuances more precisely.

Furthermore, artificial intelligence was frequently used to assist with bibliographic searches for specific results, as this provided an efficient way to determine whether certain results were already available in the literature and to identify the relevant references for proper citation.

As far as the mathematical content is concerned, however, the contribution of artificial intelligence was negligible. All the mathematical results presented in this manuscript were obtained by the author through independent mathematical work, including the development of the arguments and proofs, under the supervision of Professor Tuomas Hytönen.

\printaddress

\end{document}